\documentclass[10pt]{articleHJ}

\usepackage[centertags]{amsmath}
\usepackage{amsfonts,amssymb,amsthm}
\usepackage{caption,graphicx}
\usepackage{subcaption}
\usepackage[text={6.0in,8.6in},centering,letterpaper]{geometry}
\usepackage[numbers,comma,square,sort&compress]{natbib}

\usepackage{tabularx}

\usepackage{mathtools}

\newlength{\defbaselineskip}
\newcommand{\setlinespacing}[1]%
           {\setlength{\baselineskip}{#1 \defbaselineskip}}

\newcommand{\singlespacing}{\setlength{\baselineskip}{\defbaselineskip}}

\newcommand{\R}{\ensuremath{\mathbb{R}}}

\renewcommand{\b}{\ensuremath{\beta}}
\newcommand{\e}{\ensuremath{\varepsilon}}

\renewcommand{\d}{\ensuremath{\delta}}
\newcommand{\g}{\ensuremath{\gamma}}
\newcommand{\G}{\ensuremath{\Gamma}}
\newcommand{\s}{\ensuremath{\sigma}}

\renewcommand{\a}{\ensuremath{\alpha}}

\renewcommand{\L}{\ensuremath{\mathcal{L}}}
\renewcommand{\O}{\ensuremath{\mathcal{O}}}

\DeclarePairedDelimiter{\nrm}\lVert\rVert

\def\Re{\mathop\mathrm{Re}\nolimits}			
\newcommand{\Real}{\mathbb{R}}							
\newcommand{\abs}[1]{\left\vert#1\right\vert}			
\newcommand{\norm}[1]{\left\Vert#1\right\Vert}		
\newcommand{\sref}[1]{(\ref{#1})}                       

\newtheorem{thm}{Theorem}[section]
\newtheorem{cor}[thm]{Corollary}
\newtheorem{lem}[thm]{Lemma}
\newtheorem{prop}[thm]{Proposition}

\newtheorem*{defn*}{Definition}

 {\begin{trivlist}\item[]\textbf{Acknowledgments }}{\end{trivlist}}

\numberwithin{equation}{section}

\begin{document}

\begin{frontmatter}
\title{Stability of Traveling Waves \\ for Reaction-Diffusion Equations \\
  with Multiplicative Noise}
\journal{SIADS}

\author[LD1]{C. H. S. Hamster\corauthref{coraut}},
\corauth[coraut]{Corresponding author. }
\author[LD2]{H. J. Hupkes},
\address[LD1]{
  Mathematisch Instituut - Universiteit Leiden \\
  P.O. Box 9512; 2300 RA Leiden; The Netherlands \\
  Email:  {\normalfont{\texttt{c.h.s.hamster@math.leidenuniv.nl}}}
}
\address[LD2]{
  Mathematisch Instituut - Universiteit Leiden \\
  P.O. Box 9512; 2300 RA Leiden; The Netherlands \\ Email:  {\normalfont{\texttt{hhupkes@math.leidenuniv.nl}}}
}

\date{Version of \today}

\begin{abstract}
\singlespacing
We consider reaction-diffusion equations
that are stochastically forced by a small multiplicative noise term.
We show that spectrally stable traveling wave solutions
to the deterministic system retain their orbital stability
if the amplitude of the noise is sufficiently small.

By applying a stochastic phase-shift
together with a time-transform, we obtain a
semilinear SPDE that describes the fluctuations from the
primary wave. We subsequently develop a semigroup approach to handle the
nonlinear stability
question in a fashion that is closely related to modern
deterministic methods.

\end{abstract}



\end{frontmatter}

\section{Introduction}
\label{sec:int}

In this paper we consider stochastically perturbed versions
of a class of reaction-diffusion equations
that includes the bistable Nagumo equation
\begin{equation}
\label{eq:int:nagumo}
u_t = u_{xx} + f_{\mathrm{cub}}(u)
\end{equation}
and the FitzHugh-Nagumo equation
\begin{equation}
\label{eq:int:fzhnag:pde}
\begin{array}{lcl}
u_t & = & u_{xx} + f_{\mathrm{cub}}(u) - v
\\[0.2cm]
v_t & = & v_{xx} + \varrho [ u - \gamma v ]
\end{array}
\end{equation}
Here we take $\varrho > 0$, $\gamma > 0$ and consider
the standard bistable nonlinearity
\begin{equation}
f_{\mathrm{cub}}(u) = u ( 1 - u )( u -a ).
\end{equation}

It is well-known \cite{Fife1977,Sattinger} that \sref{eq:int:nagumo}
admits spectrally stable traveling front solutions
\begin{equation}
u(x,t) = \frac{1}{2}
  \big[ 1 + \tanh(\frac{1}{4} \sqrt{2} (x - ct) \big) \big]
\end{equation}
that travel with speed
\begin{equation}
c =  \sqrt{2}  \left(a - \frac{1}{2} \right).
\end{equation}
In addition, the existence of traveling
pulse solutions to \sref{eq:int:fzhnag:pde}
with $0 < \varrho \ll 1$
was established recently \cite{chen2015traveling}
using variational methods. Using the Maslov index,
a proof for the spectral stability of these waves
has recently been obtained in \cite{cornwell2017opening,cornwell2017existence}.

Our main results show that these spectrally
stable wave solutions survive in a suitable
sense upon adding a small pointwise
multiplicative
noise term to the underlying PDE. This
noise term is assumed to be globally Lipschitz and
must vanish for the asymptotic
values of the waves. For example,
our results
cover the scalar Stochastic Partial Differential Equation (SPDE)
\begin{equation}
\label{eq:int:SPDE:nag}
dU = \big[ U_{xx} + f_{\mathrm{cub}}(U) ] dt
+ \sigma \chi(U)  U ( 1 - U) d \beta_t
\end{equation}
together with the two-component SPDE
\begin{equation}
\label{eq:int:SPDE:fhn}
\begin{array}{lcl}
dU & = & \big[ U_{xx} + f_{\mathrm{cub}}(U) - V] dt
+ \sigma \chi(U)  U ( 1 - U) d \beta_t ,
\\[0.2cm]
dV & = & \big[V_{xx} + \varrho \big( U - \gamma V \big) \big] dt + \sigma (U - \gamma V) d \beta_t ,
\end{array}
\end{equation}
both for small $\abs{\sigma}$,
in which $(\beta_t)$
is a Brownian motion and $\chi(U)$
is a cut-off function with $\chi(U) = 1$
for $\abs{U} \le 2$. The presence
of this cut-off is required to
enforce the global Lipschitz-smoothness of the noise term.
In this regime, one can think of \sref{eq:int:SPDE:nag} and
\sref{eq:int:SPDE:fhn}
as versions of the PDEs \sref{eq:int:nagumo}-\sref{eq:int:fzhnag:pde}
where the parameters $a$ and $\varrho$ are replaced by
$a + \sigma  \dot{\beta}_t$ respectively $\varrho + \sigma  \dot{\beta}_t$.

Many additional
multi-component reaction-diffusion PDEs
such as the Gray-Scott
\cite{lee1994experimental},
Rinzel-Keller \cite{rinzel1973traveling},
Tonnelier-Gerstner
\cite{tonnelier2003piecewise}
and Lotka-Volterra systems \cite{iannelli2015introduction}
are also known to admit spectrally stable
traveling waves in the equal-diffusion
setting
\cite{hale1999stability, zemskov2010wave, gardner1982existence}.
This allows our results to be applied to these waves
after appropriately truncating the
deterministic nonlinearities
(in regimes that are far away from the interesting
dynamics).

Such cut-offs are not necessary
when considering equal-diffusion
three-component FitzHugh-Nagumo-type systems
such as those studied in \cite{ORG1998,van2010front}.
Such equations were first used by Purwins
to study the formation
of patterns during gas discharges \cite{SCHENK1997}.
However, in the equal-diffusion setting there is at present
only numerical evidence to suggest that spectrally stable waves
exist for the underlying deterministic equation.
Analytical approaches to prove such facts typically use methods from
singular perturbation theory, but these often require
the diffusive length scales to be strictly separated.

\paragraph{Noisy patterns}
Stochastic forcing of PDEs has become an important tool for modellers in a
large number of
fields, ranging from medical applications such as neuroscience
\cite{Bressloff,bressloff2015nonlinear}
and cardiology \cite{Zhang} to
finance \cite{NunnoAdvMathFinance2011}
and meteorology \cite{Climate}.
While a rather general existence theory for
solutions to SPDEs has been developed over the past decades
\cite{Chow,Concise,Gawarecki,DaPratoZab},
the study of patterns such as stripes, spots and waves in
such systems is less well-developed.

Preliminary results for specific
equations such as Ginzburg-Landau
\cite{brassesco1995brownian,funaki1995scaling}
and Swift-Hohenberg \cite{Kuske2017} are available.
K\"uhn and Gowda \cite{gowda2015early} analyzed
both these equations in the linear regime before the onset
of the Turing bifurcation. They obtained
scaling laws for the natural covariance operators
that can be used as early-warning signs to predict the appearance of patterns.

In addition, several numerical studies
have been initiated to study the impact of noise on patterns,
see e.g. \cite{vinals1991numerical,Lord2012,Shardlow}.
The results in \cite{Lord2012} relating to
\sref{eq:int:SPDE:nag} are particularly interesting
from our perspective.
Indeed, they clearly show that traveling wave solutions
persist under the stochastic forcing, but the speed decreases
linearly in $\s^2$ and the wave becomes steeper.

Rigorous results concerning the impact of stochastic forcing on deterministic waves
are still relatively scarce. However, some important contributions have already been made,
focusing on two important issues that need to be addressed.
The first of these is that one needs to identify appropriate mechanisms
to identify the phase, speed and shape of a stochastic wave.
The second issue is that one needs to control the influence of the nonlinear
terms by using the decay properties of the linear terms.

\paragraph{Phase tracking}

An appealing intuitive idea is to define the phase $\vartheta(u)$ of a solution profile $u$
relative to the deterministic traveling wave $\Phi$
by writing
\begin{equation}
\label{eq:int:def:phase}
\vartheta(u) = \mathrm{argmin}_{\vartheta \in \Real} \norm{ u - \Phi(\cdot + \vartheta) }_{L^2},
\end{equation}
which picks the closest translate of $\Phi$.
Inspired by this idea,
Stannat \cite{Stannat,stannat2014stability} 
obtained
orbital stability results for a class of
systems including \sref{eq:int:SPDE:nag} by appending
an ODE to track the position of the wave. This is done via a gradient-descent
technique, whereby the phase is updated continuously in the direction
that lowers the norm in \sref{eq:int:def:phase}. A slight drawback of this method
is that the phase is always lagging in a certain sense. In particular, it is not immediately clear
how to define a stochastic speed and relate it with its deterministic counterpart.

This gradient-descent approach has been extended
to neural field equations with additive noise \cite{Lang,LangStannat2016l2}.
In order to clarify the dynamic effects caused by the noise, the authors
employed a perturbative approach and expanded the phase of the wave and the shape of
 the perturbations
in powers of the noise strength $\sigma$. By taking the infinite update-speed limit, the authors were able
to eliminate the phase lag mentioned above.  At lowest order they roughly recovered the diffusive
wandering of the phase that was predicted
by Bressloff and Webber \cite{Bressloff}. This perturbative expansion
can be maintained on finite time intervals, which increase in length to infinity
as the noise size $\sigma$ is decreased. However, one needs
separate control on the deviations
of the  phase and the shape from the deterministic wave,
which are both required to stay small.

Inglis and MacLaurin take a directer approach in \cite{Inglis}
by using a stochastic differential equation for the phase that forces \sref{eq:int:def:phase} to hold.
For equations with additive noise, they obtain results that allow waves
to be tracked over finite time intervals. As above, this tracking time increases
to infinity as $\sigma \downarrow 0$.
The main issue
here is that global minima do not necessarily behave in a continuous fashion.
This means that \sref{eq:int:def:phase} can become multi-valued at times, leading
to sudden jumps of the phase.
However, under a (restrictive) technical condition the extension of the tracking time
can be performed uniformly in $\sigma$.

\paragraph{Nonlinear effects}
In order to control the nonlinear terms over long time intervals
one needs the linear flow to admit suitable decay properties.
We write $S(t)$ for the semigroup generated by the
linear operator $\mathcal{L}_{\mathrm{tw}}$ associated to the linearization
of the PDEs above around their traveling wave $\Phi$.
A direct consequence of the translational invariance is that
$\mathcal{L}_{\mathrm{tw}} \Phi' = 0$ and hence $S(t) \Phi' = \Phi'$ for all $t > 0$.
In order to isolate this neutral mode, we write $P$ for the spectral projection onto $\Phi'$,
together with its complement $Q = I-P$.
Assuming a standard spectral gap condition on the remainder of the spectrum of $\mathcal{L}_{\mathrm{tw}}$,
one can subsequently obtain the estimate
\begin{equation}
\label{eq:int:bnd:semigroup}
\norm{S(t) Q}_{L^2 \to L^2} \le M e^{ - \beta t}
\end{equation}
for some constants $\beta > 0$ and $M \ge 1$;
see, for example, \cite[Lem. 5.1.2]{volpert1994traveling}.

The common feature in all the approaches described above
is that they require the identity $M = 1$ to hold. In this special
case the linear flow is immediately contractive in the direction orthogonal to the
translational eigenfunction.
This identity certainly holds if one can obtain an estimate of the form
\begin{equation}
\label{eq:int:coercivity:eq}
\langle \mathcal{L}_{\mathrm{tw}} v , v \rangle \le -  \beta \norm{v}_{H^1}^2  + \kappa \norm{P v}^2
\end{equation}
for some $\kappa > 0$, since one can then use the commutation property
$P S(t)  =  S(t) P$ to compute
\begin{equation}
\frac{d}{dt} \norm{S(t) Q v}_{L^2}^2 = 2 \langle \mathcal{L}_{\mathrm{tw}} S(t) Qv, S(t) Q v \rangle_{L^2}
\le - 2 \beta \norm{S(t) Q v}_{H^1}^2
\le - 2 \beta \norm{S(t) Q v}_{L^2}^2 .
\end{equation}
In the deterministic case,  coercive estimates of this type
can be used to obtain similar differential inequalities
for the $L^2$-norm of perturbations from the phase-adjusted traveling wave.
Using the It{\^o} formula this can be generalized to the stochastic
case \cite{Stannat,stannat2014stability,Lang,LangStannat2016l2},
allowing stability estimates to be obtained that do not need
any control over the $H^1$-norm of these perturbations.
The approach developed in \cite{Inglis} proceeds directly
from \sref{eq:int:bnd:semigroup} using a renormalisation method.
Similar $L^2$-stability results can be obtained in this fashion,
again crucially using the fact that $M = 1$; see \cite[(6.15)]{Inglis}.

In light of the discussion above, a considerable effort is underway
to identify systems for which the immediate contractivity condition $M=1$ indeed holds.
This has been explicitly verified
for the Nagumo PDE
\sref{eq:int:SPDE:nag} and several classes of one-component systems
\cite{Stannat,stannat2014stability, LangStannat2016l2}. However,
these computations are very delicate and typically proceed on an ad-hoc basis.
For example, it is unclear (and doubtful) whether such a condition holds
for the FitzHugh-Nagumo PDE \sref{eq:int:SPDE:fhn}. We refer
to \cite[{\S}1]{veraar2011note}
for an informative discussion on this issue.

In the case $M > 1$ the semigroup is still eventually contractive on the range of $Q$,
but it can cause transient dynamics that grow on short timescales.
Such dynamics play an important role and need to be tracked over temporal intervals
of intermediate length. In this case the nonlinearities cannot be immediately
dominated by the linear terms
as above. To control these terms it is hence crucial to understand
the $H^1$-norm of perturbations, which poses some challenging regularity issues
in the stochastic setting.

\paragraph{Semigroup approach}

In this paper we take a step towards harnessing the power
of modern deterministic nonlinear stability techniques
for use in the stochastic setting.
In particular, inspired by the informative
expository paper \cite{Zumbrun2009}, we abandon
any attempt to describe the phase of the wave via
a priori geometric conditions. Instead,
we initiate a semigroup approach based on the
stochastic variation of constants formula. This leads to
a stochastic evolution equation for the phase
that follows naturally from technical considerations.
More specifically, we use the phase to neutralize
the dangerous non-decaying terms in our evolution equation.
Our tracking mechanism is robust and allows us to focus solely
on the behaviour of the perturbation from the phase-shifted wave.
This allows us to track solutions up to the point where this
perturbation becomes too large as a result of the stochastic forcing,
which resembles an Ornstein-Uhlenbeck process and hence is unbounded
almost certainly.
In particular, we do not need to impose restrictions
on the size of the phaseshift
as in \cite{Lang,kruger2017multiscale}.

The first main advantage of our approach is that it provides
orbital stability results without requiring the
immediate contractivity condition described above.
Indeed, we are able to track the $H^1$-norm of perturbations
and not merely the $L^2$-norm, which allows us to have $M > 1$
in \sref{eq:int:bnd:semigroup}.
This significantly broadens
the class of systems that can be understood and
aligns the relevant spectral assumptions
with those that are traditionally used in deterministic settings.

The second main advantage is that we are (in some sense) able to isolate
the drift-like contributions to the shape and speed of the wave
that are caused by the noise term. This becomes fully visible
in our analysis of \sref{eq:int:SPDE:nag}, where the noise term is
specially tailored to the deterministic wave $\Phi$ in the sense that it is proportional to
the neutral mode $\Phi'$. In this case we are able
to obtain an exponential stability result for a modified waveprofile $\Phi_{\sigma}$
that propagates with a modified speed $c_{\sigma}$
and exists for all positive time.
This allows us to rigorously
understand the changes to the waveprofile and speed that
were numerically observed for \sref{eq:int:SPDE:nag} in \cite{Lord2012}.
In general, if the $\Real^n$-orbit of the traveling wave
of an $n$-component reaction-diffusion equation contains no self-intersections,
our results allow special forcing terms to be constructed
for which the modified waves remain exponentially stable.

However, the need to use stochastic calculus
causes several delicate technical complications that are not observed in the
deterministic setting. For example, the It\^o Isometry is based on $L^2$ norms.
At times, this forces us to square the natural semigroup decay rates,
which leads to short-term regularity issues.
Indeed, the heat semigroup $S(t)$ behaves as $\norm{S(t)}_{\mathcal{L}(L^2; H^1)} \sim t^{-1/2}$,
which is in $L^1(0,1 )$ but not in $L^2(0,1)$.
This precludes us from obtaining supremum control on the $H^1$-norm of our solutions.
Instead, we obtain bounds on square integrals of the $H^1$-norm. For this reason,
we need to carefully track how the cubic behaviour of $f_{\mathrm{cub}}(u)$ propagates through
our arguments.

A second major complication is that stochastic phase-shifts lead to extra nonlinear
diffusive terms. By contrast, deterministic phase-shifts lead to extra convective terms,
which are of lower order and hence less dangerous.
As a consequence, we encounter quasi-linear equations
in our analysis that do not immediately fit into a semigroup framework.
We solve this problem by using a suitable stochastic time-transform to scale
out the extra diffusive terms. The fact that we need
the diffusion coefficients in \sref{eq:int:fzhnag:pde} to be identical
is a direct consequence of this procedure.

\paragraph{Outlook}

Let us emphasize that we view the present paper merely as a
proof-of-concept result for a pure semigroup-based approach.
For example, in a companion paper \cite{Hamster2018Uneq} we show how the severe restriction
on the diffusion coefficients of \sref{eq:int:fzhnag:pde} can be removed
by exploiting the block-structure
of the semigroup.

In addition, our results here use the variational framework
developed by Liu and R\"ockner \cite{LiuRockner} in order
to ensure that our SPDE has a well-defined global weak solution. In future work,
we intend to replace this procedure by constructing local mild solutions
directly based on fixed-point arguments.

Finally, we are interested in more delicate spectral stability scenarios,
which allow one or more branches of essential spectrum to touch the imaginary axis
in a quadratic tangency. Situations of this type are encountered
when analyzing the two-dimensional stability of traveling planar
waves \cite{BHM, KAP1997, HJHSTB2D, HJHOBST2D}
or when studying viscous shocks in the context of conservation laws
\cite{beck2010nonlinear,MasciaZumbrun02,HJHNLS}.

\paragraph{Organization}

This paper is organized as follows.
We formulate our phase-tracking mechanism and state our main results in \S\ref{sec:mr}.
In \S\ref{sec:prlm} we obtain preliminary estimates on our nonlinearities,
which are used in \S\ref{sec:var} to fit our coupled SPDE into the theory
outlined in \cite{Concise, LiuRockner}. This guarantees that our SPDE
has well-defined solutions, to which we apply a stochastic phase-shift
in \S\ref{sec:sps} followed by a stochastic time-transform in \S\ref{sec:md}.
These steps lead to a stochastic variation of constants formula.

In \S\ref{sec:swv} we develop two fixed-point arguments that capture
the modifications to the waveprofile and speed that arise from the stochastic forcing.
These modifications allow us to obtain suitable estimates on the nonlinearities
in the variation of constants formula in \S\ref{sec:fnl},
which allow us to pursue a nonlinear-stability argument in \S\ref{sec:nls}.

\paragraph{Acknowledgements.}
Hupkes acknowledges support from the Netherlands Organization
for Scientific Research (NWO) (grant 639.032.612).
Both authors wish to thank O. van Gaans and C. da Costa
for helpful discussions during the conception and writing of this paper. In addition, some of our results were inspired by the valuable comments made by two anonymous referees.

\section{Main results}
\label{sec:mr}

In this paper we are interested in the stability of traveling
wave solutions to SPDEs of the form
\begin{equation}
\label{eq:mr:main:spde}
dU =
  \big[ A_* U + f(U) \big] dt
  + \sigma g(U) d \beta_t.
\end{equation}
Here we take $U = U(x,t) \in \Real^n$ with $x \in \Real$ and $t \ge 0$.

In \S\ref{sec:mr:det} we formulate several
conditions on the nonlinearity $f$ and the diffusion operator $A_*$,
which imply that in the deterministic case $\sigma =0 $
the system \sref{eq:mr:main:spde}
has a variational structure
and admits a spectrally stable traveling wave solution.
In \S\ref{sec:mr:stoch} we impose several standard
conditions on the noise term in \sref{eq:mr:main:spde},
which guarantee that \sref{eq:mr:main:spde}
is covered by the variational framework developed in \cite{LiuRockner}.
In addition, we couple an extra SDE to our SPDE that will serve as a phase-tracking
mechanism.
Finally, in \S\ref{sec:mr:wave} and \S\ref{sec:mr:disc} we formulate
and discuss our main results
concerning the impact of the noise term on the deterministic traveling wave solutions.

\subsection{Deterministic setup}
\label{sec:mr:det}

We start here by stating our conditions on the form of $A_*$
and $f$. These conditions require $A_*$ to be a diffusion operator
with identical diffusion coefficients and
restrict the growth-rate of $f$ to be at most cubic.

\begin{itemize}
\item[(HA)]{
  For any $u \in C^2(\Real; \Real^n)$ we have
  $A_* u = \rho I_n u_{xx}$,
  in which $\rho >0$ and $I_n$ is the $n \times n$-identity matrix.
}
\item[(Hf)]{
 We have $f \in C^3(\Real^n; \Real^n)$
 and there exist $u_\pm \in \Real^n$
 for which $f(u_-) = f(u_+) = 0 $.
 In addition, there exists a constant $K_f > 0$
 so that the bound
 \begin{equation}
   \abs{ D^3 f(u) } \le K_f
 \end{equation}
 holds for all $u \in \Real^n$.
}
\end{itemize}

We now demand that the deterministic part
of \sref{eq:mr:main:spde}
has a traveling wave solution
that connects the two equilibria $u_\pm$
(which are allowed to be equal). This traveling
wave should approach these equilibria at an exponential rate.

\begin{itemize}
\item[(HTw)]{
  There exists a waveprofile $\Phi_{0} \in C^2(\Real ; \Real^n)$
  and a wavespeed $c_0 \in \Real$
  so that the function
  \begin{equation}
  \label{eq:mr:trv:wave:ansatz}
  u(x, t) = \Phi_{0}(x - c_0 t)
  \end{equation}
  satisfies the deterministic PDE
  \begin{equation}
    \label{eq:mr:det:pde}
    u_t = A_* u + f(u)
  \end{equation}
  for all $(x,t) \in \Real \times \Real$.
  In addition, there is a constant $K > 0$ together with
  exponents $\nu_\pm>0 $ so that the bound
  \begin{equation}
    \label{eq:mr:exp:apr:left}
    \abs{\Phi_0(\xi) - u_- } + \abs{\Phi_0'(\xi)} \le K e^{-\nu_- \abs{\xi} }
  \end{equation}
  holds for all $\xi \le 0$, while the bound
  \begin{equation}
    \label{eq:mr:exp:apr:right}
    \abs{\Phi_0(\xi) - u_+ } + \abs{\Phi_0'(\xi)}
      \le K e^{-\nu_+ \abs{\xi} }
  \end{equation}
  holds for all $\xi \ge 0$.
}
\end{itemize}

Throughout this paper, we will use the shorthands
\begin{equation}
L^2 = L^2(\Real;\Real^n),
\qquad
H^1 = H^1(\Real ;\Real^n),
\qquad
H^2 = H^2(\Real; \Real^n).
\end{equation}
Linearizing the deterministic PDE
\sref{eq:mr:det:pde} around
the traveling wave $(\Phi_{0} , c_0)$,
we obtain the linear operator
\begin{equation}
\mathcal{L}_{\mathrm{tw}}: H^2 \to L^2
\end{equation}
that acts as
\begin{equation}
[\mathcal{L}_{\mathrm{tw}} v](\xi) =
c_0 v'(\xi) + [A_* v](\xi) + Df\big(\Phi_0(\xi) \big)^T v(\xi) .
\end{equation}
The formal adjoint
\begin{equation}
\mathcal{L}_{\mathrm{tw}}^{\mathrm{adj}}: H^2 \to L^2
\end{equation}
of this operator acts as
\begin{equation}
[\mathcal{L}_{\mathrm{tw}}^{\mathrm{adj}} w](\xi) =
-c_0 w'(\xi) + [A_* w](\xi) + Df\big(\Phi_0(\xi) \big) w(\xi) .
\end{equation}
Indeed, one easily verifies that
\begin{equation}
\langle \mathcal{L}_{\mathrm{tw}} v , w \rangle_{L^2}
= \langle v, \mathcal{L}_{\mathrm{tw}}^{\mathrm{adj}} w \rangle_{L^2}
\end{equation}
whenever $(v,w) \in H^2 \times H^2$. Here $\langle \cdot, \cdot \rangle_{L^2}$
denotes the standard inner-product on $L^2$.

We now impose a standard spectral stability condition
on the wave. In particular, we require
that the standard translational eigenvalue
at zero is a simple eigenvalue. In addition,
the remainder of the spectrum of $\mathcal{L}_{\mathrm{tw}}$
must be strictly bounded to the left of the imaginary axis.
\begin{itemize}
\item[(HS)]{
  There exists $\beta > 0$ so that
  the operator $\mathcal{L}_{\mathrm{tw}} - \lambda$
  is invertible for all $\lambda \in \mathbb{C} \setminus \{0 \}$
  that have $\Re \lambda \ge - 2 \beta$,
  while $\mathcal{L}_{\mathrm{tw}}$ is a Fredholm operator
  with index zero. In addition, we have the identities
  \begin{equation}
  \mathrm{Ker}( \mathcal{L}_{\mathrm{tw}} ) = \mathrm{span} \{ \Phi_0' \},
  \qquad
  \mathrm{Ker}( \mathcal{L}_{\mathrm{tw}}^{\mathrm{adj}} )
    = \mathrm{span} \{ \psi_{\mathrm{tw}} \}
  \end{equation}
  for some $\psi_{\mathrm{tw}} \in H^2$ that has
  \begin{equation}
    \label{eq:mr:hs:norm:cnd:psitw}
     \langle \Phi_0' , \psi_{\mathrm{tw}} \rangle_{L^2} = 1.
  \end{equation}
}
\end{itemize}

We conclude by imposing
a standard monotonicity condition on $f$,
which ensures that the SPDE \sref{eq:mr:main:spde}
fits into the variational framework
of \cite{LiuRockner}. We remark here
that we view this condition purely as a technical
convenience, since it guarantees that solutions
to \sref{eq:mr:main:spde} do not blow up.
However, it does not a play a key role in the heart
of our computations, where we restrict our attention
to solutions that remain small in some sense.

\begin{itemize}
\item[(HVar)]{
 There exists $K_{\mathrm{var}} > 0$ so that the one-sided inequality
 \begin{equation}
    \langle f( u_A ) - f(u_B) , u_A - u_B \rangle_{\Real^n} \le K_{\mathrm{var}} \abs{ u_A - u_B}^2
 \end{equation}
 holds for all pairs $(u_A,u_B) \in \Real^n \times \Real^n$.
}
\end{itemize}

\subsection{Stochastic setup}
\label{sec:mr:stoch}

Our first condition here states
that the noise term in \sref{eq:mr:main:spde}
is driven by a standard Brownian motion.
Let us emphasize that we made this choice purely to enhance
the readability of our arguments.
Indeed, our results can easily be generalized
to the situation where the noise is driven by cylindrical $Q$-Wiener processes.
\begin{itemize}
\item[(H$\beta$)]{
  The process $(\beta_t)_{t \ge 0}$ is a Brownian motion
  with respect to the complete filtered probability space
  \begin{equation}
    \Big(\Omega, \mathcal{F}, ( \mathcal{F}_t)_{t \ge 0} , \mathbb{P} \Big).
  \end{equation}
}
\end{itemize}
We require the function $Dg$ to be globally Lipschitz
and uniformly bounded. While
the former condition is essential in our analysis
to ensure that our cut-offs only depend on $L^2$-norms,
the latter condition is only used to
fit \sref{eq:mr:main:spde}
into the framework of \cite{LiuRockner}.
\begin{itemize}
\item[(Hg)]{
  We have $g \in C^2(\Real^n; \Real^n)$
  with $g(u_-) = g(u_+) = 0$.
  In addition, there is $K_g > 0$ so that
  \begin{equation}
   \abs{Dg(u)} \le K_g
  \end{equation}
  holds for all $u \in \Real^n$, while
  \begin{equation}
    \abs{g(u_A) - g(u_B)}
    + \abs{Dg(u_A) - Dg(u_B)} \le K_g \abs{u_A - u_B}
  \end{equation}
  holds for all pairs $(u_A,u_B) \in \Real^n \times \Real^n$.
}
\end{itemize}

We remark here that it is advantageous to view
SPDEs as evolutions on Hilbert spaces, since
powerful tools are available in this setting.
However, in the case where $u_- \neq u_+$,
the waveprofile $\Phi_0$ does not lie in the natural
statespace $L^2$. In order to circumvent this problem,
we use $\Phi_0$ as a reference function
that connects $u_-$ to $u_+$, allowing
us to measure deviations from this
function in the Hilbert spaces $H^1$ and $L^2$.

In order to highlight this dual role and prevent
any confusion,
we introduce the duplicate notation
\begin{equation}
\Phi_{\mathrm{ref}} = \Phi_0
\end{equation}
and emphasize the fact that $\Phi_{\mathrm{ref}}$
remains fixed in the original frame,
unlike the wave-solution \sref{eq:mr:trv:wave:ansatz}.
We also introduce the sets
\begin{equation}
\mathcal{U}_{L^2} = \Phi_{\mathrm{ref}} + L^2,
\qquad
\mathcal{U}_{H^1} = \Phi_{\mathrm{ref}} + H^1,
\qquad
\mathcal{U}_{H^2} = \Phi_{\mathrm{ref}} + H^2,
\end{equation}
which we will use as the relevant state-spaces
to capture the solutions $U$ to \sref{eq:mr:main:spde}.

We now set out
to append a phase-tracking SDE
to \sref{eq:mr:main:spde}. In the deterministic case, we would couple the PDE to a phase-shift $\gamma$
that solves an ODE of the form
\begin{align}
    \dot\g(t)=c_0+  \O  \Big( U(t) - \Phi_0\big( \cdot - \gamma(t) \big) \Big) . 
\end{align}
By tuning the forcing function it is possible to remove the non-decaying terms in the original PDE,
which act in the direction of $\Phi_0'\big(\cdot - \gamma(t) \big)$.
This allows a nonlinear stability argument to be closed; see e.g. \cite{Zumbrun2009}.

In this paper we extend this procedure by
introducing a phase-shift $\Gamma$ that experiences the stochastic forcing
\begin{align}
    d\G= \Big[c_\s+  \O \Big( U(t) - \Phi_\sigma\big( \cdot - \Gamma(t) \big) \Big) \Big]dt
      + \O(\sigma) d\b_t.
\end{align}
By choosing the function $\Phi_\sigma$, the scalar $c_\sigma$
and the two forcing functions in an appropriate fashion,
the dangerous neutral terms can be eliminated from the original SPDE.
These are hence purely technical considerations,
but in \S\ref{sec:mr:disc} we discuss how these choices
can be related to quantities that are interesting from
an applied point of view.

In order to define our forcing functions in a fashion
that is globally Lipschitz continuous,
we introduce
the constant
\begin{equation}
\label{eq:mr:def:k:ip}
K_{\mathrm{ip}} = \big[\norm{g(\Phi_{0})}_{L^2} + 2 K_g \big] \norm{\psi_{\mathrm{tw}}}_{L^2}.
\end{equation}
In addition, we pick
two $C^\infty$-smooth non-decreasing cut-off functions
\begin{equation}
\chi_{\mathrm{low}}: \Real \to [\frac{1}{4}, \infty),
\qquad
\chi_{\mathrm{high}}: \Real \to
[- K_{\mathrm{ip}} - 1, K_{\mathrm{ip}} + 1]
\end{equation}
that satisfy the identities
\begin{equation}
\chi_{\mathrm{low}}(\vartheta) = \frac{1}{4} \hbox{ for } \vartheta \le \frac{1}{4},
\qquad
\chi_{\mathrm{low}}(\vartheta) = \vartheta \hbox{ for } \vartheta \ge \frac{1}{2},
\end{equation}
together with
\begin{equation}
\chi_{\mathrm{high}}(\vartheta) = \vartheta \hbox{ for } \abs{\vartheta} \le K_{\mathrm{ip}},
\qquad
\chi_{\mathrm{high}}(\vartheta) = \mathrm{sign}(\vartheta) \big[K_{\mathrm{ip}} + 1]
   \hbox{ for } \abs{\vartheta} \ge K_{\mathrm{ip}} + 1.
\end{equation}

For any $u \in \mathcal{U}_{H^1}$
and $\psi \in H^1$,
this allows us to introduce
the functions
\begin{equation}
\label{eq:mr:def:b:kappa}
\begin{array}{lcl}
b(u, \psi)
 & = &
  - \Big[
    \chi_{\mathrm{low}}\big(
    \langle \partial_\xi u ,
   \psi \rangle_{L^2} \big) \Big]^{-1}
     \chi_{\mathrm{high}}\big(
        \langle g(u) , \psi \rangle_{L^2}
     \big) ,
\\[0.2cm]
\kappa_{\sigma}(u, \psi)
& = & 1 + \frac{1}{2\rho} \sigma^2 b( u, \psi )^2.
\end{array}
\end{equation}
In addition, for any
$u \in \mathcal{U}_{H^1}$, $c \in \Real$ and $\psi \in H^1$
we define the expression
\begin{equation}
\label{eq:mr:def:j}
\begin{array}{lcl}
\mathcal{J}_{\sigma}(u , c, \psi)
& = &
  \kappa_{\sigma}(u , \psi)^{-1}
  \Big[
    f(u )
    +c u'
    + \sigma^2 b( u, \psi)
      \partial_\xi [g(u )]
  \Big] ,
\\[0.2cm]
\end{array}
\end{equation}
while for any
$u \in \mathcal{U}_{H^1}$, $c \in \Real$ and $\psi \in H^2$
we write
\begin{equation}
\begin{array}{lcl}
\label{eq:mr:def:a}
a_{\sigma}(u,  c, \psi )
 & = & - \kappa_{\sigma}(u,\psi)
    \Big[ \chi_{\mathrm{low}}\big( \langle \partial_\xi u, \psi \rangle_{L^2} \big) \Big]^{-1}
 \Big[
    \langle u , A_* \psi \rangle_{L^2}
    + \langle
       \mathcal{J}_{\sigma}(u , c, \psi ),
     \psi \rangle_{L^2}
 \Big] .
\end{array}
\end{equation}
Finally, we introduce the
right-shift operators
\begin{equation}
[T_{\gamma} u](\xi)
= u(\xi - \gamma)
\end{equation}
that act on any function $u: \Real \to \Real^n$.

With these ingredients in hand, we are
ready to introduce the main SPDE that we analyze
in this paper.
We formally write this SPDE
as the skew-coupled system\footnote{
Note here that formally $b(U, T_{\Gamma} \psi_{\mathrm{tw}})$ is a multiplication
operator from $\R\to\R$, hence a number.
If we generalize $\b_t$ to a cylindrical $Q-$Wiener process on a space $H$
then the term involving $b$ becomes a functional from $H$ to $\R$. Fortunately,
all the relevant technical machinery that we use
has been generalized to this setting. In particular, our approach
carries over but the notation becomes significantly more convoluted.
}
\begin{equation}
\label{eq:mr:formal:spde}
\begin{array}{lcl}
dU &= &
  \big[ A_* U + f(U) \big] dt
  + \sigma g(U) d \beta_t,
\\[0.2cm]
d\Gamma & = &
  \big[ c + a_\sigma \big( U ,c, T_{\Gamma} \psi_{\mathrm{tw}} \big)
  \big] dt
  + \sigma b\big(U , T_{\Gamma} \psi_{\mathrm{tw}} \big)
     \, d \beta_t ,
\end{array}
\end{equation}
noting that we seek solutions
with
$\big(U(t), \Gamma(t) \big) \in \mathcal{U}_{H^1} \times \Real$.
Observe that the first equation is the same as \sref{eq:mr:main:spde}.

In order to make this precise,
we introduce the spaces
\begin{equation}
\begin{array}{lcl}
\mathcal{N}^2\big( [0,T] ; (\mathcal{F}_t ) ;
   \mathcal{H} \big)
& = & \{ X \in L^2\big( [0 , T] \times \Omega ;
  dt \otimes \mathbb{P} ;  \mathcal{H} \big)
    :
\\[0.2cm]
& & \qquad \qquad X \hbox{ has a }
    (\mathcal{F}_t)\hbox{-progressively measurable version}
 \},
\end{array}
\end{equation}
where we allow
$\mathcal{H} \in \{ \Real, L^2, H^1 \}$.
We note that we follow
the convention of \cite{Concise,revuz2013continuous}
here by requiring progressive measurability
instead of the usual stronger notion of predictability.
Since we are exclusively dealing with Brownian motions,
this choice suffices to construct
stochastic integrals.

Our first result clarifies
what we mean by a solution
to \sref{eq:mr:formal:spde}.
We note that (i) and (ii) in Proposition \ref{prp:mr:main:ex} imply that $(X, \Gamma)$
is an
$L^2 \times \Real$-valued continuous $(\mathcal{F}_t)$-adapted process.
We remark that in the integral equation \sref{eq:mr:main:ex:itg:eq:x}
we interpret the diffusion operator $A_*$ as an
element of $\mathcal{L}( H^1 ; H^{-1})$,
where $H^{-1}$ is the dual of $H^1$ under the standard
embeddings
\begin{equation}
H^1 \hookrightarrow L^2 \cong \big[L^2\big]^* \hookrightarrow H^{-1} = \big[H^{1}\big]^*.
\end{equation}
We note that the set $(H^1, L^2, H^{-1})$ is commonly referred to as
a Gelfand triple; see e.g. \cite[{\S}5.9]{Evans} for a more detailed explanation.
For $(v, w) \in H^{-1} \times H^1$ we write $\langle v, w \rangle_{H^{-1}; H^1}$
to refer to the duality pairing between $H^1$ and $H^{-1}$. If in fact $v \in L^2$, then we have
\begin{equation}
\langle v, w \rangle_{H^{-1}; H^1} = \langle v, w \rangle_{L^2}.
\end{equation}


\begin{prop}[{see \S\ref{sec:var}}]
\label{prp:mr:main:ex}
Suppose that $(HA)$, $(Hf)$, $(HVar)$,
$(HTw)$, $(HS)$, $(Hg)$ and $(H\beta)$ are all satisfied
and fix $T > 0$, $c \in \Real$ and $0 \le \sigma \le 1$.
In addition, pick an initial condition
\begin{equation}
  (X_0, \Gamma_0) \in L^2 \times \Real.
\end{equation}
Then there are maps
\begin{equation}
X: [0, T] \times \Omega \to L^2,
\qquad
\Gamma: [0,T] \times \Omega \to \Real
\end{equation}
that satisfy the following properties.
\begin{itemize}
\item[(i)]{
  For almost all $\omega \in \Omega$,
  the map
  \begin{equation}
  t \mapsto \big( X(t, \omega) , \Gamma(t, \omega) \big)
  \end{equation}
  is of class $C([0,T]; L^2 \times \Real)$.
}
\item[(ii)]{
  For all $t \in [0,T]$, the map
  \begin{equation}
  \omega \mapsto \big( X(t, \omega), \Gamma(t, \omega) \big)
    \in L^2 \times \Real
  \end{equation}
  is $(\mathcal{F}_t)$-measurable.
}
\item[(iii)]{
  We have the inclusion
  \begin{equation}\label{eq:mr:inclL2}
  X \in L^6\big( \Omega, \mathbb{P} ; C([0,T]; L^2 ) \big),
  \end{equation}
  together with
  \begin{equation}
  \label{eq:prp:mr:ex:inclusions:n2:x:gamma}
  \begin{array}{lcl}
     X & \in &  \mathcal{N}^2 \big( [0 , T] ; (\mathcal{F}_t ) ; H^1 \big) ,
  \\[0.2cm]
     \Gamma & \in &  \mathcal{N}^2 \big( [0 , T] ; (\mathcal{F}_t ) ; \Real \big)
  \end{array}
  \end{equation}
  and
  \begin{equation}
  \label{eq:mr:main:ex:inclusions:g:b}
  \begin{array}{lcl}
     g(X + \Phi_{\mathrm{ref}})
      & \in &
       \mathcal{N}^2 \big( [0 , T] ; (\mathcal{F}_t ) ; L^2 \big),
     \\[0.2cm]
     b\big(X + \Phi_{\mathrm{ref}}, T_{\Gamma} \psi_{\mathrm{tw}} \big)
      & \in &
       \mathcal{N}^2 \big( [0 , T] ; (\mathcal{F}_t ) ; \Real \big) .
  \end{array}
  \end{equation}
}
\item[(iv)]{
  For almost all $\omega \in \Omega$,
  the identities
  \begin{equation}
  \label{eq:mr:main:ex:itg:eq:x}
  \begin{array}{lcl}
  X(t) & = & X_0 + \int_0^t A_*[ X(s) + \Phi_{\mathrm{ref}} ] \, ds
    + \int_0^t f\big(X(s) + \Phi_{\mathrm{ref}}\big) \, ds
   \\[0.2cm]
   & & \qquad
    + \sigma \int_0^t  g\big(X(s) + \Phi_{\mathrm{ref}}\big) d \beta_s
  \end{array}
  \end{equation}
  and
  \begin{equation}
  \begin{array}{lcl}
   \Gamma(t) & = & \Gamma_0
    + \int_0^t
      \big[ c + a_{\sigma}\big(X(s) + \Phi_{\mathrm{ref}}, c , T_{\Gamma(s)}\psi_{\mathrm{tw}} \big) \big] \, ds
  \\[0.2cm]
  & & \qquad
     + \sigma \int_0^t
        b\big(X(s) + \Phi_{\mathrm{ref}} , T_{\Gamma(s)} \psi_{\mathrm{tw}} \big)
         d \beta_s
  \end{array}
  \end{equation}
  hold\footnote{Note that this equation initially only holds as an identity in $H^{-1}$. Inclusion \sref{eq:mr:inclL2} makes that we can interpret the integrals in $L^2$. We have $X  \in   \mathcal{N}^2 \big( [0 , T] ; (\mathcal{F}_t ) ; H^1 \big)$ but this does not mean that $X(t)\in H^1$ pointwise.} for all $0 \le t \le T$.
}
\item[(v)]{Suppose that the pair $(\tilde{X}, \tilde{\Gamma} ): [0,T] \times \Omega \to L^2 \times \Real$
also satisfies (i)-(iv). Then for almost all $\omega \in \Omega$,
we have
\begin{equation}
(\tilde{X}, \tilde{\Gamma})(t) = (X, \Gamma)(t) \qquad \hbox{ for all } 0 \le t \le T.
\end{equation}
}
\end{itemize}
\end{prop}

\subsection{Wave stability}
\label{sec:mr:wave}

By inserting the traveling wave Ansatz
\sref{eq:mr:trv:wave:ansatz} into the deterministic PDE
\sref{eq:mr:det:pde},
we observe that
\begin{equation}
A_* \Phi_0 + \mathcal{J}_0(\Phi_0, c_0, \psi_{\mathrm{tw}}) = 0,
\end{equation}
which means that $a_0(\Phi_0, c_0, \psi_{\mathrm{tw}}) = 0$.
Our first result here shows that this can be extended
into a branch of profiles and speeds
for which
\begin{equation}
a_{\sigma}(\Phi_{\sigma}, c_{\sigma}, \psi_{\mathrm{tw}}) = 0.
\end{equation}
Roughly speaking,
this means that
the adjusted phase $\Gamma(t) - c t$
will (instantaneously) feel only stochastic forcing
if one takes $c = c_{\sigma}$
and $U = T_{\Gamma(t)} \Phi_{\sigma}$
in \sref{eq:mr:formal:spde}.

\begin{prop}[{see \S\ref{sec:swv}}]
\label{prp:mr:swv:ex}
Suppose that $(HA)$, $(Hf)$,
$(HTw)$, $(HS)$ and $(Hg)$
are all satisfied and pick a sufficiently large
constant $K > 0$. Then there exists $\delta_{\sigma} > 0$
so that 
for every $0 \le \sigma \le \delta_{\sigma}$,
there is a unique pair
\begin{equation}
(\Phi_{\sigma},c_{\sigma}) \in   \mathcal{U}_{H^2}  \times \Real
\end{equation}
that satisfies the system
\begin{equation}
\label{eq:mr:prop:swv:eq}
  A_* \Phi_{\sigma} +
   \mathcal{J}_{\sigma}\big(\Phi_{\sigma} , c_{\sigma}, \psi_{\mathrm{tw}} \big) = 0
\end{equation}
and admits the bound
\begin{equation}
\label{eq:mr:prp:swv:bnd:phi:sigma}
\norm{\Phi_{\sigma} - \Phi_0}_{H^2} + \abs{c_{\sigma} - c_{0} } \le K \sigma^2 .
\end{equation}
\end{prop}

We are interested in solutions
to \sref{eq:mr:formal:spde} with
an initial condition for $U$
that is close to $\Phi_{\sigma}$.
We will use the remaining degree of freedom
to pick the initial phase
$\Gamma$ in such a way that the orthogonality
condition described in the following result
is enforced.

\begin{prop}[{see \S\ref{sec:swv}}]
\label{prp:mr:phase:shift}
Suppose that (HA), (Hf), (HTw), (HS) and (Hg) are all satisfied.
Then there exist constants $\delta_0 > 0$, $\delta_{\sigma} > 0$
and $K > 0$
so that the following holds true.
For every $0 \le \sigma \le \delta_{\sigma}$
and any $u_0 \in \mathcal{U}_{L^2} $ that has
\begin{equation}
\label{eq:swv:phaseshift:small:cond:u:only}
\norm{u_0 - \Phi_{\sigma} }_{L^2} < \delta_0,
\end{equation}
there exists
$\gamma_0 \in \Real$
for which the function
\begin{equation}
v_{\gamma_0} = T_{-\gamma_0} [u_0] - \Phi_{\sigma}
\end{equation}
satisfies
the identity
\begin{equation}
\langle v_{\gamma_0} , \psi_{\mathrm{tw}} \rangle_{L^2} = 0
\end{equation}
together with the bound
\begin{equation}
\label{eq:swv:bnds:on:shift:final}
\abs{\gamma_0} + \norm{v_{\gamma_0}}_{L^2} \le
   K \norm{ u_0 - \Phi_{\sigma} }_{L^2} .
\end{equation}
If in fact $u_0 \in \mathcal{U}_{H^1}$,
then we also have the estimate
\begin{equation}
\label{eq:swv:bnds:on:shift:final:h1}
\abs{\gamma_0} + \norm{v_{\gamma_0}}_{H^1} \le
   K \norm{ u_0 - \Phi_{\sigma} }_{H^1} .
\end{equation}
\end{prop}

Let us now pick any $u_0 \in \mathcal{U}_{H^1}$
for which \sref{eq:swv:phaseshift:small:cond:u:only}
holds. We write $(X_{u_0} , \Gamma_{u_0})$
for the process described in Proposition
\ref{prp:mr:main:ex}
with the initial condition
\begin{equation}
(X_0, \Gamma_0) = \big( u_0 - \Phi_{\mathrm{ref}}, \gamma_0 \big),
\end{equation}
in which $\gamma_0$ is the initial phase
defined in Proposition \ref{prp:mr:phase:shift}.
We then define the process
\begin{equation}
\label{eq:mr:def:v:u0}
V_{u_0}(t) = T_{-\Gamma_{u_0}(t) }
  \big[ X_{u_0}(t) + \Phi_{\mathrm{ref}} \big] - \Phi_{\sigma} ,
\end{equation}
which can be thought of as the
deviation of the solution
$U$ of \sref{eq:mr:formal:spde}
from the stochastic wave
$\Phi_{\sigma}$ shifted to the position $\Gamma_{u_0}(t)$.

In order to measure the size
of the perturbation,
we pick $\e > 0$
and introduce the scalar function
\begin{equation}
N_{\e;u_0} (t) =
\norm{V_{u_0}(t)}_{L^2}^2
 + \int_0^t e^{- \e (t - s) }
    \norm{  V_{u_0}(s)}_{H^1}^2 \, ds .
\end{equation}
For each $T > 0$ we now define
a probability
\begin{equation}
p_{\e}(T, \eta, u_0) = P\Big(
 \sup_{0 \le t \le T}
 N_{\e;u_0}(t) > \eta
\Big) .
\end{equation}
Our first main result shows that
the probability that
$N_{\e;u_0}$ remains
small on timescales of order $\sigma^{-2}$ can be pushed arbitrarily close to one
by restricting the strength of the noise
and the size of the initial perturbation.

\begin{thm}[{see \S\ref{sec:nls}}]
\label{thm:mr:orbital:stb}
Suppose that $(HA)$, $(Hf)$, $(HVar)$,
$(HTw)$, $(HS)$, $(Hg)$ and $(H\beta)$ are all satisfied
and pick sufficiently small constants $\e > 0$, $\delta_0>0$,
$\delta_{\eta} > 0$ and $\delta_{\sigma} > 0$.
Then there exists a constant $K > 0$
so that for every $T> 1$, any $0 \le \sigma \le \delta_{\sigma}T^{-1/2}$,
any $u_0 \in \mathcal{U}_{H^1}$
that satisfies \sref{eq:swv:phaseshift:small:cond:u:only}
and any $0 < \eta \le \delta_{\eta}$,
we have the inequality
\begin{equation}
p_{\e}(T, \eta, u_0)
  \le   \eta^{-1} K \Big[ \norm{u_0 - \Phi_{\sigma}}_{H^1}^2
     + \sigma^2 T \Big].
\end{equation}
\end{thm}

Our second main result concerns the
special case
where the noise pushes the stochastic wave
$\Phi_{\sigma}$ in a rigid fashion.
This is the case when
\begin{equation}
\label{eq:mr:prop:g:vs:phi:prime}
g(\Phi_0) = \vartheta_0 \Phi_0'
\end{equation}
for some proportionality constant
$\vartheta_0 \in \Real$.
It is easy to verify that
\sref{eq:mr:prop:g:vs:phi:prime}
with $\vartheta_0 = - \sqrt{2}$
holds for \sref{eq:int:SPDE:nag}.

In this setting we expect the perturbation $V$ to decay
exponentially on timescales of order $\sigma^{-2}$ with large probability.
In order to formalize this, we pick small constants
$\e > 0$ and $\alpha > 0$
and introduce the scalar function
\begin{equation}
N_{\e,\alpha;u_0}(t) =
e^{\alpha t} \norm{V_{u_0}(t)}_{L^2}^2
 + \int_0^t e^{- \e (t - s) }
    e^{\alpha s} \norm{  V_{u_0}(s)}_{H^1}^2 \, ds ,
\end{equation}
together with the associated probabilities
\begin{equation}
\label{eq:mr:def:p:eps:alpha}
p_{\e,\alpha}(T,\eta, u_0) = P\Big(
 \sup_{0 \le t \le T}  N_{\e,\alpha;u_0}(t)
 > \eta
\Big) .
\end{equation}

\begin{thm}[{see \S\ref{sec:nls}}]
\label{thm:mr:exp:stb}
Suppose that $(HA)$, $(Hf)$, $(HVar)$,
$(HTw)$, $(HS)$, $(Hg)$ and $(H\beta)$ are all satisfied .
Suppose furthermore that
\sref{eq:mr:prop:g:vs:phi:prime} holds
and pick sufficiently small constants $\e > 0$, $\delta_0>0$,
$\alpha > 0$,
$\delta_{\eta} > 0$ and $\delta_{\sigma} > 0$.
Then there exists a constant $K > 0$
so that for any $T > 1$,
every $0 \le \sigma \le \delta_{\sigma}T^{-1/2}$,
any $u_0 \in \mathcal{U}_{H^1}$
that satisfies \sref{eq:swv:phaseshift:small:cond:u:only}
and any $0 < \eta \le \delta_{\eta}$,
we have the inequality
\begin{equation}
 p_{\e,\alpha}(T, \eta, u_0)
  \le   \eta^{-1} K\norm{u_0 - \Phi_{\sigma}}_{H^1}^2 .
\end{equation}
\end{thm}

\subsection{Interpretation}
\label{sec:mr:disc}

In \S\ref{sec:sps}
we show that
the pair $(V, \Gamma) = (V_{u_0} , \Gamma_{u_0} )$
defined in \S\ref{sec:mr:wave}
satisfies the SPDE
\begin{equation}
\label{eq:mr:disc:main:spde}
\begin{array}{lcl}
d V  & = & \mathcal{R}_{\sigma}(V) \, dt + \sigma \mathcal{S}_{\sigma}(V) \, d \beta_t ,
\\[0.2cm]
d \Gamma & = &
  \big[
  c_{\sigma}
+ a_{\sigma}\big(\Phi_{\sigma} + V, c_{\sigma} ,
    \psi_{\mathrm{tw}} \big) \big] \, dt
    + \sigma b \big(\Phi_{\sigma} + V ,
          \psi_{\mathrm{tw}} \big) \, d \beta_t ,
\end{array}
\end{equation}
in which the nonlinearities satisfy
the identities
\begin{equation}
\label{eq:mr:disc:ids}
a_{\sigma}\big( \Phi_{\sigma} , c_{\sigma}, \psi_{\mathrm{tw}} \big) = 0,
\qquad
\mathcal{R_{\sigma} }(0) = 0,
\qquad
\mathcal{S_{\sigma} }(0) = g(\Phi_{\sigma}) + b (\Phi_{\sigma} ) \Phi_{\sigma}' ,
\end{equation}
together with the asymptotics
\begin{equation}
\label{eq:mr:disc:asymp:nl}
D_1 a_{\sigma}\big( \Phi_{\sigma} , c_{\sigma}, \psi_{\mathrm{tw}} \big)
= \O(\sigma^2),
\qquad
D\mathcal{R}_{\sigma}(0) = \O(\sigma^2) .
\end{equation}
For our discussion here  we take $V(0) = 0$ and $\Gamma(0) = 0$,
which corresponds with the initial condition
condition $U(0) = \Phi_{\sigma}$
for the original system \sref{eq:mr:main:spde}.

The identities \sref{eq:mr:disc:ids}
imply that $V(t)$ and $\Gamma(t) - c_{\sigma} t$
experience no deterministic forcing at
$t = 0$. We now briefly discuss
the consequences of this observation
on the behaviour of \sref{eq:mr:disc:main:spde}
in the two regimes
covered by Theorems \ref{thm:mr:orbital:stb} and
\ref{thm:mr:exp:stb}.

\paragraph{Exponential stability}
Our results are easiest to interpret
in the special case
\begin{equation}
g(\Phi_0) = \vartheta_0 \Phi_0',
\end{equation}
where Theorem \ref{thm:mr:exp:stb} applies.
Remarkably, the modified profiles and speeds
$(\Phi_{\sigma} , c_{\sigma})$
can be computed explicitly in this setting. \footnote{ Similar formula can be obtained by following the formal approach in \cite{Cartwright2019}, which appeared during the revision phase of this paper.}

\begin{prop}[{see \S\ref{sec:swv}}]
\label{prp:mr:expl:waves}
Consider the setting of Proposition \ref{prp:mr:swv:ex}
and suppose that \sref{eq:mr:prop:g:vs:phi:prime} holds. Then for all sufficiently small $0 \le \sigma \le \delta_{\sigma}$
we have the identities
\begin{equation}
\label{eq:mr:id:for:mod:wave:exp:stb}
\begin{array}{lcl}
   \Phi_{\sigma}(\xi)
  & = & \Phi_0\Big(
    \big[ 1 + \frac{1}{2\rho} \sigma^2 \vartheta_{0}^2 \big]^{1/2}
      \xi \Big) ,
\\[0.2cm]
  c_{\sigma} & = &
    \big[ 1 + \frac{1}{2\rho} \sigma^2 \vartheta_{0}^2 \big]^{-1/2}
       c_0 ,
\end{array}
\end{equation}
together with
\begin{equation}
\label{eq:mr:id:for:g:expl:case}
g(\Phi_{\sigma} ) =
 \big[ 1 + \frac{1}{2\rho} \sigma^2 \vartheta_{0}^2 \big]^{-1/2}
   \vartheta_0 \Phi_{\sigma}'
 = -b\big( \Phi_{\sigma}, \psi_{\mathrm{tw} } \big)
 \Phi_{\sigma}'.
\end{equation}
\end{prop}

A direct consequence of \sref{eq:mr:id:for:g:expl:case}
is that the identity
\begin{equation}
\mathcal{S}_{\sigma}(0) = 0
\end{equation}
can be added to the list \sref{eq:mr:disc:ids}.
In particular,
we obtain the explicit solution
\begin{equation}
\begin{array}{lcl}
(V, \Gamma) & = &
\Big( 0 , c_{\sigma} t + \sigma b(\Phi_{\sigma} , \psi_{\mathrm{tw}} \big) \beta_t \Big)
\\[0.2cm]
& = &
\Big( 0 , c_{\sigma} t - \sigma \big[ 1 + \frac{1}{2 \rho} \sigma^2 \vartheta_0^2]^{-1/2} \vartheta_0 \beta_t \Big)
\end{array}
\end{equation}
for the system \sref{eq:mr:disc:main:spde}.
This corresponds to the solution
\begin{equation}
U(t) = \Phi_{\sigma}(\cdot + \Gamma(t))
\end{equation}
for \sref{eq:mr:formal:spde}, which exists
for all $t \ge 0$.

We hence see that the shape $\Phi_{\sigma}$
of the stochastic profile remains fixed,
while the phase $\Gamma(t)$ of the wave performs
a scaled Brownian motion around the position $c_{\sigma} t$.
Since the identities \sref{eq:mr:id:for:mod:wave:exp:stb}
imply that the waveprofile is steepened
while the speed is slowed down,
our results indeed confirm the numerical observations
from \cite{Lord2012} that were discussed in \S\ref{sec:int}.

Any small perturbation in the $V$ component will decay exponentially fast with
high probability
on account of Theorem \ref{thm:mr:exp:stb}. Intuitively,
the leading order behaviour for $V$ resembles a geometric Brownian motion,
as the noise term is proportional to $V$ while the deterministic forcing
leads to exponential decay. In particular, we expect
that our approach can keep track of the wave for timescales that
are far longer than the $\O(\sigma^{-2})$ bounds stated in our results.

\paragraph{Orbital stability}
In general we have $\mathcal{S}_{\sigma}(0) \neq 0$,
which prevents us from
solving \sref{eq:mr:disc:main:spde} explicitly.
Indeed, Theorem \ref{thm:mr:orbital:stb} states that
$V(t)$ will remain small with high probability,
but the stochastic forcing will preclude it from
converging to zero. However, our construction does guarantee
that $\langle V(t), \psi_{\mathrm{tw}} \rangle = 0$
as long as $V$ stays small. Since $\langle \Phi_0' , \psi_{\mathrm{tw}} \rangle = 1$,
this still allows us to interpret
$\G(t)$ as the position of the wave.
In particular, if the expression
$t^{-1} \Gamma(t)$ converges in a suitable sense as $t \to \infty$
then it is natural to use this limit as a proxy for the notion of a wavespeed.

In order to explore this, we introduce the formal expansion
\begin{equation}
V(t) = \sigma V^{(1)}_{\sigma}(t) + \O (\sigma^2)
\end{equation}
and use the mild formulation
developed in \S\ref{sec:md}
to obtain
\begin{equation}
\label{eq:mr:process:for:v}
V^{(1)}_{\sigma}(t) = \int_0^t S (t - s) \mathcal{S}_{\sigma}(0) \, d \beta_s .
\end{equation}
Here $S$ denotes the semigroup generated by $\mathcal{L}_{\mathrm{tw}}$,
which by construction decays exponentially when applied to $\mathcal{S}_{\sigma}(0)$.
In particular, for any
bilinear map $B: H^1 \times H^1 \to \Real$ we can use the It{\^o} isometry
to obtain
\begin{equation}
\begin{array}{lcl}
E B \big[
  V^{(1)}_{\sigma}(t) ,V^{(1)}_{\sigma}(t)
\big]
  & = &  \int_0^t  B \big[ S(t-s) \mathcal{S}_{\sigma}(0) , S(t-s) \mathcal{S}_{\sigma}(0) \big] \, ds
\\[0.2cm]
  & = &  \int_0^t  B \big[ S(s) \mathcal{S}_{\sigma}(0) , S(s) \mathcal{S}_{\sigma}(0) \big] \, ds ,
\end{array}
\end{equation}
which converges in the limit $t \to \infty$.

Introducing the formal expansion
\begin{equation}
\Gamma(t) = c_{\sigma} t + \sigma \Gamma^{(1)}_{\sigma}(t)
+ \sigma^2 \Gamma^{(2)}_{\sigma}(t) + \O(\sigma^3),
\end{equation}
the first bound in \sref{eq:mr:disc:asymp:nl}
implies that
\begin{equation}
\begin{array}{lcl}
\Gamma^{(1)}_{\sigma}(t)
 & = &  b(\Phi_{\sigma}, \psi_{\mathrm{tw}}) \beta_t
 \end{array}
\end{equation}
together with
\begin{equation}
\begin{array}{lcl}
\Gamma^{(2)}_{\sigma}(t)
  & = &
     \frac{1}{2}\int_0^t D^2_1 a_{\sigma}
     \big(\Phi_{\sigma} , c_{\sigma}, \psi_{\mathrm{tw}}  \big)
      \big[V_\s^{(1)}(s), V_\s^{(1)}(s) \big] \, ds
\\[0.2cm]
 & & \qquad
   +  D_1b\big(\Phi_{\sigma} , \psi_{\mathrm{tw}} \big)
   \left[\int_0^t V_\s^{(1)}(s) \, d \beta_s\right] .
\end{array}
\end{equation}
Since $E V^{(1)}_{\sigma}(t) = 0$ we obtain
\begin{equation}
 E \Gamma^{(1)}_{\sigma}(t) =  0 
\end{equation}
together with
\begin{equation}
\begin{array}{lcl}
E \Gamma^{(2)}_{\sigma}(t)
  & = & \frac{1}{2}
\int_0^t \int_{0}^s D^2_1 a_{\sigma}
     \big(\Phi_{\sigma} , c_{\sigma}, \psi_{\mathrm{tw}}  \big)
      \big[ S(s') \mathcal{S}_{\sigma}(0), S(s') \mathcal{S}_{\sigma}(0)  \big]
      \, ds' \, ds
\\[0.2cm]
& = & \frac{1}{2} \int_0^t (t - s) D^2_1 a_{\sigma}
     \big(\Phi_{\sigma} , c_{\sigma}, \psi_{\mathrm{tw}}  \big)
      \big[ S(s) \mathcal{S}_{\sigma}(0), S(s) \mathcal{S}_{\sigma}(0)  \big]
       \, ds .
\\[0.2cm]
\end{array}
\end{equation}
Upon writing
\begin{equation}
\label{eq:mr:c:inft:order:two}
c_{\infty}^{(2)} = c_{\sigma} +
\frac{\sigma^2}{2} \int_0^\infty
D^2_1 a_{\sigma}
     \big(\Phi_{\sigma} , c_{\sigma}, \psi_{\mathrm{tw}}  \big)
      \big[S(s) \mathcal{S}_{\sigma}(0) , S(s) \mathcal{S}_{\sigma}(0) \big] \, ds,
\end{equation}
we hence conjecture that the expected limiting value
of the wavespeed behaves as
$c_{\infty}^{(2)} +\O (\sigma^3)$.
Since $c_{\sigma} = c_0 + \O(\sigma^2)$ this would mean that the stochastic contributions
to the wavespeed are second order in $\sigma$.

We remark that computations of this kind resemble the multi-scale
approach initiated by Lang in \cite{Lang} and Stannat
and Kr\"uger in \cite{kruger2017multiscale}.
However, our approach does allow us to
consider limiting expressions such as \sref{eq:mr:c:inft:order:two},
for which one needs the exponential decay of the semigroup.
Indeed, \sref{eq:mr:process:for:v}
resembles a mean-reverting Ornstein-Uhlenbeck process,
which has a variation that can be globally bounded in time,
despite the fact that the individual paths are unbounded.

As above, we expect
to be able to track the wave for timescales that
are longer than the $\O(\sigma^{-2})$ bounds stated in our results.
The key issue is that the mild version of the Burkholder-Davis-Gundy
inequality that we use is not able to fully incorporate the mean-reverting
effects of the semigroup.
We emphasize that even the standard scalar Ornstein-Uhlenbeck process
requires sophisticated probabilistic machinery to uncover
statistics concerning the behaviour of the running supremum
\cite{ricciardi1988first,alili2005representations}.
We plan to explore these issues
in more detail in a forthcoming paper.
For the moment however, we note that our 
initial numerical experiments seem to confirm that the
expression \sref{eq:mr:c:inft:order:two} indeed captures the leading order stochastic correction
to the wavespeed.

\section{Preliminary estimates}
\label{sec:prlm}

In this section we derive several preliminary estimates
for the functions $f$,  $g$, $\mathcal{J}_0$, $b$ and $\kappa_{\sigma}$.
We will write the arguments
$(u, \overline{c}) \in \mathcal{U}_{H^1} \times \Real$
as
\begin{equation}
u = \Phi + v, \qquad \qquad \overline{c} = c + d,
\end{equation}
in which we take $(\Phi, c) \in \mathcal{U}_{H^1} \times \Real$
and $(v, d) \in H^1 \times \Real$.
We do not restrict ourselves to the case where $(\Phi,c) =(\Phi_0, c_0)$,
but impose the following condition.
\begin{itemize}
\item[(hPar)]{
 The conditions (HTw) and (HS) hold
 and the pair $(\Phi,c) \in \mathcal{U}_{H^1} \times \Real$
 satisfies the bounds
 \begin{equation}
    \norm{\Phi - \Phi_{0} }_{H^1} \le
      \min\{ 1 ,[4 \norm{\psi_{\mathrm{tw}}}_{L^2}]^{-1} \},
    \qquad \qquad
    \abs{c - c_0} \le 1.
  \end{equation}
}
\end{itemize}
In \S\ref{sec:prlm:fg} we obtain global and Lipschitz bounds
for the functions $f$ and $g$. These bounds
are subsequently used in \S\ref{sec:prlm:rest}
to analyze the auxiliary functions
$\mathcal{J}_0$, $b$ and $\kappa_{\sigma}$.
Throughout this paper we use the convention that all numbered constants appearing
in proofs are strictly positive and have the same dependencies as the
constants appearing in the statement of the result.

\subsection{Bounds on $f$ and $g$}
\label{sec:prlm:fg}

The conditions (Hf) and (Hg) allow us to obtain standard cubic
bounds on $f$
and globally Lipschitz bounds on $g$. We also consider
expressions of the form $\partial_{\xi} g(u)$,
which give rise to quadratic bounds.

\begin{lem}
\label{lem:prlm:f:bnds}
Suppose that (Hf) and (hPar) are satisfied.
Then there exists a constant $K > 0$, which does not
depend on the pair $(\Phi, c)$, so that the following holds true.
For any $v \in H^1$ and $\psi \in H^1$
we have the bounds
\begin{equation}
\label{eq:prlm:f:glb:bnds}
\begin{array}{lcl}
 \norm{f(\Phi + v)}_{L^2} & \le &
   K \big[1 + \norm{v}_{H^1}^2 \norm{v}_{L^2} \big] ,
\\[0.2cm]
 \abs{\langle f(\Phi + v), \psi \rangle_{L^2} }
   & \le &
     K \big[1  + \norm{v}_{H^1} \norm{v}_{L^2}^2 \big]
       \norm{\psi}_{H^1} ,
\\[0.2cm]
\end{array}
\end{equation}
while for any set of pairs $(v_A, v_B) \in H^1 \times H^1$
and $(\psi_A, \psi_B) \in H^1 \times H^1$,
the expressions
\begin{equation}
\begin{array}{lcl}
\Delta_{AB} f & = & f(\Phi + v_A) - f(\Phi + v_B) ,
\\[0.2cm]
\Delta_{AB}\langle f, \cdot \rangle_{L^2}
  & = & \langle f(\Phi + v_A ) , \psi_A \rangle_{L^2}
  - \langle f(\Phi + v_B) , \psi_B \rangle_{L^2}
\end{array}
\end{equation}
satisfy the estimates
\begin{equation}
\label{eq:prlm:f:lip:bnds}
\begin{array}{lcl}
\norm{ \Delta_{AB} f }_{L^2}
  & \le &
  K \norm{v_A-v_B}_{L^2}
\\[0.2cm]
& & \qquad
  + K
  \Big( \norm{v_A}_{H^1} \norm{v_A}_{L^2}
    + \norm{v_B}_{H^1} \norm{v_B}_{L^2} \Big)
   \norm{v_A-v_B}_{H^1} ,
\\[0.2cm]
\abs{\Delta_{AB}\langle f, \cdot \rangle_{L^2} }
 & \le &
 K \norm{v_A-v_B}_{L^2} \norm{\psi_A}_{H^1}
\\[0.2cm]
& & \qquad
  + K \norm{v_A-v_B}_{H^1} \big( \norm{v_A}_{L^2}^2 + \norm{v_B}_{L^2}^2 \big)
    \norm{\psi_A}_{H^1}
\\[0.2cm]
& & \qquad
  + K \Big[ 1
    + \norm{v_B}_{H^1} \norm{v_B}_{L^2}^2 \Big]
    \norm{\psi_A - \psi_B}_{H^1} .
\end{array}
\end{equation}
\end{lem}
\begin{proof}
Exploiting (Hf) we obtain
\begin{equation}
\label{eq:prlm:f:bnd:d2:f}
\abs{D^2 f(u)} \le C_1[ 1 + \abs{u} ],
\end{equation}
together with
\begin{equation}
\abs{Df(u)} \le C_1 [ 1 + \abs{u}^2 ]
\end{equation}
for all $u \in \Real^n$.
In particular, (hPar) yields
the pointwise Lipschitz bound
\begin{equation}
\label{eq:prlm:f:pointwise:lip}
\abs{f(\Phi + v_A) - f(\Phi + v_B) }
\le C_2[ 1 +
    \abs{v_A}^2 +  \abs{v_B}^2
   ] \abs{v_A - v_B}.
\end{equation}
Using the Sobolev embedding
$\norm{\cdot}_{\infty} \le C_3 \norm{\cdot}_{H^1}$
this immediately implies the first
estimate in \sref{eq:prlm:f:lip:bnds}.
Applying this estimate with $v_A = 0$
and $v_B = \Phi_0 - \Phi$
we find
\begin{equation}
\label{eq:prlm:f:apriori:f:phi}
\begin{array}{lcl}
\norm{f(\Phi)}_{L^2}
 & \le & \norm{f( \Phi_0)}_{L^2} + \norm{f(\Phi) - f(\Phi_0)}_{L^2}
\\[0.2cm]
 & \le & C_4.
\end{array}
\end{equation}
Exploiting
\begin{equation}
\norm{f(\Phi + v)}_{L^2}
\le \norm{f(\Phi)}_{L^2}
  + \norm{f(\Phi +v) - f(\Phi)}_{L^2},
\end{equation}
we hence
obtain
\begin{equation}
\begin{array}{lcl}
 \norm{f(\Phi + v)}_{L^2} & \le &
   C_5 \big[1 + \norm{v}_{L^2} + \norm{v}_{H^1}^2 \norm{v}_{L^2} \big] .
\\[0.2cm]
\end{array}
\end{equation}
The first estimate
in \sref{eq:prlm:f:glb:bnds}
now follows by noting
that
$\norm{v}_{L^2}
\le \norm{v}_{H^1}^2 \norm{v}_{L^2}$
for $\norm{v}_{L^2} \ge 1$.

Turning to the inner products,
\sref{eq:prlm:f:pointwise:lip}
allows us to compute
\begin{equation}
\label{eq:prlm:f:inn:product:a}
\begin{array}{lcl}
\abs{\langle
  f(\Phi + v_A) - f(\Phi + v_B)  , \psi_A
\rangle_{L^2} }
&
\le & C_2 \norm{v_A - v_B}_{L^2} \norm{\psi_A}_{L^2}
\\[0.2cm]
& & \qquad
+ C_2 \big[ \norm{v_A}_{L^2}^2 + \norm{v_B}_{L^2}^2 \big]
  \norm{v_A - v_B}_{H^1} \norm{\psi_A}_{H^1} .
\\[0.2cm]
\end{array}
\end{equation}
Exploiting \sref{eq:prlm:f:apriori:f:phi},
the second estimate
in \sref{eq:prlm:f:glb:bnds}
hence follows from the bound
\begin{equation}
\begin{array}{lcl}
\abs{ \langle f(\Phi + v) , \psi \rangle_{L^2} }
& \le&
\abs{ \langle f(\Phi), \psi \rangle_{L^2} }
+ \abs{ \langle f(\Phi + v) - f(\Phi),
   \psi \rangle_{L^2} },
\end{array}
\end{equation}
using a similar observation
as above to eliminate
the $\norm{v}_{L^2} \norm{\psi}_{L^2}$
term.
Finally, the second estimate
in \sref{eq:prlm:f:lip:bnds}
can be obtained by applying
\sref{eq:prlm:f:inn:product:a}
and
\sref{eq:prlm:f:glb:bnds}
to the splitting
\begin{equation}
\begin{array}{lcl}
\abs{ \langle f(\Phi + v_A), \psi_A \rangle_{L^2}
- \langle f(\Phi + v_B ) , \psi_B \rangle_{L^2} }
& \le &
  \abs{ \langle f(\Phi + v_A) - f(\Phi + v_B) , \psi_A \rangle_{L^2} }
\\[0.2cm]
& & \qquad
  + \abs{ \langle f(\Phi + v_B), \psi_A - \psi_B \rangle_{L^2} } .
\\[0.2cm]
\end{array}
\end{equation}
\end{proof}

\begin{lem}
\label{lem:prlm:ests:g}
Suppose that (Hg) and (hPar) are satisfied.
Then there exists a constant $K > 0$, which does not
depend on the pair $(\Phi, c)$, so that the following holds true.
For any $v \in H^1$
we have the bounds
\begin{equation}
 \label{eq:prlm:g:l2:bnds}
    \begin{array}{lcl}
      \norm{g(\Phi + v)}_{L^2}  & \le &
        \norm{g(\Phi_{0}) }_{L^2}
               + K_g (1 + \norm{v}_{L^2} )
\\[0.2cm]
& \le &
        K [1 + \norm{v}_{L^2}] ,
     \\[0.2cm]
     \norm{\partial_\xi g(\Phi + v)}_{L^2}  & \le &
        K [1 + \norm{v}_{H^1}] ,
    \end{array}
\end{equation}
while for any pair $(v_A, v_B ) \in H^1 \times H^1$
we have the estimates
\begin{equation}
  \label{eq:lem:prlm:g:lipschitz}
    \begin{array}{lcl}
      \norm{g(\Phi + v_A) - g(\Phi + v_B) }_{L^2} & \le &
        K \norm{v_A - v_B}_{L^2} ,
       \\[0.2cm]
       \norm{\partial_\xi [ g(\Phi + v_A) - g(\Phi + v_B) ]}_{L^2}
       & \le &
        K \big[ 1 + \norm{v_A}_{H^1}  \big] \norm{v_A - v_B}_{H^1} .
    \end{array}
\end{equation}
\end{lem}
\begin{proof}
The Lipschitz estimate on $g$ implies that
\begin{equation}
\norm{ g(\Phi + v_A) - g(\Phi + v_B) }_{L^2}
 \le K_g \norm{v_A - v_B}_{L^2}.
\end{equation}
Applying this inequality with $v_A = v$
and $v_B = \Phi_0 - \Phi$ we obtain
\begin{equation}
\begin{array}{lcl}
\norm{g( \Phi + v) }_{L^2} & \le &
 \norm{ g(\Phi_0 ) }_{L^2}
  + K_{g}
   \big[ \norm{\Phi - \Phi_{0} }_{L^2} + \norm{v}_{L^2}
       \big] ,
\\[0.2cm]
\end{array}
\end{equation}
which in view of (hPar)
yields the first line
of \sref{eq:prlm:g:l2:bnds}.

The uniform bound
\begin{equation}
\abs{ Dg(\Phi + v) } \le K_g
\end{equation}
together with the identity
\begin{equation}
\partial_{\xi} g(\Phi + v) =
  Dg( \Phi  + v) \big( \Phi' + v' \big)
\end{equation}
immediately imply
the second estimate
in \sref{eq:prlm:g:l2:bnds}.
Finally, using
\begin{equation}
\abs{ Dg(\Phi + v)  - Dg(\Phi + w)}
  \le K_{g} \abs{ v - w}
\end{equation}
and the identity
\begin{equation}
\begin{array}{lcl}
 \partial_\xi \big[ g(\Phi + v_A)  - g(\Phi + v_B) \big]
& = & \big[ Dg( \Phi  + v_A) - Dg(\Phi + v_B) \big]
   \big( \Phi' + v_A' \big)
\\[0.2cm]
& & \qquad + Dg(\Phi + v_B) \big( v_A' - v_B' \big) ,
\end{array}
\end{equation}
we obtain
\begin{equation}
\begin{array}{lcl}
\norm{
  \partial_\xi \big[ g(\Phi + v_A)  - g(\Phi + v_B) \big]
}_{L^2}
& \le & K_g \norm{v_A - v_B}_\infty
   \big[ \norm{\Phi'}_{L^2} + \norm{v_A'}_{L^2} \big]
\\[0.2cm]
& & \qquad + K_g \norm{ v_A' - v_B'}_{L^2}.
\end{array}
\end{equation}
The second estimate
in \sref{eq:lem:prlm:g:lipschitz}
now follows easily.
\end{proof}

\begin{lem}
\label{lem:prlm:g:bnds:ip}
Suppose that (Hg) and (hPar) are satisfied.
Then there exists a constant $K > 0$, which does not
depend on the pair $(\Phi, c)$, so that the following holds true.
For any $v \in H^1$ and $\psi \in H^1$ we have the bounds
\begin{equation}
 \label{eq:prlm:g:ests:g:ip}
    \begin{array}{lcl}
      \abs{\langle g(\Phi + v) , \psi \rangle_{L^2} }  & \le &
        K [1 + \norm{v}_{L^2}]\norm{\psi}_{L^2},
     \\[0.2cm]
     \abs{\langle \partial_\xi g(\Phi + v), \psi \rangle_{L^2} }  & \le &
        K [1 + \norm{v}_{L^2}]\norm{\psi}_{H^1},
    \end{array}
\end{equation}
while for any set of pairs $(v_A, v_B) \in H^1 \times H^1$
and $(\psi_A, \psi_B) \in H^1 \times H^1$
we have the estimates
\begin{equation}
\label{eq:prlm:g:diff:ests:g:ip}
    \begin{array}{lcl}
      \abs{ \langle g(\Phi + v_A) , \psi_A \rangle_{L^2}
      - \langle g(\Phi + v_B) , \psi_B \rangle_{L^2} } & \le &
        K \norm{v_A-v_B}_{L^2} \norm{\psi_A}_{L^2}
    \\[0.2cm]
    & & \qquad
        + K [1 + \norm{v_B}_{L^2}] \norm{\psi_A - \psi_B}_{L^2} ,
       \\[0.2cm]
       \abs{ \langle \partial_\xi [ g(\Phi + v_A) ] , \psi_A \rangle_{L^2}
         - \langle  \partial_\xi [ g(\Phi + v_B) ] , \psi_B \rangle_{L^2} }
       & \le &
        K \norm{v_A-v_B}_{L^2} \norm{\psi_A}_{H^1}
\\[0.2cm]
& & \qquad
        + K [1 + \norm{v_B}_{L^2}] \norm{\psi_A - \psi_B}_{H^1} .
    \end{array}
\end{equation}
\end{lem}
\begin{proof}
The estimates \sref{eq:prlm:g:ests:g:ip}
follow immediately
from
the bound $\norm{g(\Phi + v)}_{L^2}
\le K \big[ 1 + \norm{v}_{L^2} \big]$.
The first bound in
\sref{eq:prlm:g:diff:ests:g:ip}
can be obtained
from Lemma \ref{lem:prlm:ests:g}
by noting that
\begin{equation}
\begin{array}{lcl}
 \abs{  \langle g(\Phi + v_A) , \psi_A \rangle_{L^2}
      - \langle g(\Phi + v_B) , \psi_B \rangle_{L^2}
 }
 & \le &
     \abs{
        \langle g(\Phi + v_A)- g(\Phi + v_B) ,
          \psi_A \rangle_{L^2}
      }
\\[0.2cm]
& & \qquad
      + \abs{ \langle g(\Phi + v_B), \psi_A - \psi_B \rangle_{L^2} }.
\end{array}
\end{equation}
The final bound can be obtained
by transferring the derivative to the
functions $\psi_A$ and $\psi_B$.
\end{proof}

\subsection{Bounds on $\mathcal{J}_0$, $b$ and $\kappa_\sigma$}
\label{sec:prlm:rest}

We are now ready to obtain global and Lipschitz bounds
on the functions $\mathcal{J}_0$, $b$ and $\kappa_{\sigma}$.
In addition, we show that it suffices to impose an a priori bound
on $\norm{v}_{L^2}$ in order to avoid hitting the cut-offs in the definition of $b$.
This is crucial for the estimates in \S\ref{sec:nls}, where we have
uniform control on $\norm{v}_{L^2}$, but only
an integrated form of control on $\norm{v}_{H^1}$.

\begin{lem}
\label{lem:prlm:j0:bnds}
Suppose that (Hf) and (hPar) are satisfied.
Then there exists a constant $K > 0$, which does not
depend on the pair $(\Phi, c)$, so that the following holds true.
For any $(v, d) \in H^1 \times \Real$ and $\psi \in H^1$
we have the bounds\footnote{ We are dropping the third argument of $\mathcal{J}_0$ here since it is irrelevant when $\sigma = 0$.}
\begin{equation}
\label{eq:prlm:j0:glb:bnds}
\begin{array}{lcl}
 \norm{\mathcal{J}_0(\Phi + v, c + d)}_{L^2} & \le &
   K (1 + \abs{d} ) \big[
     1 + \norm{v}_{H^1} + \norm{v}^2_{H^1} \norm{v}_{L^2}
   \big] ,
\\[0.2cm]
 \abs{\langle \mathcal{J}_0(\Phi + v, c + d), \psi
    \rangle_{L^2} }
   & \le &
      K (1 + \abs{d} ) \big[
     1 +  \norm{v}_{H^1} \norm{v}_{L^2}^2
   \big]
       \norm{\psi}_{H^1} .
\\[0.2cm]
\end{array}
\end{equation}
In addition, for any set
of pairs
$(v_A, v_B) \in H^1 \times H^1$,
$(d_A, d_B) \in \Real \times \Real$
and $(\psi_A, \psi_B) \in H^1 \times H^1$,
the expressions
\begin{equation}
\begin{array}{lcl}
\Delta_{AB} \mathcal{J}_0
  & = &
    \mathcal{J}_0(\Phi + v_A, c +  d_A)
  - \mathcal{J}_0(\Phi + v_B, c + d_B) ,
\\[0.2cm]
\Delta_{AB} \langle \mathcal{J}_0, \cdot \rangle_{L^2}
  & = &
    \langle \mathcal{J}_0(\Phi + v_A, c +  d_A) , \psi_A \rangle_{L^2}
  - \langle \mathcal{J}_0(\Phi + v_B, c +  d_B) ,
    \psi_B \rangle_{L^2}
\end{array}
\end{equation}
satisfy the estimates
\begin{equation}
\label{eq:prlm:j0:lip:bnds}
\begin{array}{lcl}
\norm{\Delta_{AB} \mathcal{J}_0 }_{L^2}
  & \le &
  K [\norm{v_A}_{H^1} \norm{v_A}_{L^2}
    + \norm{v_B}_{H^1} \norm{v_B}_{L^2} ] \norm{v_A-v_B}_{H^1}
\\[0.2cm]
& & \qquad  + [1 + \norm{v_A}_{H^1} ]
  \abs{d_A - d_B }
\\[0.2cm]
& & \qquad
  + K( 1 + \abs{d_B} ) \norm{v_A-v_B}_{H^1} ,
\\[0.2cm]
\abs{ \Delta_{AB} \langle \mathcal{J}_0, \cdot \rangle_{L^2} }
 & \le &
 K [\norm{v_A}_{L^2}^2 +  \norm{v_B}_{L^2}^2 ] \norm{v_A-v_B}_{H^1}
   \norm{\psi_A}_{H^1}
\\[0.2cm]
& & \qquad  + [1 + \norm{v_A}_{L^2} ] \abs{d_A - d_B }
   \norm{\psi_A}_{H^1}
\\[0.2cm]
& & \qquad
  + K( 1 + \abs{d_B} ) \norm{v_A-v_B}_{L^2} \norm{\psi_A}_{H^1}
\\[0.2cm]
& & \qquad
  + K (1 + \abs{d_B})
  \big[1
     + \norm{v_B}_{H^1} \norm{v_B}_{L^2}^2 \big]
  \norm{\psi_A - \psi_B}_{H^1} .
\end{array}
\end{equation}
\end{lem}
\begin{proof}
We first note that
the terms in
\sref{eq:prlm:f:glb:bnds}-\sref{eq:prlm:f:lip:bnds}
can be absorbed
in
\sref{eq:prlm:j0:glb:bnds}-\sref{eq:prlm:j0:lip:bnds},
so it suffices
to study the function
\begin{equation}
\mathcal{J}_{0;II}(u, \overline{c})
 = \overline{c} u'.
\end{equation}
Recalling that (hPar) implies
\begin{equation}
\abs{c} + \norm{\Phi'}_{L^2} \le C_1,
\end{equation}
we find
\begin{equation}
\norm{
\mathcal{J}_{0;II}(\Phi + v, c + d)
}_{L^2} \le
C_2 ( 1 + \abs{d} )( 1 + \norm{v}_{H^1} ) ,
\end{equation}
together with
\begin{equation}
\abs{
\langle \mathcal{J}_{0;II}(\Phi + v, c + d), \psi \rangle_{L^2}
} \le C_2 ( 1 + \abs{d} )( 1 + \norm{v}_{L^2} )
  \norm{\psi}_{H^1},
\end{equation}
which can be absorbed in
\sref{eq:prlm:j0:glb:bnds}.

In addition, writing
\begin{equation}
\begin{array}{lcl}
\Delta_{AB} \mathcal{J}_{0;II}
  & = &
    \mathcal{J}_{0;II}(\Phi + v_A, c +  d_A)
  - \mathcal{J}_{0;II}(\Phi + v_B, c + d_B) ,
\\[0.2cm]
\Delta_{AB} \langle \mathcal{J}_{0;II}, \cdot \rangle_{L^2}
  & = &
    \langle \mathcal{J}_{0;II}(\Phi + v_A, c +  d_A) , \psi_A \rangle_{L^2}
  - \langle \mathcal{J}_{0;II}(\Phi + v_B, c +  d_B) ,  \psi_B \rangle_{L^2} ,
\end{array}
\end{equation}
we compute
\begin{equation}
\Delta_{AB} \mathcal{J}_{0;II}
= (d_A - d_B ) (\Phi' + v_A')
+ (c + d_B) (v_A' - v_B') .
\end{equation}
This yields
\begin{equation}
\norm{
  \Delta_{AB} \mathcal{J}_{0;II}
}_{L^2}
\le C_3 \abs{d_A - d_B} ( 1 + \norm{v_A}_{H^1} )
+ C_3 (1 + \abs{d_B} )\norm{v_A - v_B}_{H^1},
\end{equation}
which establishes the first estimate
in \sref{eq:prlm:j0:lip:bnds}.

In a similar fashion, we obtain
\begin{equation}
\abs{
 \langle
  \Delta_{AB} \mathcal{J}_{0;II},
  \psi \rangle_{L^2}
}
\le C_3 \abs{d_A - d_B} ( 1 + \norm{v_A}_{L^2} )
    \norm{\psi}_{H^1}
+ C_3 (1 + \abs{d_B} )\norm{v_A - v_B}_{L^2}
    \norm{\psi}_{H^1}.
\end{equation}
The remaining estimate
now follows from the inequality
\begin{equation}
\begin{array}{lcl}
\abs{
  \Delta_{AB} \langle \mathcal{J}_{0;II}, \cdot \rangle_{L^2}
}
& \le &
\abs{\langle \Delta_{AB} \mathcal{J}_{0;II}
 , \psi_A \rangle_{L^2} }
\\[0.2cm]
& & \qquad
  + \abs{
    \langle \mathcal{J}_{0;II}(\Phi + v_B, c + d_B ),
       \psi_A - \psi_B
    \rangle_{L^2} } .
\end{array}
\end{equation}
\end{proof}

\begin{lem}
\label{lem:prlm:bnds:deriv:u:vs:psi}
Assume that (hPar) is satisfied.
Then there exists a constant $K > 0$, which does not
depend on the pair $(\Phi, c)$, so that the following holds true.
For any $v \in H^1$
and $\psi \in H^1$
we have the bound
\begin{equation}
\abs{\langle \partial_\xi (\Phi + v) , \psi \rangle_{L^2} }
\le K \big[ 1 + \norm{v}_{L^2} \big] \norm{\psi}_{H^1},
\end{equation}
while for any set of pairs
$(v_A, v_B) \in H^1 \times H^1$
and $(\psi_A, \psi_B) \in H^1 \times H^1$
we have the estimate
\begin{equation}
\begin{array}{lcl}
\abs{\langle \partial_\xi [\Phi + v_A] , \psi_A \rangle_{L^2}
 - \langle \partial_\xi [\Phi + v_B] , \psi_B \rangle_{L^2}
}
& \le &  K\norm{v_A - v_B}_{L^2 } \norm{\psi_A}_{H^1}
\\[0.2cm]
& & \qquad
  + K\big[1 + \norm{v_B}_{L^2} \big] \norm{\psi_A-\psi_B}_{H^1}.
\end{array}
\end{equation}
\end{lem}
\begin{proof}
The desired bounds follow
from the identity
\begin{equation}
\abs{\langle \partial_\xi (\Phi + v) , \psi \rangle_{L^2} }
= \abs{\langle \Phi + v , \partial_\xi \psi \rangle_{L^2} } ,
\end{equation}
together with
the estimate
\begin{equation}
\begin{array}{lcl}
\abs{\langle \partial_\xi [\Phi + v_A] , \psi_A \rangle_{L^2}
 - \langle \partial_\xi [\Phi + v_B] , \psi_B \rangle_{L^2}
}
& \le & \abs{\langle v_A - v_B , \partial_\xi \psi_A \rangle_{L^2} }
+ \abs{ \langle \partial_\xi \Phi , \psi_A - \psi_B \rangle_{L^2} }
\\[0.2cm]
& & \qquad
+ \abs{ \langle  v_B , \partial_{\xi} [\psi_A - \psi_B]
  \rangle_{L^2} } .
\\[0.2cm]
\end{array}
\end{equation}
\end{proof}

\begin{lem}
\label{lem:prlm:bnds:b}
Suppose that (Hg)
and (hPar) are satisfied.
Then there exist constants $K_b > 0$
and $K > 0$,
which do not
depend on the pair $(\Phi, c)$, so that the following holds true.
For any $v \in H^1$ and $\psi \in H^1$
we have the bound
\begin{equation}
  \label{eq:prlm:b:glb:bnd}
  \begin{array}{lcl}
      \abs{b(\Phi + v, \psi)}  & \le &
        K_b,
     \\[0.2cm]
  \end{array}
\end{equation}
while for any set of pairs
$(v_A , v_B) \in H^1 \times H^1$
and $(\psi_A, \psi_B) \in H^1 \times H^1$
we have the estimate
\begin{equation}
\label{eq:prlm:b:lipschitz}
    \begin{array}{lcl}
      \abs{b(\Phi + v_A , \psi_A) - b(\Phi + v_B, \psi_B) }
        & \le &
           K \norm{v_A-v_B}_{L^2}
             \norm{\psi_A}_{H^1}
   \\[0.2cm]
   & & \qquad
           + K \big[ 1 +  \norm{v_B}_{L^2} \big]
              \norm{\psi_A - \psi_B}_{H^1}.
    \end{array}
\end{equation}
\end{lem}
\begin{proof}
The uniform bound \sref{eq:prlm:b:glb:bnd}
follows directly from the properties of the cut-off functions.
Upon introducing the function
\begin{equation}
\tilde{b}(x,y)
= - \chi_{\mathrm{low}}( x)^{-1} \chi_{\mathrm{high}}(y),
\end{equation}
the global Lipschitz smoothness of the cut-off functions
implies that
\begin{equation}
\abs{
  \tilde{b}(x_A ,y_A)
 - \tilde{b}(x_B , y_B)
}
\le C_1 \big[ \abs{x_B - x_A}
+ \abs{y_B - y_A} \big].
\end{equation}
Using the identity
\begin{equation}
b(u, \psi)
= \tilde{b}\Big( \langle \partial_\xi u , \psi \rangle_{L^2} ,
    \langle g(u), \psi \rangle_{L^2} \Big),
\end{equation}
the desired bound \sref{eq:prlm:b:lipschitz}
follows from Lemma's
\ref{lem:prlm:g:bnds:ip}
and
\ref{lem:prlm:bnds:deriv:u:vs:psi}.
\end{proof}

\begin{lem}
\label{lem:prlm:b:without:cutoffs}
Assume that (Hg) and (hPar) are satisfied.
Then for any $v \in H^1$ that has
\begin{equation}
  \norm{v}_{L^2} \le \min \{1 , [4 \norm{\psi_{\mathrm{tw}}}_{H^1} ]^{-1} \},
\end{equation}
we have the identity
\begin{equation}
\label{eq:prlm:id:for:b:without:cutoffs}
b(\Phi + v, \psi_{\mathrm{tw}} ) =
- \big[\langle \partial_\xi [\Phi + v] , \psi_{\mathrm{tw}} \rangle_{L^2} \big]^{-1}
\langle g(\Phi + v) , \psi_{\mathrm{tw}} \rangle_{L^2}.
\end{equation}
\end{lem}
\begin{proof}
Using \sref{eq:prlm:g:l2:bnds}
and recalling
the definition
\sref{eq:mr:def:k:ip},
we find that
\begin{equation}
\abs{\langle g(\Phi + v), \psi_{\mathrm{tw}} \rangle_{L^2} }
\le
\Big[ \norm{g(\Phi_0)}_{L^2} + 2 K_g \Big]
  \norm{\psi_{\mathrm{tw}}}_{L^2}
= K_{\mathrm{ip}}.
\end{equation}
In addition, we note that
(hPar) and the normalisation
\sref{eq:mr:hs:norm:cnd:psitw}
imply that
\begin{equation}
 \langle \partial_\xi \Phi  , \psi_{\mathrm{tw}} \rangle_{L^2}
  =  \langle \partial_\xi \Phi_{0} ,  \psi_{\mathrm{tw}}
         \rangle_{L^2}
  + \langle \partial_\xi [ \Phi - \Phi_{0} ],
        \psi_{\mathrm{tw}} \rangle_{L^2}
  \ge 1 - \frac{1}{4} = \frac{3}{4}.
\end{equation}
This allows us to estimate
\begin{equation}
\langle \partial_\xi (\Phi + v) , \psi_{\mathrm{tw}}
 \rangle_{L^2}
\ge
\frac{3}{4} -
\langle v, \partial_\xi \psi_{\mathrm{tw}} \rangle_{L^2}
 \ge
\frac{3}{4} - \frac{1}{4}
= \frac{1}{2},
\end{equation}
which shows that the cut-off functions
do not modify their arguments.
\end{proof}

\begin{lem}
\label{lem:prlm:bnds:kappa}
Suppose that $(Hg)$ and $(hPar)$
are satisfied.
Then there exists a constant $K_{\kappa} > 0$, which does not
depend on the pair $(\Phi, c)$, so that 
for any
$0 \le \sigma \le 1$,
any
$v \in H^1$
and any $\psi\in H^1$,
we have the bound
\begin{equation}
\label{eq:prlm:kappa:glb:ests}
\abs{ \kappa_{\sigma}(\Phi + v, \psi) }
+ \abs{ \kappa_{\sigma}(\Phi + v, \psi)^{-1} }
+ \abs{ \kappa_{\sigma}(\Phi + v, \psi)^{-1/2} }
  \le K_{\kappa} .
\end{equation}
\end{lem}
\begin{proof}
This follows directly from the bound
\begin{equation}
1 \le \kappa_{\sigma}(\Phi + v, \psi) \le
  1 + \frac{1}{2\rho} \sigma^2 K_b^2.
\end{equation}
\end{proof}

In order to state our final result,
we introduce the functions
\begin{equation}
\label{eq:prlm:defs:nu}
\begin{array}{lcl}
\nu^{(1)}_{\sigma}(u, \psi) & = & \kappa_{\sigma}(u, \psi) - 1 ,
\\[0.2cm]
\nu^{(-1)}_{\sigma}(u, \psi) & =  &
    \kappa_{\sigma}(u, \psi)^{-1} - 1 ,
\\[0.2cm]
\nu^{(-1/2)}_{\sigma}(u, \psi) & =  &
    \kappa_{\sigma}(u, \psi)^{-1/2} - 1,
\end{array}
\end{equation}
which isolate the $\sigma$-dependence in $\kappa_{\sigma}$.

\begin{lem}
\label{lem:prlm:nu}
Suppose that $(Hg)$ and $(hPar)$
are satisfied and pick
$\vartheta \in \{-1, -\frac{1}{2} ,1 \}$.
Then there exist constants $K_{\nu} > 0$ and $K > 0$, which do not
depend on the pair $(\Phi, c)$, so that the following holds true.
For any
$0 \le \sigma \le 1$,
any
$v \in H^1$
and any $\psi\in H^1$
we have the bound
\begin{equation}
\label{eq:prlm:nu:abs:bnd}
\abs{ \nu_{\sigma}^{(\vartheta)}(\Phi + v, \psi) }
  \le \sigma^2 K_\nu ,
\end{equation}
while for any
$0 \le \sigma \le 1$
and any
set of pairs
$(v_A, v_B) \in H^1 \times H^1$
and $(\psi_A, \psi_B) \in H^1 \times H^1$,
we have the estimate
\begin{equation}
\label{eq:prlm:lip:bnds:nu}
\begin{array}{lcl}
\abs{
  \nu^{(\vartheta)}_{\sigma} ( \Phi + v_A, \psi_A)
    - \nu^{(\vartheta)}_{\sigma} ( \Phi + v_B, \psi_B)
} & \le &
  K \sigma^2
   \norm{v_A-v_B}_{L^2} \norm{\psi_A}_{H^1}
\\[0.2cm]
& & \qquad
 + K \sigma^2 \big[ 1 +  \norm{v_B}_{L^2} \big]
   \norm{\psi_A -\psi_B}_{H^1}.
\end{array}
\end{equation}
\end{lem}
\begin{proof}
As a preparation, we observe that
for any $x \ge 0$ and $y \ge 0$
we have the inequality
\begin{equation}
\abs{\frac{1}{1 + x} - \frac{1}{1+y}}
 = \frac{\abs{y - x} } {(1 + x)(1+y)}
 \le \abs{y-x} ,
\end{equation}
together with
\begin{equation}
\abs{\frac{1}{\sqrt{1 + x}}
  - \frac{1}{\sqrt{1+y}}}
 = \frac{\abs{y - x} }
   {\sqrt{(1 + x)(1+y)}(\sqrt{1 + x} + \sqrt{1 + y})}
 \le \frac{1}{2}\abs{y-x}.
\end{equation}
Applying these bounds with $y = 0$,
we obtain
\begin{equation}
\abs{\nu^{(\vartheta)}_{\sigma}(\Phi + v, \psi) }
 \le \frac{1}{2 \rho} \sigma^2 \abs{ b(\Phi + v, \psi)}^2
\le \frac{1}{2 \rho} \sigma^2 K_b^2,
\end{equation}
which yields \sref{eq:prlm:nu:abs:bnd}.
In addition,
we may compute
\begin{equation}
\begin{array}{lcl}
\abs{
  \nu^{(\vartheta)}_{\sigma} ( \Phi + v_A, \psi_A)
    - \nu^{(\vartheta)}_{\sigma} ( \Phi + v_B, \psi_B)
}
& \le &
  \frac{1}{2 \rho} \sigma^2
\abs{
  b(\Phi + v_A, \psi_A)^2 - b(\Phi + v_B, \psi_B)^2
}
\\[0.2cm]
& = & \frac{1}{2 \rho} \sigma^2
   \abs{b(\Phi + v_A, \psi_A) + b(\Phi + v_B, \psi_B) }
\\[0.2cm]
& & \qquad \qquad \times
       \abs{b(\Phi + v_A, \psi_A) - b(\Phi + v_B, \psi_B) }.
\\[0.2cm]
\end{array}
\end{equation}
In particular,
the bounds \sref{eq:prlm:lip:bnds:nu}
follow from Lemma \ref{lem:prlm:bnds:b}.
\end{proof}

\section{Variational solution}
\label{sec:var}

In this section we set out to establish
Proposition \ref{prp:mr:main:ex}.
Our strategy is to fit the
first component of \sref{eq:mr:formal:spde}
into the framework of \cite{LiuRockner}.
Indeed, the conditions (H1)-(H4) in \cite{LiuRockner} 
are explicitly verified in Lemma \ref{lem:var:cond:from:liu:rockner} below.
The second line of \sref{eq:mr:formal:spde} can subsequently be treated
as an SDE for $\Gamma$ with random coefficients.
In Lemma \ref{lem:var:conds:thm:rock:prevot} below we
show that this SDE fits into the framework that was developed in
\cite[Chapter 3]{Concise} to handle such equations.

\begin{lem}
\label{lem:var:cond:from:liu:rockner}
Suppose that (HA), (Hf), (HTw), (HS), (HVar)
and (Hg) are all satisfied.
Then there exist constants $K > 0$
and $\vartheta > 0$ so that
the following properties hold true.
\begin{itemize}
\item[(i)]{
  For any triplet $(v_A, v_B, v) \in  H^1 \times H^1 \times H^1$,
  the map
  \begin{equation}
      s \mapsto \langle A_* [v_A + s v_B ] , v \rangle_{H^{-1};H^1}
      + \langle f( \Phi_{\mathrm{ref}} + v_A + s v_B ) , v \rangle_{L^2}
  \end{equation}
  is continuous.
}
\item[(ii)]{
 For every pair $(v_A, v_B) \in H^1 \times H^1$,
 we have the inequality
 \begin{equation}
 \begin{array}{lcl}
   K \norm{v_A - v_B}_{L^2}^2
    & \ge &
     2 \langle A_* (v_A - v_B) , v_A - v_B \rangle_{H^{-1};H^1}
     \\[0.2cm]
     & & \qquad
     +2 \langle f(\Phi_{\mathrm{ref}} + v_A) - f(\Phi_{\mathrm{ref}} + v_B) ,
        v_A - v_B \rangle_{L^2}
     \\[0.2cm]
     & & \qquad
     + \norm{g(\Phi_{\mathrm{ref}}+ v_A) - g(\Phi_{\mathrm{ref}} + v_B)}_{L^2}^2 .
    \end{array}
 \end{equation}
}
\item[(iii)]{
  For any $v \in H^1$ we have
  the inequality
  \begin{equation}
     2 \langle A_* v, v \rangle_{H^{-1} ; H^1 }
     + 2 \langle f(\Phi_{\mathrm{ref}} + v), v \rangle_{L^2}
     + \norm{g(\Phi_{\mathrm{ref}} + v)}_{L^2}^2
     + \vartheta \norm{v}_{H^1}^2
     \le K \big[ 1 + \norm{v}_{L^2}^2 \big].
  \end{equation}
}
\item[(iv)]{
  For any $v \in H^1$ we have the bound
  \begin{equation}
     \norm{A_* v}_{H^{-1}}^2
     + \norm{ f(\Phi_{\mathrm{ref}} + v) }_{H^{-1}}^2
     \le K \big[ 1 + \norm{v}_{H^1}^2 \big]
      \big[ 1 + \norm{v}_{L^2}^4 \big] .
  \end{equation}
}
\end{itemize}
\end{lem}
\begin{proof}
Item (i) follows from the linearity of $A_*$ and
the Lipschitz bound \sref{eq:prlm:f:lip:bnds}.
In addition,
writing
\begin{equation}
\mathcal{I} =
\langle f(\Phi_{\mathrm{ref}} + v_A)
  - f(\Phi_{\mathrm{ref}} + v_B), v_A - v_B \rangle_{L^2},
\end{equation}
(HVar) implies the one-sided inequality
\begin{equation}
\begin{array}{lcl}
\mathcal{I}
&  = &
\langle f(\Phi_{\mathrm{ref}} + v_A) - f(\Phi_{\mathrm{ref}} + v_B),
     \Phi_{\mathrm{ref}} + v_A - (\Phi_{\mathrm{ref}} +v_B) \rangle_{L^2}
\\[0.2cm]
& \le & C_1 \norm{v_A - v_B}_{L^2}^2.
\end{array}
\end{equation}
Item (ii) hence follows from
the Lipschitz bound
\sref{eq:lem:prlm:g:lipschitz}
together with the bound
\begin{equation}
\label{eq:var:azero:ineq}
\langle A_* v, v \rangle_{H^{-1}, H^1}
\le - \rho \norm{v}_{H^1}^2.
\end{equation}
A second consequence of (HVar)
is that
\begin{equation}
\begin{array}{lcl}
\langle f(\Phi_{\mathrm{ref}} + v), v \rangle_{L^2}
& = &
  \langle f(\Phi_{\mathrm{ref}} + v) - f(\Phi_{\mathrm{ref}}) , (\Phi_{\mathrm{ref}} + v) - \Phi_{\mathrm{ref}} \rangle_{L^2}
\\[0.2cm]
& & \qquad
  + \langle f(\Phi_{\mathrm{ref}}) , v \rangle_{L^2}
\\[0.2cm]
& \le &
  C_1 \norm{v}_{L^2}^2 + \norm{f(\Phi_{\mathrm{ref}})}_{L^2} \norm{v}_{L^2}
\\[0.2cm]
& \le &
  C_2 \big[ 1 + \norm{v}_{L^2}^2 \big] .
\end{array}
\end{equation}
In particular, we may obtain
(iii) by combining
\sref{eq:var:azero:ineq}
with \sref{eq:prlm:g:l2:bnds}.

Finally,
for any $v \in H^1$ and $\psi \in H^1$
we may use
\sref{eq:prlm:f:glb:bnds}
to compute
\begin{equation}
\begin{array}{lcl}
\langle f(\Phi_{\mathrm{ref}} + v) ,
   \psi \rangle_{H^{-1} ; H^1}
& = &
  \langle f(\Phi_{\mathrm{ref}} + v) ,
   \psi \rangle_{ L^2}
\\[0.2cm]
& \le &
 C_3
 \Big[ 1
   + \norm{v}_{H^1} \norm{v}_{L^2}^2
 \Big] \norm{\psi}_{H^1} .
\\[0.2cm]
\end{array}
\end{equation}
In other words, we see that
\begin{equation}
\norm{f(\Phi_{\mathrm{ref}} + v)}_{H^{-1}} \le
C_3
\big[ 1 + \norm{v}_{H^1} \norm{v}_{L^2}^2  \big]
\le C_3 ( 1 + \norm{v}_{H^1} )
  ( 1 + \norm{v}_{L^2}^2  ),
\end{equation}
which yields (iv).
\end{proof}

\begin{lem}
\label{lem:var:bnds:a}
Suppose that (HA), (Hf),
(Hg) and (hPar) are all satisfied.
Then there exists a constant $K > 0$, which does not
depend on the pair $(\Phi, c)$, so that the following
properties hold true for any $0 \le \sigma \le 1$.
\begin{itemize}
\item[(i)]{
For any $v \in H^1$
and any $\psi \in H^2$ with
$\norm{\psi}_{H^2} \le 2 \norm{\psi_{\mathrm{tw}}}_{H^2}$,
we have the bound
\begin{equation}
\abs{a_{\sigma}(\Phi + v , c,\psi) }
\le K \Big[1
   + \norm{v}_{H^1} \norm{v}_{L^2}^2
     \Big]  .
\end{equation}
}
\item[(ii)]{For any $v \in H^1$
and any pair  $(\psi_A, \psi_B) \in H^2 \times H^2$
for which $\norm{\psi_A}_{H^2} \le 2 \norm{\psi_{\mathrm{tw}}}_{H^2}$
and $\norm{\psi_B}_{H^2} \le 2 \norm{\psi_{\mathrm{tw}}}_{H^2}$,
the difference
\begin{equation}
  \Delta_{AB} a_{\sigma} =
  a_{\sigma}(\Phi + v , c,\psi_A )
  - a_{\sigma}(\Phi + v , c,\psi_B)
\end{equation}
satisfies the bound
\begin{equation}
\abs{ \Delta_{AB} a_{\sigma} }
 \le
K \Big[ 1  + \norm{v}_{H^1} (1 + \norm{v}_{L^2}^3) \Big]
  \norm{\psi_A -\psi_B}_{H^1} .
\end{equation}
}
\end{itemize}
\end{lem}
\begin{proof}
We first compute
\begin{equation}
\begin{array}{lcl}
\kappa_{\sigma}(u,\psi)
\langle \mathcal{J}_{\sigma}(u, c, \psi) ,
   \psi \rangle_{L^2}
& = &
  \langle
    f(u )
    +c u'
    +\sigma^2 b( u, \psi)
      \partial_\xi [g(u )]
   , \psi \rangle_{L^2}
\\[0.2cm]
& = &
  \langle \mathcal{J}_0(u , c) , \psi \rangle_{L^2}
    + \sigma^2 b(u, \psi)
    \langle \partial_\xi [g(u)], \psi \rangle_{L^2}.
\end{array}
\end{equation}
Upon defining
\begin{equation}
\begin{array}{lcl}
  \mathcal{E}_{I}(u,c,\psi) & = &
    \langle \mathcal{J}_0(u, c), \psi \rangle_{L^2},
\\[0.2cm]
  \mathcal{E}_{II}(u,\psi) & = &
      \sigma^2 b(u, \psi)
      \langle \partial_\xi g(u), \psi \rangle_{L^2},
\\[0.2cm]
  \mathcal{E}_{III}(u,\psi) & = &
    \kappa_{\sigma}(u, \psi) \langle u, A_* \psi \rangle_{L^2},
\end{array}
\end{equation}
we hence see that
\begin{equation}
a_\sigma(u, c, \psi) =
- \Big[ \chi_{\mathrm{low}}\big(\langle \partial_\xi u, \psi \rangle_{L^2} \big) \Big]^{-1}
\big[ \mathcal{E}_{I}(u,c,\psi)
  + \mathcal{E}_{II}(u,\psi)
  + \mathcal{E}_{III}(u,\psi)\big].
\end{equation}
For $\# \in \{I, II, III\}$,
we define
\begin{equation}
\Delta_{AB} \mathcal{E}_{\#}
= \mathcal{E}_{\#}(\Phi + v , c, \psi_A)
 - \mathcal{E}_{\#}(\Phi + v , c, \psi_B) .
\end{equation}
We note that
Lemma's
\ref{lem:prlm:g:bnds:ip},
\ref{lem:prlm:j0:bnds}
and \ref{lem:prlm:bnds:b}
yield the bounds
\begin{equation}
\begin{array}{lcl}
\abs{\mathcal{E}_{I}(\Phi +v , c, \psi )}
  & \le & C_1
  \big[1
     + \norm{v}_{H^1} \norm{v}_{L^2}^2 \big] ,
\\[0.2cm]
\abs{\mathcal{E}_{II}(\Phi + v, \psi)}
 & \le & C_1
    \big[1 + \norm{v}_{L^2} \big] ,
\\[0.2cm]
\end{array}
\end{equation}
together with
\begin{equation}
\begin{array}{lcl}
\abs{\Delta_{AB} \mathcal{E}_{I}}
  & \le & C_1 \big[ 1  + \norm{v}_{H^1} \norm{v}_{L^2}^2 \big]
      \norm{\psi_A - \psi_B}_{H^1} ,
\\[0.2cm]
\abs{\Delta_{AB} \mathcal{E}_{II}}
 & \le &   C_1  [ 1+ \norm{v}_{L^2}]^2
  \norm{\psi_A - \psi_B}_{H^1}
\\[0.2cm]
& & \qquad
  + C_1  [1 + \norm{v}_{L^2} ]
    \norm{\psi_A - \psi_B}_{H^1} .
\\[0.2cm]
\end{array}
\end{equation}

A direct estimate using
the a-priori bound on $\norm{\psi}_{H^2}$
and \sref{eq:prlm:kappa:glb:ests}
yields
\begin{equation}
\begin{array}{lcl}
\abs{\mathcal{E}_{III}(\Phi + v, \psi)}
 & \le & K_{\kappa}
 \big[
    \abs{ \langle \Phi, A_* \psi \rangle_{L^2} }
    + \abs{ \langle v , A_* \psi \rangle_{L^2} }
 \big]
\\[0.2cm]
& \le &
 C_2 \big[1 + \norm{v}_{L^2} \big].
\end{array}
\end{equation}
By transferring one of the derivatives in $A_*$,
we also obtain using Lemma \ref{lem:prlm:nu} the bound
\begin{equation}
\begin{array}{lcl}
\abs{\Delta\mathcal{E}_{III}}
 & \le &
  \abs{\kappa_{\sigma}(\Phi + v , \psi_A)
  - \kappa_{\sigma}(\Phi + v, \psi_B) }
    \abs{ \langle \Phi + v , A_* \psi_A \rangle_{L^2} }
\\[0.2cm]
& & \qquad
+ \abs{\kappa_{\sigma}(\Phi + v, \psi_B) }
   \abs{ \langle \Phi + v ,
      A_* [\psi_A - \psi_B] \rangle_{L^2}
   }
\\[0.2cm]
& \le &
 C_3 (1 + \norm{v}_{L^2})^2 \norm{\psi_A - \psi_B}_{H^1}
\\[0.2cm]
& & \qquad
    + C_3 \big[ 1 + \norm{v}_{H^1} \big]
      \norm{\psi_A - \psi_B}_{H^1} .
\end{array}
\end{equation}

Upon
writing
\begin{equation}
\begin{array}{lcl}
\mathcal{E}(u, c, \psi)
 & = &
 \mathcal{E}_{I}(u, c, \psi)
 +\mathcal{E}_{II}(u,  \psi)
 +\mathcal{E}_{II}(u,  \psi) ,
\\[0.2cm]
\Delta_{AB} \mathcal{E}
& = &
\Delta_{AB} \mathcal{E}_{I}
+ \Delta_{AB} \mathcal{E}_{II}
+ \Delta_{AB} \mathcal{E}_{III},
\end{array}
\end{equation}
we hence conclude that
\begin{equation}
\begin{array}{lcl}
\abs{\mathcal{E}(\Phi +v, c , \psi) }
& \le & C_4 \big[
1 + \norm{v}_{H^1} \norm{v}_{L^2}^2
\big] ,
\\[0.2cm]
\abs{\Delta_{AB} \mathcal{E}}
& \le &
C_4
\big[
1 + \norm{v}_{H^1} \big]
\big[ 1 + \norm{v}_{L^2}^2 \big]
\norm{\psi_A - \psi_B}_{H^1} .
\end{array}
\end{equation}
Item (i) follows
immediately from the first bound,
since $\chi_{\mathrm{low}}(\cdot)^{-1}$
is globally bounded.
To obtain (ii),
we compute
\begin{equation}
\begin{array}{lcl}
 \abs{\Delta_{AB} a_{\sigma}}
 & \le &
   C_5 \abs{ \langle \partial_\xi (\Phi + v) , \psi_A \rangle_{L^2}
     - \langle \partial_\xi (\Phi + v) , \psi_B \rangle_{L^2} }
     \abs{\mathcal{E}(\Phi + v, \psi_A) }
\\[0.2cm]
& & \qquad
 +  C_5 \abs{\Delta_{AB} \mathcal{E} }
\\[0.2cm]
& \le &
  C_6 [1 + \norm{v}_{L^2} ]
   \big[
1 + \norm{v}_{H^1} \norm{v}_{L^2}^2
\big] \norm{\psi_A - \psi_B}_{H^1}
\\[0.2cm]
& & \qquad
  + C_6 \big[
1 + \norm{v}_{H^1} \big]
\big[ 1 + \norm{v}_{L^2}^2 \big]
\norm{\psi_A - \psi_B}_{H^1}
\\[0.2cm]
& \le &
C_7 \Big[ 1 +  \norm{v}_{H^1} (1 + \norm{v}_{L^2}^3) \Big]
  \norm{\psi_A -\psi_B}_{H^1},
\end{array}
\end{equation}
in which we used several estimates of the form
\begin{equation}
\norm{v}_{L^2} \le C_8 \big[ 1 + \norm{v}_{L^2}^4 \big]
\le C_8 \big[ 1 + \norm{v}_{H^1} \norm{v}_{L^2}^3 \big].
\end{equation}
\end{proof}

Upon introducing the shorthands
\begin{equation}
\begin{array}{lcl}
  p(v, \gamma) & = &
    c + a_{\sigma}( \Phi_{\mathrm{ref}} + v , c, T_{\gamma} \psi_{\mathrm{tw}} ) ,
\\[0.2cm]
  q(v, \gamma) & = & b( \Phi_{\mathrm{ref}} + v,  T_{\gamma} \psi_{\mathrm{tw}} ) ,
\\[0.2cm]
\end{array}
\end{equation}
the second line of
\sref{eq:mr:formal:spde}
can be written as
\begin{equation}
\label{eq:var:sode:gamma}
d \Gamma = p\big(X(t), \Gamma(t)\big) dt
   + \sigma q\big(X(t), \Gamma(t)\big) d \beta_t.
\end{equation}
Taking the view-point that $X(t) = X(t, \omega)$
is known upon picking $\omega \in \Omega$,
\sref{eq:var:sode:gamma} can be viewed as an SDE
with random coefficients. Our next result
relates directly to the conditions of
\cite[Thm 3.1.1]{Concise}, which is specially
tailored for equations of this type.

\begin{lem}
\label{lem:var:conds:thm:rock:prevot}
Suppose that (HA), (Hf), (HTw), (HS)
and (Hg) are all satisfied and fix
$c \in \Real$ together with
$0 \le \sigma \le 1$.
Then there exists $K > 0$ so that
the following properties are satisfied.
\begin{itemize}
\item[(i)]{
For any $v \in H^1$
and any pair $(\gamma_A, \gamma_B) \in \Real^2$,
we have the inequality
\begin{equation}
\begin{array}{lcl}
K \big[ 1 +  \norm{v}_{H^1}(1 +  \norm{v}^3_{L^2} ) \big]
     \abs{ \gamma_A - \gamma_B }^2
& \ge &
 2 \big[ \gamma_A - \gamma_B \big] \big[ p( v, \gamma_A) - p(v, \gamma_B) \big]
\\[0.2cm]
& & \qquad
   + \abs{q(v,  \gamma_A) - q(v, \gamma_B) }^2 .
\end{array}
\end{equation}
}
\item[(ii)]{
For any $v \in H^1$ and $\gamma \in \Real$,
we have the inequality
\begin{equation}
 2  \gamma p(v,\gamma)
   + \abs{q(v, \gamma)}^2
   \le K \big[ 1  + \norm{v}_{H^1} \norm{v}^2_{L^2} \big]
   \big[ 1 + \gamma^2 \big] .
\end{equation}
}
\item[(iii)]{
For any $v \in H^1$ and $\gamma \in \Real$,
we have the bound
\begin{equation}
  \abs{p(v,\gamma)} + \abs{q(v , \gamma)}^2
  \le K \big[ 1  + \norm{v}_{H^1} \norm{v}^2_{L^2} \big] .
\end{equation}
}
\end{itemize}
\end{lem}
\begin{proof}
The exponential decay of $\psi'_{\mathrm{tw}}$ and $\psi''_{\mathrm{tw}}$ implies that
\begin{equation}
\label{eq:var:diff:psi:tw}
\norm{T_{\gamma_A} \psi_{\mathrm{tw}} - T_{\gamma_B} \psi_{\mathrm{tw}} }_{H^1}
\le C_1 \abs{\gamma_A - \gamma_B}.
\end{equation}
Using
Lemma's \ref{lem:prlm:bnds:b} and \ref{lem:var:bnds:a},
we hence find the bounds
\begin{equation}
\begin{array}{lcl}
\abs{p(v, \gamma)}
& \le &
C_2 \big[ 1 +  \norm{v}_{H^1} \norm{v}^2_{L^2} \big] ,
\\[0.2cm]
\abs{q(v, \gamma)}
& \le & K_b,
\end{array}
\end{equation}
together with
\begin{equation}
\begin{array}{lcl}
\abs{p(v, \gamma_A) - p(v , \gamma_B)}
& \le &
C_3 \big[ 1  + \norm{v}_{H^1}(1 +  \norm{v}^3_{L^2} ) \big]
\abs{ \gamma_A - \gamma_B } ,
\\[0.2cm]
\abs{ q(v, \gamma_A) - q(v, \gamma_B) }
& \le & C_3 [1 + \norm{v}_{L^2}] \abs{\gamma_A - \gamma_B } .
\end{array}
\end{equation}
Items (i), (ii) and (iii) can now be verified directly.
\end{proof}

\begin{proof}[Proof of Proposition \ref{prp:mr:main:ex}]
The existence of the $dt \otimes \mathbb{P}$
version of $X$ that
is $(\mathcal{F}_t)$-progressively
measurable as a map into $H^1$,
follows from \cite[Ex. 4.2.3]{Concise}.

We remark that the
conditions (H1) through (H4) appearing
in \cite{LiuRockner} correspond
directly with items (i)-(iv) of Lemma \ref{lem:var:cond:from:liu:rockner}.
In particular, we may apply the main result
from this paper with $\alpha = 2$ and $\beta = 4$
to verify the remaining statements concerning $X$.

Finally, we note that items (i)-(iii) of
Lemma \ref{lem:var:conds:thm:rock:prevot}
allow us to apply
\cite[Thm. 3.1.1]{Concise},
provided that the function
\begin{equation}
t \mapsto \big[ 1  + \norm{X(t)}_{H^1}(1 +  \norm{X(t)}^3_{L^2} )\big]
\end{equation}
is integrable on $[0,T]$ for almost all $\omega \in \Omega$. This however
follows directly from the inclusions
\begin{equation}
X \in L^6\big(\Omega, \mathbb{P} ; C([0,T];L^2) \big)
\cap \mathcal{N}^2\big( [0,T]; (\mathcal{F}_t); H^1 \big),
\end{equation}
allowing us to verify the statements concerning $\Gamma$.
The remaining inclusions \sref{eq:mr:main:ex:inclusions:g:b}
follow directly from
the bounds in Lemma \ref{lem:prlm:ests:g}
and \ref{lem:prlm:bnds:b}.
\end{proof}

\section{The stochastic phase-shift}
\label{sec:sps}

In this section we consider the process $(X, \Gamma)$
described in Proposition \ref{prp:mr:main:ex} and define
the new process
\begin{equation}
\label{eq:sps:def:V}
V(t) = T_{-\Gamma(t)} [ X(t) + \Phi_{\mathrm{ref}}]
  - \Phi
\end{equation}
for some $\Phi \in \mathcal{U}_{H^1}$.
In addition, we introduce the
nonlinearity
\begin{equation}
\begin{array}{lcl}
\mathcal{R}_{\sigma;\Phi, c}
(v)
& = &
 \kappa_{\sigma}(\Phi + v , \psi_{\mathrm{tw}} )
A_* [\Phi + v]
\\[0.2cm]
& & \qquad
 + f(\Phi + v)
+ \sigma^2  b(\Phi + v, \psi_{\mathrm{tw}} ) \partial_\xi[ g(\Phi + v) ]
\\[0.2cm]
& & \qquad
+ \Big[c + a_{\sigma}\big(\Phi + v , c , \psi_{\mathrm{tw}}
\big) \Big] [\Phi' + v']
\\[0.2cm]
& = &
 \kappa_{\sigma}(\Phi + v , \psi_{\mathrm{tw}} )
\Big[
A_* [\Phi + v]
+ \mathcal{J}_{\sigma}(\Phi + v,
c , \psi_{\mathrm{tw}} )
\Big]
+ a_{\sigma}\big(\Phi + v , c , \psi_{\mathrm{tw}}
\big) [\Phi' + v'] ,
\end{array}
\end{equation}
together with
\begin{equation}
\begin{array}{lcl}
\mathcal{S}_{\Phi}(v)
& = &
     g( \Phi + v )
     + b( \Phi + v , \psi_{\mathrm{tw}}) [\Phi' + v'].
\\[0.2cm]
\end{array}
\end{equation}
Our main result
states that the shifted process $V$ can be interpreted as a weak solution
to the SPDE
\begin{equation}
\begin{array}{lcl}
d V & = &
\mathcal{R}_{\sigma; \Phi, c}(V) \, dt
 + \sigma \mathcal{S}_{\Phi}(V) d\b_t .
\end{array}
\end{equation}

\begin{prop}
\label{prp:sps:props:v}
Consider the setting
of Proposition
\ref{prp:mr:main:ex} and suppose that (hPar)
is satisfied.
Then the map
\begin{equation}
V: [0, T] \times \Omega \to L^2
\end{equation}
defined by \sref{eq:sps:def:V}
satisfies the following properties.
\begin{itemize}
\item[(i)]{
  For almost all $\omega \in \Omega$,
  the map
  $t \mapsto V(t,\omega)$ 
  is of class $C([0,T]; L^2)$.
}
\item[(ii)]{
  For all $t \in [0,T]$, the map
  $\omega \mapsto V(t, \omega) 
    \in L^2$ 
  is $(\mathcal{F}_t)$-measurable.
}

\item[(iii)]{
  We have the inclusion
  \begin{equation}
  \label{eq:sps:props:v:n2:h1}
  \begin{array}{lcl}
     V & \in &  \mathcal{N}^2
       \big( [0 , T] ; (\mathcal{F}_t ) ; H^1 \big)
  \end{array}
  \end{equation}
  together with
  \begin{equation}
  \label{eq:sps:props:sw:n2:h1}
  \begin{array}{lcl}
     \mathcal{S}_{\sigma;\Phi}(V)
      & \in &
       \mathcal{N}^2 \big( [0 , T] ; (\mathcal{F}_t ) ; L^2 \big).
     \end{array}
  \end{equation}
}
\item[(iv)]{
  For almost all $\omega \in \Omega$, we have the inclusion
  \begin{equation}
    \mathcal{R}_{\sigma;\Phi, c}\big( V(\cdot, \omega) \big)
      \in L^1([0,T]; H^{-1}).
  \end{equation}
}
\item[(v)]{
  For almost all $\omega \in \Omega$,
  the identity
  \begin{equation}
  \label{eq:sps:int:eg:for:v}
  \begin{array}{lcl}
  V(t) & = & V(0) + \int_0^t
     \mathcal{R}_{\sigma;\Phi, c}\big( V(s) \big)\, ds
    + \sigma \int_0^t
        \mathcal{S}_{\Phi}\big(V(s)\big)
       d \beta_s
  \end{array}
  \end{equation}
  holds for all $0 \le t \le T$. 
}
\end{itemize}
\end{prop}

Taking derivatives of translation operators typically requires extra regularity of the underlying function,
which prevents us from applying an It\^o formula directly to \sref{eq:sps:def:V}.
In order to circumvent this technical issue,
we pick a test function $\zeta \in C_c^\infty(\Real;\Real^n)$
and consider the two maps
\begin{equation}
\phi_{1;\zeta}: H^{-1} \times \Real \to \Real,
\qquad
\qquad
\phi_{2;\zeta}:  \Real \to \Real
\end{equation}
that act as
\begin{equation}
\begin{array}{lcl}
\phi_{1;\zeta}\big( x , \gamma \big)
  & = &
   \langle x, T_{\gamma} \zeta
       \rangle_{H^{-1} ; H^1 },
\\[0.2cm]
\phi_{2;\zeta}\big( \gamma \big)
 & = &
   \langle T_{-\gamma} \Phi_{\mathrm{ref}} - \Phi , \zeta
       \rangle_{H^{-1} ; H^1 }
\\[0.2cm]
 & = &
     \langle T_{-\gamma} \Phi_{\mathrm{ref}} - \Phi , \zeta
       \rangle_{L^2}.
\end{array}
\end{equation}
These two maps do have sufficient smoothness for our purposes here.

\begin{lem}
Consider the setting
of Proposition
\ref{prp:mr:main:ex}.
Then for almost all $\omega \in \Omega$
the identity
\begin{equation}
\label{eq:sps:ito:i}
\begin{array}{lcl}
\phi_{1;\zeta}\big(X(t), \Gamma(t) \big)
 & = & \phi_{1;\zeta}\big( X(0), \Gamma(0) \big)
\\[0.2cm]
& & \qquad
  + \int_0^t
  \langle A_* [X(s) + \Phi_{\mathrm{ref}} ]  + f(X(s) + \Phi_{\mathrm{ref}}),
    T_{\Gamma(s)} \zeta \rangle_{H^{-1}; H^1} \, ds
\\[0.2cm]
& & \qquad
  - \int_0^t
  \big[
    c + a_{\sigma}\big( X(s) + \Phi_{\mathrm{ref}},c, T_{\Gamma(s)} \psi_{\mathrm{tw}} \big)
  \big]
    \langle X(s) , T_{\Gamma(s)} \zeta'
       \rangle_{L^2} \, d s
\\[0.2cm]
& & \qquad
- \frac{1}{2} \sigma^2
\int_0^t
2 b\big(X(s) + \Phi_{\mathrm{ref}}, T_{\Gamma(s)} \psi_{\mathrm{tw}} \big)
  \langle g(X(s) + \Phi_{\mathrm{ref}}) , T_{\Gamma(s)} \zeta' \rangle_{L^2}
  \, d s
\\[0.2cm]
& & \qquad
  + \frac{1}{2} \sigma^2 \int_0^t
     b\big(X(s) + \Phi_{\mathrm{ref}}, T_{\Gamma(s)} \psi_{\mathrm{tw}} \big)^2
   \langle X(s) , T_{\Gamma(s)} \zeta'' \rangle_{L^2}
    \, d s
\\[0.2cm]
& & \qquad
+ \sigma \int_0^t
  \langle g(X(s) + \Phi_{\mathrm{ref}}) , T_{\Gamma(s)} \zeta \rangle_{L^2} \, d \beta_s
\\[0.2cm]
& & \qquad
  - \sigma \int_0^t
   b\big( X(s) + \Phi_{\mathrm{ref}}  , T_{\Gamma(s)} \psi_{\mathrm{tw}} \big)
    \langle X(s) , T_{\Gamma(s)} \zeta'
       \rangle_{L^2} \, d \beta_s
\end{array}
\end{equation}
holds for all $0 \le t \le T$.
\end{lem}
\begin{proof}
We note that $\phi_{1;\zeta}$ is $C^2$-smooth,
with derivatives given by
\begin{equation}
D \phi_{1;\zeta}(x, \gamma)[y, \beta]
 = \langle y , T_{\gamma} \zeta \rangle_{H^{-1}; H^1}
  - \beta \langle x , T_{\gamma} \zeta'  \rangle_{H^{-1}; H^1} ,
\end{equation}
together with
\begin{equation}
D^2 \phi_{1;\zeta}(x, \gamma)[y, \beta][y, \beta]
 = -2\beta \langle y , T_{\gamma} \zeta' \rangle_{H^{-1}; H^1}
   + \beta^2 \langle x , T_{\gamma} \zeta'' \rangle_{H^{-1}; H^1} .
\end{equation}
Applying a standard It\^o formula such as
\cite[Thm. 1]{DaPratomild} with $S = I$,
the result readily follows.
\end{proof}

\begin{lem}
Consider the setting
of Proposition
\ref{prp:mr:main:ex}.
Then for almost all $\omega \in \Omega$
the identity
\begin{equation}
\label{eq:sps:ito:ii}
\begin{array}{lcl}
\phi_{2;\zeta}\big( \Gamma(t) \big)
 & = & \phi_{2;\zeta}( \Gamma(0) \big)
\\[0.2cm]
& & \qquad
    - \int_0^t
  \big[
    c + a_{\sigma}\big( X(s) + \Phi_{\mathrm{ref}}, c, T_{\Gamma(s)} \psi_{\mathrm{tw}} \big)
  \big]
    \langle \Phi_{\mathrm{ref}}, T_{\Gamma(s)} \zeta'
       \rangle_{L^2} \, d s
\\[0.2cm]
& & \qquad
  + \frac{1}{2} \sigma^2 \int_0^t
     b\big(X(s) + \Phi_{\mathrm{ref}}, T_{\Gamma(s)} \psi_{\mathrm{tw}} \big)^2
   \langle \Phi_{\mathrm{ref}} , T_{\Gamma(s)} \zeta'' \rangle_{L^2}
    \, d s
\\[0.2cm]
& & \qquad
  - \sigma \int_0^t
   b\big( X(s)  + \Phi_{\mathrm{ref}}, T_{\Gamma(s)} \psi_{\mathrm{tw}} \big)
    \langle \Phi_{\mathrm{ref}} , T_{\Gamma(s)} \zeta'
       \rangle_{L^2} \, d \beta_s
\end{array}
\end{equation}
holds for all $0 \le t \le T$.
\end{lem}
\begin{proof}
We note that $\phi_{2;\zeta}$ is $C^2$-smooth,
with derivatives given by
\begin{equation}
D \phi_{2;\zeta}(\gamma)[\beta]
 = - \beta \langle \Phi_{\mathrm{ref}} , T_{\gamma} \zeta'  \rangle_{L^2} ,
\end{equation}
together with
\begin{equation}
D^2 \phi_{2;\zeta}( \gamma)[ \beta][ \beta]
 = \beta^2 \langle \Phi_{\mathrm{ref}} , T_{\gamma} \zeta'' \rangle_{L^2}.
\end{equation}
The result again follows from the It\^o formula.
\end{proof}

\begin{cor}
\label{cor:sps:weak:form:with:zeta}
Consider the setting
of Proposition
\ref{prp:mr:main:ex},
suppose that (hPar) is satisfied
and pick a test-function $\zeta \in C^\infty_c(\Real ;\Real^n)$.
Then for almost all $\omega \in \Omega$,
the map $V$
defined by \sref{eq:sps:def:V}
satisfies
the identity
\begin{equation}
\label{eq:sps:itg:weak:formulation}
  \begin{array}{lcl}
  \langle V(t) , \zeta \rangle_{L^2}
   & = & \langle V(0) , \zeta \rangle_{L^2}
    + \int_0^t
     \langle \mathcal{R}_{\sigma;\Phi, c}\big( V(s) \big)  , \zeta \rangle_{H^{-1}; H^1}   \, ds
    + \sigma \int_0^t
        \langle \mathcal{S}_{\Phi}\big(V(s)\big) , \zeta \rangle_{L^2}
       d \beta_s
  \end{array}
  \end{equation}
for all $0 \le t \le T$.
\end{cor}
\begin{proof}
For any $\gamma \in \Real$,
we have the identities
\begin{equation}
a_{\sigma}( u, c, T_{\gamma} \psi )
=a_{\sigma} (T_{-\gamma} u , c, \psi ),
\qquad
b( u,  T_{\gamma} \psi )
= b(T_{-\gamma} u , \psi ),
\end{equation}
together with the commutation relations
\begin{equation}
T_{\gamma} f(u) = f (T_{\gamma} u),
\qquad
T_{\gamma} g(u) = g (T_{\gamma} u),
\qquad
T_{\gamma} A_* u = A_* T_{\gamma} u.
\end{equation}
By construction, we also have
\begin{equation}
\langle V(t), \zeta \rangle_{L^2}
 = \phi_{1;\zeta}\big(X(t), \Gamma(t) \big)
   + \phi_{2; \zeta}\big( \Gamma(t) \big) ,
\end{equation}
together with
\begin{equation}
T_{-\Gamma(s)} [X(s) + \Phi_{\mathrm{ref}}] = \Phi + V(s).
\end{equation}
The derivatives in \sref{eq:sps:ito:i}
and \sref{eq:sps:ito:ii} can now be
transferred from $\zeta$
to yield \sref{eq:sps:itg:weak:formulation}.

We emphasise that the identity
\begin{equation}
\frac{1}{2} \sigma^2 b(\Phi + V(s)  , \psi_{\mathrm{tw}})^2
[X'' + \Phi_{\mathrm{ref}}'']
= \frac{1}{2\rho} \sigma^2 b(\Phi + V(s)  , \psi_{\mathrm{tw}})^2 A_*
  [X(s) + \Phi_{\mathrm{ref}} ]
\end{equation}
is a crucial ingredient in this computation.
This is where we use the requirement in (HA) that all
the diffusion coefficients in $A_*$ are equal.
\end{proof}

\begin{proof}[Proof of Proposition \ref{prp:sps:props:v}]
Items (i) and (ii) follow immediately
from items (i) and (ii) of Proposition \ref{prp:mr:main:ex}.
Turning to (iii), notice first that we have the isometry
\begin{equation}
\norm{ T_{\gamma} x }_{H^1} = \norm{x}_{H^1}.
\end{equation}
Observe in addition that
\begin{equation}
\norm{ T_{\gamma} \Phi_{\mathrm{ref}} - \Phi }_{H^1}
\le
  \norm{ T_{\gamma} \Phi_{\mathrm{ref}} - \Phi_{\mathrm{ref}} }_{H^1}
   + \norm{\Phi_{\mathrm{ref}} - \Phi}_{H^1}
\le
C_1 \big[ 1 + \abs{\gamma} \big],
\end{equation}
since $\Phi_{\mathrm{ref}}'$ and $\Phi_{\mathrm{ref}}''$
decay exponentially.
In particular,
the inclusion \sref{eq:sps:props:v:n2:h1}
follows from the corresponding
inclusions \sref{eq:prp:mr:ex:inclusions:n2:x:gamma}
for the pair $(X, \Gamma)$.
The second inclusion \sref{eq:sps:props:sw:n2:h1}
now follows immediately from
the bounds in Lemma's \ref{lem:prlm:ests:g}
and \ref{lem:prlm:bnds:b}.

Using Lemma's \ref{lem:prlm:f:bnds},
\ref{lem:prlm:ests:g},
\ref{lem:prlm:bnds:b}
and \ref{lem:var:bnds:a},
we obtain the bound
\begin{equation}
\label{eq:sps:bnd:comp:on:r:sigma}
\begin{array}{lcl}
\norm{ \mathcal{R}_{\sigma;\Phi,c}(v) }_{H^{-1}}
& \le &
  C_2 K_{\kappa} \big[ 1 + \norm{v}_{H^1} \big]
\\[0.2cm]
& & \qquad
+ C_2 \big[ 1 + \norm{v}_{H^1}^2 \norm{v}_{L^2} \big]
\\[0.2cm]
& & \qquad
+ C_2 \sigma^2 K_b \big[ 1 + \norm{v}_{H^1} \big]
\\[0.2cm]
& & \qquad
 + \big[ 1 + \norm{v}_{H^1} \norm{v}_{L^2}^2 \big]
   \big[ 1 + \norm{v}_{H^1} \big].
\end{array}
\end{equation}
Since items (i) and (iii) imply that
\begin{equation}
\sup_{0 \le t \le T} \norm{V(t, \omega)}_{L^2}
+ \int_0^T \norm{V(t, \omega)}_{H^1}^2 \, d t
< \infty
\end{equation}
for almost all $\omega \in \Omega$,
item (iv) follows from the standard
bound
\begin{equation}
\int_0^T
  \norm{V(t, \omega)}_{H^1} \, dt
\le \sqrt{T} \Big[\int_0^T
  \norm{V(t, \omega)}^2_{H^1} \, dt \Big]^{1/2}.
\end{equation}

Finally, we note that items (iii) and (iv)
imply that the integrals in
\sref{eq:sps:int:eg:for:v} are well-defined.
In view of Corollary
\ref{cor:sps:weak:form:with:zeta},
we can apply
a standard diagonalisation argument
involving the separability of $L^2$
and the density of test-functions
to conclude that (v) holds.
\end{proof}

\section{The stochastic time transform}
\label{sec:md}

We note that \sref{eq:sps:int:eg:for:v}
can be interpreted as a quasi-linear
equation due to the presence
of the $\kappa_{\sigma} A_*$ term.
In this section we transform
our problem to a semilinear form
by rescaling time, using the fact that
$\kappa_{\sigma}$ is a scalar.
In addition, we investigate
the impact of this transformation
on the probabilities \sref{eq:mr:def:p:eps:alpha}.

Recalling the map $V$
defined in Proposition
\ref{prp:sps:props:v},
we write
\begin{equation}
\tau_\Phi(t, \omega) = \int_0^t
  \kappa_{\sigma}\big( \Phi + V(s, \omega) , \psi_{\mathrm{tw}} \big)  \, d s .
\end{equation}
Using Lemma \ref{lem:prlm:bnds:kappa}
we see that
$t \mapsto \tau_\Phi(t)$ is a continuous
strictly increasing $(\mathcal{F}_t)$-adapted
process that satisfies
\begin{equation}
\label{eq:mild:ineqs:tau:phi}
t \le \tau_\Phi(t) \le K_{\kappa} t
\end{equation}
for $0 \le t \le T$.
In particular, we can define a map
\begin{equation}
\label{eq:mild:defn:t:phi:spaces}
t_{\Phi}: [0, T] \times \Omega \to [0, T]
\end{equation}
for which
\begin{equation}
\label{eq:mild:defn:t:phi}
\tau_{\Phi}(t_{\Phi}(\tau, \omega), \omega) = \tau.
\end{equation}

We now introduce the
time-transformed map
\begin{equation}
\overline{V}: [0,T] \times \Omega \to L^2
\end{equation}
that acts as
\begin{equation}
\label{eq:mld:def:ovl:v}
\begin{array}{lcl}
\overline{V}(\tau, \omega)
& = & V\big( t_\Phi(\tau, \omega) , \omega \big) .
\\[0.2cm]
\end{array}
\end{equation}
Before stating our main results,
we first investigate
the effects of this transformation
on the terms appearing in
\sref{eq:sps:int:eg:for:v}.

\begin{lem}
Consider the setting of Proposition \ref{prp:mr:main:ex} and suppose that (hPar) is satisfied.
Then the map $t_{\Phi}$ defined in \sref{eq:mild:defn:t:phi:spaces}
satisfies the following properties.
\begin{itemize}
\item[(i)]{
For every $0 \le \tau \le T$,
the random variable $\omega \mapsto t_{\Phi}(\tau, \omega)$
is an
$(\mathcal{F}_t)$-stopping time.
}
\item[(ii)]{
The map $\tau \mapsto t_{\Phi}(\tau, \omega)$
is  continuous and strictly increasing for all $\omega \in \Omega$.
}
\item[(iii)]{
For any $0 \le \tau \le T$ and $\omega \in \Omega$
we have the bounds
\begin{equation}
K_{\kappa}^{-1} \tau \le  t_\Phi(\tau, \omega) \le \tau.
\end{equation}
}
\item[(iv)]{
  For every $0 \le t \le T$,
  the identity
  \begin{equation}
    t_{\Phi}( \tau_\Phi(t, \omega), \omega ) =t
   \end{equation}
   holds on the set $\{\omega: \tau_\Phi(t, \omega) \le T\}$.
}
\end{itemize}
\end{lem}
\begin{proof}
On account of the identity
\begin{equation}
\{ \omega: t_{\Phi}(\tau, \omega) \le t \}
= \{ \omega: \tau_\Phi(t, \omega) \ge \tau \}
\end{equation}
and the fact that the latter set
is in $\mathcal{F}_t$, we may
conclude that $t_{\Phi}(\tau)$
is an
$(\mathcal{F}_t)$-stopping time.
The remaining properties follow
directly from
\sref{eq:mild:ineqs:tau:phi}-\sref{eq:mild:defn:t:phi}.
\end{proof}

\begin{lem}
\label{lem:md:timechange}
Consider the setting of Proposition \ref{prp:mr:main:ex},
recall the maps $(t_{\Phi}, \overline{V} )$
defined by  \sref{eq:mild:defn:t:phi:spaces}
and \sref{eq:mld:def:ovl:v}
and  suppose that (hPar) is satisfied.
Then there exists a filtration
$(\overline{\mathcal{F}}_{\tau})_{\tau \ge 0}$
together with a $(\overline{\mathcal{F}}_{\tau})$-Brownian motion
$(\overline{\beta}_{\tau})_{\tau \ge 0}$
so that for any
$H \in \mathcal{N}^2([0,T] ; (\mathcal{F}_t) ; L^2)$,
the process
\begin{equation}
\overline{H}(\tau, \omega) = H\big( t_\Phi(\tau,\omega), \omega \big)
\end{equation}
satisfies the following properties.
\begin{itemize}
\item[(i)]{
  We have the inclusion
  \begin{equation}
    \overline{H} \in \mathcal{N}^2( [0, T] ;
      (\overline{\mathcal{F}}_{\tau}) ; L^2) ,
  \end{equation}
  together with the bound
  \begin{equation}
    \label{eq:mild:l2:transform:bnd}
    E \int_0^T \norm{\overline{H}(\tau)}_{L^2}^2 \, d \tau
    \le K_{\kappa} E \int_0^T \norm{H(t)}_{L^2}^2 \, d t .
  \end{equation}
}
\item[(ii)]{
  For almost all $\omega \in \Omega$, the identity
  \begin{equation}
    \label{eq:mild:time:transf:st:int:transf}
    \int_0^{t_\Phi(\tau)} H(s) \, d \beta_s
    = \int_0^\tau \overline{H}(\tau')
      \kappa_{\sigma}(\Phi + \overline{V}(\tau') , \psi_{\mathrm{tw}} )^{-1/2}
        \, d \overline{\beta}_{\tau'}
  \end{equation}
  holds for all $0 \le \tau \le T$.
}
\end{itemize}
\end{lem}
\begin{proof}
Following \cite[{\S}1.2.3]{jeanblanc2009mathematical},
we write
\begin{equation}
\overline{\mathcal{F}}_{\tau}
= \{ A \in \cup_{t \ge 0} \mathcal{F}_t : A \cap
\{ t_{\Phi}(\tau) \le t \} \in \mathcal{F}_t
\hbox{ for all  } t \ge 0 \} .
\end{equation}
The fact that $\overline{H}$
is $(\overline{\mathcal{F}}_{\tau})$-progressively
measurable can be established following
the proof of \cite[Lem. 10.8(c)]{jacod2006calcul}.
In addition, we note that
for almost all $\omega \in \Omega$ the path
\begin{equation}
t \mapsto \norm{V(t, \omega)}_{L^2}^2
\end{equation}
is in $L^1([0,T])$, which allows us to
apply the deterministic substitution rule
to obtain
\begin{equation}
\int_{0}^{t_{\Phi}(\tau)}
\norm{V(s)}_{L^2}^2 \, ds
=
\int_0^{\tau}
\norm{V\big(t_{\Phi}(\tau')\big)}^2_{L^2}
\partial_\tau t_{\Phi}(\tau') \, d \tau' .
\end{equation}

We now note that
\begin{equation}
\begin{array}{lcl}
\partial_\tau t_{\Phi}(\tau')
 & = & \big[ \partial_t \tau_{\Phi}\big( t_{\Phi}(\tau') \big) \big]^{-1}
\\[0.2cm]
 & = &
    \kappa_{\sigma}\big(\Phi+ V\big(t_{\Phi}(\tau')\big) , \psi_{\mathrm{tw}} \big)^{-1}
\\[0.2cm]
  & = &
    \kappa_{\sigma}\big(\Phi+ \overline{V}(\tau') , \psi_{\mathrm{tw}} \big)^{-1}.
\end{array}
\end{equation}

In particular, we see that
\begin{equation}
\label{eq:md:low:bnd:deriv:t}
\abs{ \partial_\tau t_{\Phi}(\tau') } \ge K_{\kappa}^{-1}
\end{equation}

and hence
\begin{equation}
\int_0^{\tau}
\norm{V\big(t_{\Phi}(\tau')\big)}^2_{L^2} \, d \tau'
\le
K_{\kappa}
\int_{0}^{t_{\Phi}(\tau)}
\norm{V(s)}_{L^2}^2 \, ds .
\end{equation}

The bound \sref{eq:mild:l2:transform:bnd}
now follows from $t_{\Phi}(T,\omega) \le T$.

To obtain (ii),
we introduce the Brownian-motion $(\overline{\beta}_\tau)_{\tau \ge 0}$
that is given by
\begin{equation}
\overline{\beta}_\tau
= \int_0^\tau \frac{1}{\sqrt{\partial_\tau t_\Phi(\tau')}}
    \,   d \beta_{t_\Phi(\tau')} .
\end{equation}
For any test-function $\zeta \in C^\infty_c(\Real;\Real^n)$
and $0 \le t \le T$,
the proof of \cite[Lem. 5.1.3.5]{jeanblanc2009mathematical}
implies that for almost all $\omega \in \Omega$
the identity
\begin{equation}
\begin{array}{lcl}
\int_0^{t_\Phi(\tau) }
  \langle H(s) , \zeta \rangle_{L^2} d \beta_s
& = &
  \int_0^\tau
  \langle H\big(t_{\Phi}(\tau') \big) , \zeta \rangle_{L^2}
     \sqrt{\partial_\tau t_{\Phi}(\tau')}
     \, d \overline{\beta}_{\tau'}
\\[0.2cm]
& = &
  \int_0^\tau
  \langle \overline{H}(\tau')  , \zeta \rangle_{L^2}
     \kappa_{\sigma}\big(\Phi+ \overline{V}(\tau') , \psi_{\mathrm{tw}} \big)^{-1/2}
     \, d \overline{\beta}_{\tau'}
\end{array}
\end{equation}
holds for all $0 \le \tau \le T$.
Since (i) and (ii) together imply that the right-hand side
of \sref{eq:mild:time:transf:st:int:transf} is well-defined
as a stochastic-integral,
a standard diagonalisation argument involving
the separability of $L^2$
shows that both sides must be equal
for almost all $\omega \in \Omega$.
\end{proof}

In order to formulate the time-transformed SPDE,
we introduce the nonlinearity
\begin{equation}
\label{eq:md:overline:r:sigma}
\begin{array}{lcl}
\overline{\mathcal{R}}_{\sigma;\Phi, c}
(v)
& = &
  \kappa_{\sigma}(\Phi + v , \psi_{\mathrm{tw}})^{-1}
    \mathcal{R}_{\sigma;\Phi,c}(v)
  - \mathcal{L}_{\mathrm{tw}} v
\\[0.2cm]
& = & A_* [\Phi + v]
+ \mathcal{J}_{\sigma}(\Phi + v,
c , \psi_{\mathrm{tw}} )
+ \kappa_{\sigma}(\Phi + v,\psi_{\mathrm{tw}})^{-1}
a(\Phi + v , c, \psi_{\mathrm{tw}})
[\Phi' + v']
\\[0.2cm]
& & \qquad
  - \mathcal{L}_{\mathrm{tw}} v,
\end{array}
\end{equation}
together with
\begin{equation}
\label{eq:md:overline:s:sigma}
\begin{array}{lcl}
\overline{\mathcal{S}}_{\sigma; \Phi}(v)
& = &
\kappa_{\sigma}(\Phi + v,\psi_{\mathrm{tw}})^{-1/2}
\mathcal{S}_{\Phi}(v)
\\[0.2cm]
& = &
\kappa_{\sigma}(\Phi + v,\psi_{\mathrm{tw}})^{-1/2}
\Big[
  g(\Phi + v)
    + b(\Phi + v, \psi_{\mathrm{tw}})[\Phi' + v']
\Big] .
\\[0.2cm]
\end{array}
\end{equation}

\begin{prop}
\label{prp:md:props:overline:v}
Consider the setting
of Proposition \ref{prp:mr:main:ex}
and suppose that (hPar) is satisfied.
Then the map
\begin{equation}
\overline{V}: [0, T] \times \Omega \to L^2
\end{equation}
defined by the transformations
\sref{eq:sps:def:V}
and
\sref{eq:mld:def:ovl:v}
satisfies the following properties.
\begin{itemize}
\item[(i)]{
  For almost all $\omega \in \Omega$,
  the map $\tau \mapsto \overline{V}(\tau; \omega)$
  is of class $C\big([0,T]; L^2 \big)$.
}
\item[(ii)]{
  For all $\tau \in [0, T]$, the map
  $\omega \mapsto \overline{V}(\tau,\omega)$
  is $(\overline{\mathcal{F}}_{\tau})$-measurable.
}
\item[(iii)]{
  We have the inclusion
  \begin{equation}
     \begin{array}{lcl}
        \overline{V} \in \mathcal{N}^2\big([0,T];
          (\overline{\mathcal{F}})_{\tau} ; H^1 \big) ,
     \end{array}
  \end{equation}
  together with
  \begin{equation}
     \begin{array}{lcl}
        \overline{\mathcal{S}}_{\sigma;\Phi}(\overline{V})
          \in \mathcal{N}^2\big([0,T];
          (\overline{\mathcal{F}})_{\tau} ; L^2 \big) .
     \end{array}
  \end{equation}
}
\item[(iv)]{
  For almost all $\omega \in \Omega$, we have the inclusion
  \begin{equation}
    \overline{\mathcal{R}}_{\sigma;\Phi, c}\big( \overline{V}(\cdot, \omega) \big)
      \in L^1([0,T]; L^2).
  \end{equation}
}
\item[(v)]{
  For almost all $\omega \in \Omega$, the identity
  \begin{equation}
  \begin{array}{lcl}
    \overline{V}(\tau)
     & = & \overline{V}(0)
     + \int_0^\tau \Big[ \mathcal{L}_{\mathrm{tw}} \overline{V}(\tau')
       +  \overline{\mathcal{R}}_{\sigma;\Phi,c}\big(\overline{V}(\tau')\big)
       \Big]
       \, d \tau'
\\[0.2cm]
& & \qquad \qquad \qquad
     + \sigma \int_0^\tau \overline{\mathcal{S}}_{\sigma;\Phi}\big(\overline{V}(\tau')\big)
       \, d \beta_{\tau'}
  \end{array}
  \end{equation}
  holds for all $0 \le t \le T$.
}
\item[(vi)]{
  For almost all $\omega \in \Omega$,
  the identity
  \begin{equation}
  \label{eq:md:int:eq:mild:formulation}
  \begin{array}{lcl}
    \overline{V}(\tau)
     & = & S(\tau) \overline{V}(0)
     + \int_0^\tau S(\tau - \tau')
         \overline{\mathcal{R}}_{\sigma;\Phi,c}
           \big(\overline{V}(\tau')\big)
       \, d \tau'
   \\[0.2cm]
   & & \qquad
     + \sigma \int_0^\tau S(\tau - \tau')
         \overline{\mathcal{S}}_{\sigma;\Phi}
           \big(\overline{V}(\tau')\big)
              \, d \beta_{\tau'}
  \end{array}
  \end{equation}
  holds for all $\tau \in [0, T]$, in which
  \begin{equation}
    S: [0, \infty) \to \mathcal{L}(L^2;L^2)
  \end{equation}
  denotes the analytic semigroup generated by $\mathcal{L}_{\mathrm{tw}}$.
}
\end{itemize}
\end{prop}
\begin{proof}
Items (i)-(iii) follow by
applying (i) of Lemma \ref{lem:md:timechange}
to the maps $V$, $\partial_\xi V$
and using the definition
\sref{eq:md:overline:s:sigma}.
Item (iv) can be obtained from the computation
\sref{eq:sps:bnd:comp:on:r:sigma},
noting that the $A_* v$ contribution
is no longer present.

Item (v) can be obtained
by applying the stochastic time-transform
\sref{eq:mild:time:transf:st:int:transf}
and the deterministic time-transform
\begin{equation}
\begin{array}{lcl}
\int_0^{t_{\Phi}(\tau)}
  \mathcal{R}_{\sigma;\Phi,c}\big( V(s) \big) \, ds
& = & \int_0^{\tau}
  \overline{\mathcal{R}}_{\sigma;\Phi,c}
    \big( \overline{V}(\tau') \big)
    \big[\kappa_{\sigma}\big(\Phi + \overline{V}(\tau'),
       \psi_{\mathrm{tw}} \big) \big]^{-1} \, d \tau'
\end{array}
\end{equation}
to the integral equation
\sref{eq:sps:int:eg:for:v}.

Turning to (vi), we note that
$A_*$ generates a standard diagonal
heat-semigroup,
which is obviously analytic.
Noting that
\begin{equation}
\mathcal{L}_{\mathrm{tw}} - A_*
 \in \mathcal{L}( H^1 ; L^2)
\end{equation}
and recalling the interpolation estimate
\begin{equation}
\norm{ v}_{H^1} \le
C_1 \norm{ v}_{H^2}^{1/2}
    \norm{v}_{L^2}^{1/2},
\end{equation}
we may apply \cite[Prop 3.2.2(iii)]{lorenzi2004analytic}
to conclude that also $\mathcal{L}_{\mathrm{tw}}$
generates an analytic semigroup.
We may now
apply \cite[Prop 6.3]{DaPratoZab}
and the computation in the proof of
\cite[Prop 4.1.4]{lorenzi2004analytic}
to conclude the integral
identity
\sref{eq:md:int:eq:mild:formulation}.
\end{proof}

We now introduce the scalar functions
\begin{equation}
\label{eq:md:def:n}
\begin{array}{lcl}
N_{\e,\alpha}(t)
 & = &
e^{\alpha t}\norm{V(t)}_{L^2}^2
 + \int_0^t e^{- \e (t - s) }
      e^{\alpha s}\norm{  V(s)}_{H^1}^2 \, ds ,
\\[0.2cm]
\overline{N}_{\e,\alpha}(\tau)
& = &
e^{\alpha \tau}\norm{ \overline{V}(\tau)}_{L^2}^2
 + \int_0^\tau e^{- \e (\tau - \tau') }
    e^{\alpha \tau'} \norm{  \overline{V}(\tau')}_{H^1}^2 \, d\tau' ,
\end{array}
\end{equation}
together with the associated probabilities
\begin{equation}
\begin{array}{lcl}
p_{\e,\alpha}(T,\eta)
 &= & P\Big(
 \sup_{0 \le t \le T}  N_{\e,\alpha}(t)
 > \eta
\Big) ,
\\[0.2cm]
\overline{p}_{\e,\alpha}(T,\eta) & = & P\Big(
 \sup_{0 \le \tau \le T}  \overline{N}_{\e,\alpha}(\tau)
 > \eta
\Big) .
\end{array}
\end{equation}
Our second main result shows that
these two sets of probabilities
can be effectively compared with each other.

\begin{prop}
\label{prp:md:prob:after:time:transf}
Consider the setting of Proposition \ref{prp:mr:main:ex}
and recall the maps $V$ and $\overline{V}$
defined by \sref{eq:sps:def:V}
and
\sref{eq:mld:def:ovl:v}.
Then we have the bound
\begin{equation}
p_{\e,\alpha}(T,\eta)
\le \overline{p}_{ K_{\kappa}^{-1} \e, \alpha}
  \big( K_{\kappa} T , K_{\kappa}^{-1} \eta \big) .
\end{equation}
\end{prop}
\begin{proof}
We note that
\begin{equation}
e^{\alpha t}\norm{ V(t)}_{L^2}
= e^{\alpha t } \norm{
  \overline{V}\big( \tau_{\Phi}(t) \big) }_{L^2}
\le
  e^{\alpha \tau_{\Phi}(t) }
  \norm{  \overline{V}\big( \tau_{\Phi}(t) \big) }_{L^2} ,
\end{equation}
which implies that
\begin{equation}
\mathrm{sup}_{0 \le t \le  T}
 e^{\alpha t}\norm{ V(t)}_{L^2}^2
 \le \sup_{0 \le \tau \le K_{\kappa} T}
    e^{\alpha \tau} \norm{ \overline{V}(\tau)}_{L^2}^2 .
\end{equation}
In addition, we compute
\begin{equation}
\int_0^t e^{-\e(t -s)}
  e^{\alpha s}\norm{  V(s)}_{H^1}^2 \, ds
= \int_0^{\tau_\Phi(t)}
    e^{-\e(t - t_\Phi(\tau'))}
       e^{\alpha t_\Phi(\tau')} \norm{\overline{V}(\tau')}_{H^1}^2
      \kappa_{\sigma}\big(\Phi + \overline{V}(\tau'), \psi_{\mathrm{tw}} \big)^{-1}
       \, d \tau' .
\end{equation}
Using \sref{eq:md:low:bnd:deriv:t}
we obtain the estimate
\begin{equation}
t - t_\Phi(\tau') = t_\Phi\big(\tau_\Phi(t)\big) - t_\Phi(\tau')
= \int_{\tau'}^{\tau_\Phi(t)} \partial_\tau t_\Phi(\tau'') \, d \tau''
\ge K_{\kappa}^{-1} \abs{\tau_\Phi(t) - \tau'} ,
\end{equation}
which yields
\begin{equation}
\int_0^t e^{-\e(t -s)}
  e^{\alpha s}\norm{  V(s)}_{H^1}^2 \, ds
\le K_{\kappa} \int_0^{\tau_\Phi(t)}
    e^{-K_{\kappa}^{-1} \e(\tau_\Phi(t) - \tau') }
      e^{\alpha \tau'}\norm{\overline{V}(\tau')}_{H^1}^2
       \, d \tau'  .
\end{equation}
In particular, we conclude that
\begin{equation}
\sup_{0 \le t \le T} \int_0^t e^{-\e(t -s)}
  e^{\alpha s}\norm{ V(s)}_{H^1}^2 \, ds
\le
\sup_{0 \le \tau \le K_{\kappa} T }
 K_{\kappa} \int_0^{\tau}
    e^{-K_{\kappa}^{-1} \e(\tau - \tau') }
      e^{\alpha \tau'} \norm{\overline{V}(\tau')}_{H^1}^2
       \, d \tau' .
\end{equation}
This yields the implication
\begin{equation}
\sup_{0 \le \tau \le K_{\kappa} T }
 \overline{N}_{K_{\kappa}^{-1} \e , \alpha}(\tau) \le
   K_{\kappa}^{-1} \eta
\Rightarrow
\sup_{0 \le t \le  T }
  N_{\e,\alpha}(t) \le \eta ,
\end{equation}
from which the desired inequality immediately follows.
\end{proof}

\section{The stochastic wave}
\label{sec:swv}

In this section we set out to construct the
branch of modified waves $(\Phi_{\sigma}, c_{\sigma})$
and analyze the phase condition
\begin{equation}
\label{eq:swv:phase:cond}
\langle T_{-\gamma_0} [u_0] - \Phi_{\sigma} , \psi_{\mathrm{tw}} \rangle_{L^2} = 0
\end{equation}
for $u_0 \approx \Phi_{\sigma}$.
In particular, we establish
Propositions \ref{prp:mr:swv:ex},
\ref{prp:mr:phase:shift} and \ref{prp:mr:expl:waves}.

A key role in our analysis is reserved
for the function
\begin{equation}
\label{eq:swv:def:m:sigma}
\begin{array}{lcl}
\mathcal{M}_{\sigma;\Phi, c}(v , d)
& = &
\mathcal{J}_{\sigma}
  (\Phi + v, c + d , \psi_{\mathrm{tw}})
 -  \mathcal{J}_0(\Phi, c)
\\[0.2cm]
& & \qquad
  -   d \Phi'_{0}   + [ A_*- \mathcal{L}_{\mathrm{tw}} ] v ,
\\[0.2cm]
\end{array}
\end{equation}
defined for $(\Phi, c) \in \mathcal{U}_{H^1} \times \Real$
and $(v, d) \in H^1 \times \Real$.
Indeed, we will construct a solution to
\begin{equation}
\label{eq:swv:st:wv:eqn:to:solve}
A_* \Phi_{\sigma}
  + \mathcal{J}_{\sigma}(\Phi_{\sigma}, c_{\sigma}, \psi_{\mathrm{tw}}) = 0
\end{equation}
by writing
\begin{equation}
\Phi_{\sigma} = \Phi_0 + v,
\qquad
c_{\sigma} = c_0 + d.
\end{equation}
Using the fact that the pair $(\Phi_0, c_0)$
is a solution to \sref{eq:swv:st:wv:eqn:to:solve}
for $\sigma = 0$,
one readily verifies
that the pair $(v,d) \in H^2 \times \Real$
must satisfy the system
\begin{equation}
 d \Phi'_{0}   + \mathcal{L}_{\mathrm{tw}} v
= - \mathcal{M}_{\sigma;\Phi_0, c_{0}}(v, d).
\end{equation}
In addition, the function
$\mathcal{M}_{\sigma;\Phi_{\sigma}, c_{\sigma}}$
will be used in \S\ref{sec:fnl}
to obtain bounds
on the nonlinearity
$\overline{\mathcal{R}}_{\sigma;\Phi_{\sigma}, c_{\sigma} }$.

In \S\ref{sec:swv:bnd:m:sigma} we obtain
global and Lipschitz bounds on
$\mathcal{M}_{\sigma;\Phi, c}$.
These bounds are subsequently
used in \S\ref{sec:swv:fixp} to setup
two fixed-point constructions that
provide solutions
to
\sref{eq:swv:phase:cond}
and \sref{eq:swv:st:wv:eqn:to:solve}.

\subsection{Bounds for $\mathcal{M}_{\sigma}$}
\label{sec:swv:bnd:m:sigma}

In order to streamline our estimates,
it is convenient
to decompose the function
$\mathcal{J}_{\sigma}$
as
\begin{equation}
\begin{array}{lcl}
\mathcal{J}_{\sigma}(u , \overline{c}, \psi_{\mathrm{tw}})
& = &
  \kappa_{\sigma}(u , \psi_{\mathrm{tw}})^{-1}
  \Big[
    f( u )
    + \overline{c} u'
    + \sigma^2 b( u, \psi_{\mathrm{tw}})
      \partial_\xi [g(u )]
  \Big]
\\[0.2cm]
& = &
    \mathcal{J}_0(u, \overline{c})
    + \mathcal{E}_{\sigma; I}(u, \overline{c})
    + \mathcal{E}_{\sigma; II}(u).
\end{array}
\end{equation}
Here we have introduced the function
\begin{equation}
\begin{array}{lcl}
\mathcal{E}_{\sigma;I}(u, \overline{c})
& = &
\nu^{(-1)}_{\sigma}(u, \psi_{\mathrm{tw}}) [f(u) + \overline{c} u']
\\[0.2cm]
& = &
\nu^{(-1)}_{\sigma}(u, \psi_{\mathrm{tw}}) \mathcal{J}_0(u , \overline{c}) ,
\end{array}
\end{equation}
together with
\begin{equation}
\begin{array}{lcl}
\mathcal{E}_{\sigma;II}(u)
& = &
 \sigma^2
 \kappa_{\sigma}(u, \psi_{\mathrm{tw}})^{-1}
 b(u , \psi_{\mathrm{tw}}) \partial_\xi [g(u)]
\end{array}
\end{equation}
where $\nu_\s^{-1}$ is as defined in (\ref{eq:prlm:defs:nu}).\\
This decomposition
allows us to rewrite
\sref{eq:swv:def:m:sigma}
in the intermediate form
\begin{equation}
\begin{array}{lcl}
\mathcal{M}_{\sigma;\Phi, c}(v , d)
& = &
 \mathcal{M}_{0;\Phi, c}
         (v , d)
    + \mathcal{E}_{\sigma;I}(\Phi + v, c+ d)
    + \mathcal{E}_{\sigma;II}(\Phi + v) .
\end{array}
\end{equation}

We now make a final splitting
\begin{equation}
\begin{array}{lcl}
\mathcal{M}_{0;\Phi,c}(v, d)
 & = & \mathcal{J}_0(\Phi + v, c + d )
- \mathcal{J}_0(\Phi , c  )
- Df(\Phi_0) v - c_{0} v'
 - d \Phi'_{0}
\\[0.2cm]
& = &
  \mathcal{N}_{I;f,\Phi} (v)
+ \mathcal{N}_{II;\Phi, c}(v,d),
\end{array}
\end{equation}
in which we have introduced the function
\begin{equation}
\mathcal{N}_{I;f,\Phi}(v) =
  f(\Phi + v) - f(\Phi) - Df(\Phi) v,
\end{equation}
together with
\begin{equation}
\begin{array}{lcl}
\mathcal{N}_{II; \Phi, c}(v, d)
& = &
   d v'
+ \big[Df(\Phi) -  Df(\Phi_0 )   \big] v
 + (c - c_{0}) v'
 + d [\Phi' - \Phi'_{0} ] .
\end{array}
\end{equation}
We hence arrive at the convenient final expression
\begin{equation}
\label{eq:swv:id:for:m:sigma}
\begin{array}{lcl}
\mathcal{M}_{\sigma;\Phi, c}(v , d)
& = &
  \mathcal{N}_{I;f,\Phi}(v)
  +  \mathcal{N}_{II;\Phi, c}(v,d)
   + \mathcal{E}_{\sigma;I}(\Phi + v, c+ d)
    + \mathcal{E}_{\sigma;II}(\Phi + v)
\end{array}
\end{equation}
and set out to analyze each of these
terms separately.

\begin{lem}
Suppose that (Hf) and (hPar) are satisfied.
Then there exists $K > 0$ so that
for any $v \in H^1$ we have the bound
\begin{equation}
\label{eq:swv:glb:bnd:nf}
\begin{array}{lcl}
 \norm{\mathcal{N}_{I;f,\Phi}( v)}_{L^2} & \le &
   K \big[1 + \norm{v}_{H^1} \big]
     \norm{v}_{H^1} \norm{v}_{L^2} ,
\\[0.2cm]
\end{array}
\end{equation}
while for any pair $(v_A, v_B) \in H^1 \times H^1$
we have the estimates
\begin{equation}
\label{eq:swv:lip:bnd:nf}
\begin{array}{lcl}
\norm{\mathcal{N}_{I;f,\Phi}(v_A) - \mathcal{N}_{I;f,\Phi}(v_B) }_{L^2}
  & \le &
  K \big[ 1 + \norm{v_A}_{H^1} + \norm{v_B}_{H^1} \big]
    \big[ \norm{v_A}_{H^1}  + \norm{v_B}_{H^1} \big]
\\[0.2cm]
& & \qquad \qquad \times
     \norm{v_A-v_B}_{L^2} ,
\\[0.2cm]
\abs{\langle \mathcal{N}_{I;f,\Phi}( v_A )
  -  \mathcal{N}_{I;f,\Phi}(v_B), \psi_{\mathrm{tw}} \rangle_{L^2} }
 & \le &
  K \big[ 1 + \norm{v_A}_{H^1} + \norm{v_B}_{H^1} \big]
    \big[ \norm{v_A}_{L^2}  + \norm{v_B}_{L^2} \big]
\\[0.2cm]
& & \qquad \qquad \times
     \norm{v_A-v_B}_{L^2} .
\end{array}
\end{equation}
\end{lem}
\begin{proof}
Using \sref{eq:prlm:f:bnd:d2:f} and (hPar) we obtain the pointwise
bound
\begin{equation}
\abs{\mathcal{N}_{I;f,\Phi}(v) }
\le C_1  [1 + \abs{v} ] \abs{v}^2,
\end{equation}
from which \sref{eq:swv:glb:bnd:nf} easily follows.
In addition, we may compute
\begin{equation}
\begin{array}{lcl}
\mathcal{N}_{I;f,\Phi}(v_A)
 - \mathcal{N}_{I;f,\Phi}(v_B)
& = & f(\Phi + v_A) - f(\Phi + v_B) - Df(\Phi + v_B)(v_A - v_B)
\\[0.2cm]
& & \qquad
+ \big( Df(\Phi + v_B) - Df(\Phi) \big) (v_A - v_B)
\\[0.2cm]
& = & \mathcal{N}_{I;f, \Phi + v_B}(v_A - v_B)
  + \big(Df(\Phi + v_B) - Df(\Phi) \big) (v_A - v_B).
\end{array}
\end{equation}
Applying \sref{eq:prlm:f:bnd:d2:f} and (hPar)
a second time, we obtain the pointwise bound
\begin{equation}
\begin{array}{lcl}
\abs{\mathcal{N}_{I;f,\Phi}(v_A)
 - \mathcal{N}_{I;f,\Phi}(v_B)}
& \le & C_2 [1 + \abs{v_A} + \abs{v_B} \big]
              \abs{ v_A-v_B}^2
\\[0.2cm]
& & \qquad
  + C_2 [1 + \abs{v_B} ] \abs{v_B} \abs{v_A -v_B}
\\[0.2cm]
& \le &
  C_3 [1 + \abs{v_A} + \abs{v_B} \big]
      \big[ \abs{v_A} + \abs{v_B} \big]
              \abs{ v_A-v_B},
\end{array}
\end{equation}
from which the estimates in \sref{eq:swv:lip:bnd:nf}
can be readily obtained.
\end{proof}

\begin{lem}
\label{lem:swv:mzero}
Suppose that (Hf) and (hPar) are satisfied.
Then there exists $K > 0$ so that
for any $(v, d) \in H^1 \times \Real$
we have the bound
\begin{equation}
\label{eq:swv:mzero:glb:bnd}
\begin{array}{lcl}
 \norm{\mathcal{N}_{II;\Phi, c}(v, d)}_{L^2}
  & \le &
    K \big[
   \abs{c - c_0} + \norm{ \Phi- \Phi_0}_{H^1}
   + \abs{d}
  \big] \big[
    \norm{v}_{H^1} + \abs{d}
  \big] ,
\end{array}
\end{equation}
while for any set of pairs $(v_A, v_B) \in H^1 \times H^1$
and $(d_A, d_B) \in \Real^2$
the expression
\begin{equation}
\Delta_{AB} \mathcal{N}_{II;\Phi, c}
 =  \mathcal{N}_{II;\Phi, c}(v_A, d_A)
 - \mathcal{N}_{II;\Phi, c}(v_B, d_B)
\end{equation}
satisfies the estimates
\begin{equation}
\label{eq:swv:mzero:lip:bnd}
\begin{array}{lcl}
\norm{ \Delta_{AB} \mathcal{N}_{II;\Phi, c}  }_{L^2}
  & \le &
 K \big[
  \norm{v_A}_{H^1} + \abs{d_B}
   +  \norm{\Phi - \Phi_0}_{H^1}
  + \abs{c - c_0} \big]
\\[0.2cm]
& & \qquad \qquad \times
  \big[
     \norm{v_A-v_B}_{H^1}
      + \abs{d_A - d_B}
  \big] ,
\\[0.2cm]
\abs{\langle \Delta_{AB} \mathcal{N}_{II;\Phi, c} , \psi_{\mathrm{tw}} \rangle_{L^2} }
 & \le &
 K \big[
  \norm{v_A}_{L^2} + \abs{d_B}
   +  \norm{\Phi - \Phi_0}_{L^2}
  + \abs{c - c_0} \big]
\\[0.2cm]
& & \qquad \qquad \times
  \big[
     \norm{v_A-v_B}_{L^2}
      + \abs{d_A  - d_B}
  \big] .
\\[0.2cm]
\end{array}
\end{equation}
\end{lem}
\begin{proof}
In view of (hPar), we obtain
the pointwise bound
\begin{equation}
\abs{\mathcal{N}_{II;\Phi, c}(v, d)}
\le
\big[ \abs{d} + \abs{c -c_0} \big] \abs{v'}
+ C_1 \abs{\Phi - \Phi_0}
  \abs{v}
  + \abs{\Phi' - \Phi_0'}
    \abs{d},
\end{equation}
from which \sref{eq:swv:mzero:glb:bnd}
follows.
In addition, we obtain the pointwise bound
\begin{equation}
\begin{array}{lcl}
\abs{\Delta_{AB} \mathcal{N}_{II;\Phi, c} }
& \le &
\abs{d_A - d_B } \abs{v_A'}
+ \big[
 \abs{d_B}
+ \abs{c - c_0} \big] \abs{v_A' - v_B'}
\\[0.2cm]
& & \qquad
+ K \abs{\Phi - \Phi_0 }
  \big[
  \abs{v_A-v_B}
  \big]
  + \abs{\Phi' - \Phi_0'}
   \abs{d_A - d_B}
\end{array}
\end{equation}
from which
\sref{eq:swv:mzero:lip:bnd} follows.
\end{proof}

\begin{lem}
Suppose that (Hf), (Hg) and (hPar) are satisfied.
Then there exists $K > 0$ so that
for any $0 \le \sigma \le 1$
and $(v, d) \in H^1 \times \Real$,
we have the bound
\begin{equation}
\label{eq:swv:ei:glb:bnds}
\begin{array}{lcl}
 \norm{\mathcal{E}_{\sigma;I}(\Phi + v, c + d)}_{L^2} & \le &
 K \sigma^2 (1 + \abs{d})
 \big[
   1 + \norm{v}_{H^1} + \norm{v}_{H^1}^2 \norm{v}_{L^2}
 \big] ,
\\[0.2cm]
\end{array}
\end{equation}
while for any $0 \le \sigma \le 1$
and any set of pairs
$(v_A, v_B) \in H^1 \times H^1$
and $(d_A, d_B) \in \Real^2$,
the expression
\begin{equation}
\Delta_{AB} \mathcal{E}_{\sigma;I}
=  \mathcal{E}_{\sigma;I}
       (\Phi + v_A, c + d_A )
  - \mathcal{E}_{\sigma;I}
       (\Phi + v_B, c + d_B)
\end{equation}
satisfies the estimates
\begin{equation}
\label{eq:swv:ei:lip:bnds}
\begin{array}{lcl}
\norm{ \Delta_{AB} \mathcal{E}_{\sigma;I} }_{L^2}
  & \le &
K \sigma^2
(1 + \abs{d_A} )
 \big[
   1 + \norm{v_A}_{H^1} + \norm{v_A}_{H^1}^2 \norm{v_A}_{L^2}
 \big] \norm{v_A-v_B}_{L^2}
\\[0.2cm]
& & \qquad
  +K \sigma^2
  \big[1 + \abs{d_B} + \norm{v_A}_{H^1} \norm{v_A}_{L^2}
    + \norm{v_B}_{H^1} \norm{v_B}_{L^2} \big]
    \norm{v_A-v_B}_{H^1}
\\[0.2cm]
& & \qquad  + K\sigma^2 \big[1 + \norm{v_A}_{H^1} \big]
   \abs{d_A - d_B } ,
\\[0.2cm]
\abs{\langle
\Delta_{AB} \mathcal{E}_{\sigma;I}
, \psi_{\mathrm{tw}} \rangle_{L^2} }
 & \le &
K \sigma^2
(1 + \abs{d_A} + \abs{d_B} )
 \big[
   1  + \norm{v_A}_{H^1} \norm{v_A}_{L^2}^2
 \big] \norm{v_A-v_B}_{L^2}
\\[0.2cm]
& & \qquad
  +K \sigma^2
  \big[  \norm{v_A}_{L^2}^2 + \norm{v_B}_{L^2}^2 \big]
    \norm{v_A-v_B}_{H^1}
\\[0.2cm]
& & \qquad
+ K\sigma^2 [1 + \norm{v_A}_{L^2} ] \abs{d_A - d_B }.
\end{array}
\end{equation}
\end{lem}
\begin{proof}
The bound \sref{eq:swv:ei:glb:bnds}
follows directly
from Lemmas \ref{lem:prlm:j0:bnds}
and \ref{lem:prlm:nu}.
In addition,
these results allow us
to compute
\begin{equation}
\begin{array}{lcl}
\norm{
 \Delta_{AB} \mathcal{E}_{\sigma;I}
}_{L^2}
& \le &
\abs{\nu^{(-1)}_{\sigma}(\Phi + v_A,\psi_{\mathrm{tw}})
  - \nu^{(-1)}_{\sigma}(\Phi + v_B,\psi_{\mathrm{tw}}) }
 \norm{ \mathcal{J}_0 (\Phi + v_A , c + d_A) }_{L^2}
\\[0.2cm]
& & \qquad
  + \abs{\nu^{(-1)}_{\sigma}(\Phi + v_B, \psi_{\mathrm{tw}}) }
  \norm{ \mathcal{J}_0(\Phi + v_A , c + d_A)
      - \mathcal{J}_0(\Phi + v_B , c + d_B) }_{L^2}
\\[0.2cm]
& \le &
C_1 \sigma^2
 \norm{v_A-v_B}_{L^2}
(1 + \abs{d_A} )
 \big[
   1 + \norm{v_A}_{H^1} + \norm{v_A}_{H^1}^2 \norm{v_A}_{L^2}
 \big]
\\[0.2cm]
& & \qquad
  + C_1 \sigma^2 \big[\norm{v_A}_{H^1} \norm{v_A}_{L^2}
    + \norm{v_B}_{H^1} \norm{v_B}_{L^2} \big] \norm{v_A-v_B}_{H^1}
\\[0.2cm]
& & \qquad  + C_1\sigma^2 [1 + \norm{v_A}_{H^1} ]
    \abs{d_A - d_B }
\\[0.2cm]
& & \qquad
  + C_1 \sigma^2 ( 1 + \abs{d_B} ) \norm{v_A-v_B}_{H^1},
\\[0.2cm]
\end{array}
\end{equation}
together with
\begin{equation}
\begin{array}{lcl}
\abs{ \langle
 \Delta_{AB} \mathcal{E}_{\sigma;I}
, \psi_{\mathrm{tw}} \rangle_{L^2} }
& \le &
\abs{\nu^{(-1)}_{\sigma}(\Phi + v_A, \psi_{\mathrm{tw}})
   - \nu^{(-1)}_{\sigma}(\Phi + v_B, \psi_{\mathrm{tw}}) }
 \abs{ \langle \mathcal{J}_0 (\Phi + v_A , c + d_A) , \psi_{\mathrm{tw}} \rangle_{L^2} }
\\[0.2cm]
& & 
  + \abs{\nu^{(-1)}_{\sigma}(\Phi + v_B, \psi_{\mathrm{tw}}) }
  \abs{ \langle \mathcal{J}_0(\Phi + v_A , c + d_A)
      - \mathcal{J}_0(\Phi + v_B , c + d_B) , \psi_{\mathrm{tw}} \rangle_{L^2} }
\\[0.2cm]
& \le &
C_2 \sigma^2
 \norm{v_A-v_B}_{L^2}
(1 + \abs{d_A} )
 \big[
   1  + \norm{v_A}_{H^1} \norm{v_A}_{L^2}^2
 \big]
\\[0.2cm]
& & 
  + C_2 \sigma^2
    [ \norm{v_A}_{L^2}^2 +  \norm{v_B}_{L^2}^2 ] \norm{v_A-v_B}_{H^1}
\\[0.2cm]
& & 
+ C_2 \sigma^2 [1 + \norm{v_A}_{L^2} ] \abs{d_A - d_B }
\\[0.2cm]
& & 
  + C_2 \sigma^2 ( 1 + \abs{d_B} ) \norm{v_A-v_B}_{L^2} .
\\[0.2cm]
\end{array}
\end{equation}
These terms can all be absorbed
by the expressions in \sref{eq:swv:ei:lip:bnds}.
\end{proof}

\begin{lem}
Suppose that (Hg) and (hPar) are satisfied.
Then there exists $K > 0$
so that
for any $0 \le \sigma \le 1$
and $v\in H^1$
we have the bound
\begin{equation}
\label{eq:swv:glb:bnd:e:ii}
\begin{array}{lcl}
 \norm{\mathcal{E}_{\sigma;II}(\Phi + v)}_{L^2} & \le &
  K \sigma^2
    \big[ 1 + \norm{v}_{H^1} \big],
\\[0.2cm]
\end{array}
\end{equation}
while for any $0 \le \sigma \le 1$
and any pair $(v_A, v_B) \in H^1 \times H^1$
the expression
\begin{equation}
\Delta_{AB} \mathcal{E}_{\sigma;II}
 = \mathcal{E}_{\sigma;II}
       (\Phi + v_A  )
  - \mathcal{E}_{\sigma;II}
       (\Phi + v_B)
\end{equation}
satisfies the estimates
\begin{equation}
\label{eq:swv:lip:bnd:e:ii}
\begin{array}{lcl}
\norm{\Delta_{AB} \mathcal{E}_{\sigma;II} }_{L^2}
  & \le &
K \sigma^2
   \big[ 1 + \norm{v_A}_{H^1} \big]
     \norm{v_A-v_B}_{H^1} ,
\\[0.2cm]
\abs{\langle
\Delta_{AB} \mathcal{E}_{\sigma;II}
, \psi_{\mathrm{tw}} \rangle_{L^2} }
 & \le &
   K \sigma^2
    \big[ 1 + \norm{v_A}_{L^2} \big]
    \norm{v_A-v_B}_{L^2}.
\end{array}
\end{equation}
\end{lem}
\begin{proof}
The bound \sref{eq:swv:glb:bnd:e:ii}
follows directly
from Lemmas
\ref{lem:prlm:ests:g},
\ref{lem:prlm:bnds:b}
and \ref{lem:prlm:bnds:kappa}.
In addition,
we may compute
\begin{equation}
\begin{array}{lcl}
\norm{\Delta_{AB} \mathcal{E}_{\sigma;II}
}_{L^2}
& \le &
\sigma^2 \abs{\nu^{(-1)}_{\sigma}(\Phi + v_A, \psi_{\mathrm{tw}})
   - \nu^{(-1)}_{\sigma}(\Phi + v_B, \psi_{\mathrm{tw}})
} K_b \norm{\partial_\xi [g(\Phi + v_A) ] }_{L^2}
\\[0.2cm]
& & \qquad
  + \sigma^2 K_{\nu}
  \abs{b(\Phi + v_A, \psi_{\mathrm{tw}})
      - b(\Phi + v_B, \psi_{\mathrm{tw}}) }
  \norm{\partial_\xi [g(\Phi + v_A) ] }_{L^2}
\\[0.2cm]
& & \qquad
  + \sigma^2 K_{\kappa}
      K_b
      \norm{\partial_{\xi}[g(\Phi + v_A) - g(\Phi + v_B) ] }_{L^2}
\\[0.2cm]
& \le &
C_1 \sigma^2
 \norm{v_A-v_B}_{L^2} \big[ 1 + \norm{v_A}_{H^1} \big]
\\[0.2cm]
& & \qquad
+ C_1 \sigma^2
 \norm{v_A-v_B}_{L^2}
\big[ 1 + \norm{v_A}_{H^1} \big]
\\[0.2cm]
& & \qquad
 + C_1 \sigma^2
 \big[ 1 + \norm{v_A}_{H^1}
 \big]
 \norm{v_A-v_B}_{H^1} ,
\\[0.2cm]
\end{array}
\end{equation}
together with
\begin{equation}
\begin{array}{lcl}
\abs{
 \langle \Delta_{AB} \mathcal{E}_{\sigma;II}
 , \psi_{\mathrm{tw}} \rangle_{L^2}
}
& \le &
\sigma^2 \abs{\nu^{(-1)}_{\sigma}(\Phi + v_A, \psi_{\mathrm{tw}})
   - \nu^{(-1)}_{\sigma}(\Phi + v_B, \psi_{\mathrm{tw}})
} K_b \abs{ \langle\partial_\xi [g(\Phi + v_A) ] , \psi_{\mathrm{tw}}
  \rangle_{L^2} }
\\[0.2cm]
& & \qquad
  + \sigma^2 K_{\nu}
  \abs{b(\Phi + v_A, \psi_{\mathrm{tw}})
      - b(\Phi + v_B, \psi_{\mathrm{tw}}) }
  \abs{\langle
    \partial_\xi [g(\Phi + v_A) ], \psi_{\mathrm{tw}}
    \rangle_{L^2} }
\\[0.2cm]
& & \qquad
  + \sigma^2 K_{\nu}
      K_b
      \abs{ \langle
      \partial_{\xi}[g(\Phi + v_A) - g(\Phi + v_B) ],
        \psi_{\mathrm{tw}} \rangle_{L^2} }
\\[0.2cm]
& \le &
C_2 \sigma^2
 \norm{v_A-v_B}_{L^2}
 \big[ 1 + \norm{v_A}_{L^2} \big]
\\[0.2cm]
& & \qquad
+ C_2 \sigma^2
 \norm{v_A-v_B}_{L^2}
\big[ 1 + \norm{v_A}_{L^2} \big]
\\[0.2cm]
& & \qquad
 + C_2 \sigma^2
  \norm{v_A-v_B}_{L^2}.
\\[0.2cm]
\end{array}
\end{equation}
These expressions can be absorbed
into the bounds \sref{eq:swv:lip:bnd:e:ii}.
\end{proof}

\begin{cor}
\label{cor:swv:msigma:bnds}
Suppose that (Hf), (Hg)
and (hPar)
are satisfied.
Then there exists
$K > 0$ so that the following
holds true.
For any
$0 \le \sigma \le 1$
and any
$(v, d) \in H^1 \times \Real$
that has $\abs{d} \le 1$,
we have the estimate
\begin{equation}
\label{eq:swv:msigma:glb:bnd}
\begin{array}{lcl}
 \norm{\mathcal{M}_{\sigma;\Phi, c}(v, d)}_{L^2} & \le &
   K \big[1 + \norm{v}_{H^1} \big]
     \norm{v}_{H^1}  \norm{v}_{L^2}
\\[0.2cm]
& & \qquad
  +  K \big[
   \abs{c - c_0} + \norm{ \Phi- \Phi_0}_{H^1}
   + \abs{d}
  \big] \big[
    \norm{v}_{H^1} + \abs{d}
  \big]
\\[0.2cm]
& & \qquad
  + K \sigma^2
      \big[ 1 + \norm{v}_{H^1} \big] .
\\[0.2cm]
\end{array}
\end{equation}
In addition, for any $0 \le \sigma \le 1$
and any set of pairs
$(v_A, v_B) \in H^1 \times H^1$
and $(d_A, d_B) \in \Real^2$
for which $\abs{d_A} \le 1 $
and $\abs{d_B} \le 1$,
the expression
\begin{equation}
\Delta_{AB} \mathcal{M}_{\sigma;\Phi, c}
 = \mathcal{M}_{\sigma;\Phi, c}(v_A, d_A)
 - \mathcal{M}_{\sigma;\Phi, c}(v_B, d_B)
\end{equation}
satisfies the estimates
\begin{equation}
\label{eq:swv:msigma:lip:bnd}
\begin{array}{lcl}
\norm{\Delta_{AB} \mathcal{M}_{\sigma;\Phi, c} }_{L^2}
  & \le &
  K \big[ 1 + \norm{v_A}_{H^1} + \norm{v_B}_{H^1} \big]
    \big[ \norm{v_A}_{H^1}  + \norm{v_B}_{H^1} \big]
     \norm{v_A-v_B}_{L^2}
\\[0.2cm]
& & \qquad
+ K \big[
  \sigma^2 + \norm{v_A}_{H^1} + \abs{d_B}
   +  \norm{\Phi - \Phi_0}_{H^1}
  + \abs{c - c_0} \big]
\\[0.2cm]
& & \qquad \qquad \times
  \big[
     \norm{v_A-v_B}_{H^1}
      + \abs{d_A - d_B}
  \big]
\\[0.2cm]
& & \qquad
+ K \sigma^2
    \norm{v_A}_{H^1}^2 \norm{v_A}_{L^2}
 \norm{v_A-v_B}_{L^2}
\\[0.2cm]
& & \qquad
  + K \sigma^2
  \Big[ \norm{v_A}_{H^1} \norm{v_A}_{L^2} + \norm{v_B}_{H^1} \norm{v_B}_{L^2} \Big]
    \norm{v_A-v_B}_{H^1} ,
\\[0.2cm]
\abs{\langle \Delta_{AB} \mathcal{M}_{\sigma;\Phi, c}, \psi_{\mathrm{tw}} \rangle_{L^2} }
 & \le &
  K \big[ 1 + \norm{v_A}_{H^1} + \norm{v_B}_{H^1} \big]
    \big[ \norm{v_A}_{L^2}  + \norm{v_B}_{L^2} \big]
     \norm{v_A-v_B}_{L^2}
\\[0.2cm]
& & \qquad
+ K \big[
  \sigma^2 +
  \norm{v_A}_{L^2} + \abs{d_B}
   +  \norm{\Phi - \Phi_0}_{L^2}
  + \abs{c - c_0} \big]
\\[0.2cm]
& & \qquad \qquad \times
  \big[
     \norm{v_A-v_B}_{L^2}
      + \abs{d_A - d_B}
  \big]
\\[0.2cm]
& & \qquad
+ K \sigma^2
    \norm{v_A}_{H^1} \norm{v_A}_{L^2}^2
    \norm{v_A-v_B}_{L^2}
\\[0.2cm]
& & \qquad
  + K \sigma^2
  \Big[  \norm{v_A}_{L^2}^2 + \norm{v_B}_{L^2}^2 \Big]
    \norm{v_A-v_B}_{H^1} .
\\[0.2cm]
\end{array}
\end{equation}
\end{cor}
\begin{proof}
In view of the identity
\sref{eq:swv:id:for:m:sigma}
it suffices to note that
the terms
\sref{eq:swv:glb:bnd:nf},
\sref{eq:swv:mzero:glb:bnd},
\sref{eq:swv:ei:glb:bnds}
and \sref{eq:swv:glb:bnd:e:ii}
can be absorbed
in \sref{eq:swv:msigma:glb:bnd},
while the expressions
\sref{eq:swv:lip:bnd:nf},
\sref{eq:swv:mzero:lip:bnd},
\sref{eq:swv:ei:lip:bnds} and
\sref{eq:swv:lip:bnd:e:ii}
can be absorbed in \sref{eq:swv:msigma:lip:bnd}.
\end{proof}

\begin{cor}
\label{cor:swv:msigma:bnds:phi:zero}
Suppose that (Hf) and (Hg)
are satisfied.
Then there exists
$K > 0$ so that the following
holds true.
For any
$0 \le \sigma \le 1$
and any
$(v, d) \in H^1 \times \Real$
that has $\norm{v}_{H^1} \le 1$
together with $\abs{d} \le 1$,
we have the estimate
\begin{equation}
\label{eq:swv:msigma:glb:bnd:short}
\begin{array}{lcl}
 \norm{\mathcal{M}_{\sigma;\Phi_0, c_0}(v, d)}_{L^2} & \le &
    K  \big[ \norm{v}_{L^2} + \abs{d}   \big]
        \big[
    \norm{v}_{H^1} + \abs{d}
  \big]
  + K \sigma^2 .
\\[0.2cm]
\end{array}
\end{equation}
In addition, for any $0 \le \sigma \le 1$
and any set of pairs
$(v_A, v_B) \in H^1 \times H^1$
and $(d_A, d_B) \in \Real^2$
for which the bounds
\begin{equation}
\norm{v_{A}}_{H^1} \le 1,
\qquad
\abs{d_{A}} \le 1,
\qquad
\norm{v_{B}}_{H^1} \le 1,
\qquad
\abs{d_B} \le 1
\end{equation}
hold,
the expression
\begin{equation}
\Delta_{AB} \mathcal{M}_{\sigma;\Phi_0, c_0}
 = \mathcal{M}_{\sigma;\Phi_0, c_0}(v_A, d_A)
 - \mathcal{M}_{\sigma;\Phi_0, c_0}(v_B, d_B)
\end{equation}
satisfies the estimate
\begin{equation}
\label{eq:swv:msigma:lip:bnd:short}
\begin{array}{lcl}
\norm{\Delta_{AB} \mathcal{M}_{\sigma;\Phi_0, c_0} }_{L^2}
  & \le &
    K \big[
  \sigma^2 + \norm{v_A}_{H^1} + \norm{v_B}_{H^1} +  \abs{d_B}
    \big]
  \big[
     \norm{v_A-v_B}_{H^1}
      + \abs{d_A - d_B}
  \big] .
\\[0.2cm]
\end{array}
\end{equation}
\end{cor}
\begin{proof}
These bounds can easily be obtained by
simplifying the corresponding expressions from
Corollary \ref{cor:swv:msigma:bnds}.
\end{proof}

\subsection{Fixed-point constructions}
\label{sec:swv:fixp}

As a final preparation before
setting up our fixed-point problems,
we need to control the higher order effects
that arise when translating the adjoint eigenfunction $\psi_{\mathrm{tw}}$.
In particular, for any $\gamma \in \Real$
we introduce the function
\begin{equation}
\mathcal{N}_{\mathrm{tw}}(\gamma)
= T_{\gamma} \psi_{\mathrm{tw}} - \psi_{\mathrm{tw}} + \gamma \psi_{\mathrm{tw}}'
\end{equation}
and obtain the following bounds.
\begin{lem}
\label{lem:swv:shifts:ntw}
Suppose that (HTw) and (HS) hold. Then there exists $K > 0$
so that for any $\gamma \in \Real$ we have the bound
\begin{equation}
\norm{ \mathcal{N}_{\mathrm{tw}}(\gamma) }_{L^2} \le K \gamma^2,
\end{equation}
while for any pair $(\gamma_A, \gamma_B) \in \Real^2$
we have the estimate
\begin{equation}
\norm{ \mathcal{N}_{\mathrm{tw}}(\gamma_A)
  - \mathcal{N}_{\mathrm{tw}}(\gamma_B) }_{L^2} \le
   K \big[ \abs{\gamma_A } + \abs{\gamma_B} \big]
     \abs{\gamma_A - \gamma_B}.
\end{equation}
\end{lem}
\begin{proof}
In view of \sref{eq:var:diff:psi:tw},
we have the a priori bound
\begin{equation}
\norm{\mathcal{N}_{\mathrm{tw}}(\gamma)}_{L^2} \le C_1 \big[ 1 + \abs{\gamma} \big],
\end{equation}
together with
\begin{equation}
\norm{ \mathcal{N}_{\mathrm{tw}}(\gamma_A)
  - \mathcal{N}_{\mathrm{tw}}(\gamma_B) }_{L^2}
 \le C_1  \abs{\gamma_A-\gamma_B}.
\end{equation}
In particular, we can restrict our attention to the situation
where $\abs{\gamma} \le 1$ and $\abs{\gamma_A} + \abs{\gamma_B} \le 1$.
In this case we obtain the pointwise bounds
\begin{equation}
\abs{\mathcal{N}_{\mathrm{tw}}(\gamma)(\xi) } \le \frac{1}{2} \gamma^2
  \mathrm{sup}_{\xi - 1 \le \xi' \le \xi + 1} \abs{\psi''_{\mathrm{tw}}(\xi') }
\end{equation}
together with
\begin{equation}
\abs{\mathcal{N}_{\mathrm{tw}}(\gamma_A)(\xi)
  - \mathcal{N}_{\mathrm{tw}}(\gamma_B)(\xi) } \le
  \big[\mathrm{sup}_{\xi - 1 \le \xi' \le \xi + 1} \abs{\psi''_{\mathrm{tw}}(\xi') } \big]
  \big[ \frac{1}{2} (\gamma_A - \gamma_B)^2
    + \abs{\gamma_B} \abs{\gamma_A - \gamma_B}
  \big] .
\end{equation}
The desired bounds now follow from the exponential decay of $\psi''_{\mathrm{tw}}$.
\end{proof}

\begin{proof}[Proof of Proposition \ref{prp:mr:swv:ex}]
As a consequence of (HS), there exists a bounded
linear map
\begin{equation}
\mathcal{L}_{\mathrm{inv}}: L^2 \to H^2 \times \Real
\end{equation}
so that for any $h \in L^2$,
the pair $(v, d) = \mathcal{L}_{\mathrm{inv}} h$
is the unique solution in $H^2 \times \Real$
to the problem
\begin{equation}
\label{eq:swv:lin:op:def:eqn}
\mathcal{L}_{\mathrm{tw}} v
= h - \Phi_0' d .
\end{equation}
Indeed, we take $d = \langle h , \psi_{\mathrm{tw}} \rangle_{L^2}$,
which in view of the normalisation
\sref{eq:mr:hs:norm:cnd:psitw}
ensures that the right-hand side
of \sref{eq:swv:lin:op:def:eqn} is in the range
of $\mathcal{L}_{\mathrm{tw}}$.

It now suffices to find a solution
to the fixed-point problem
\begin{equation}
(v, d) = -\mathcal{L}_{\mathrm{inv}}
\mathcal{M}_{\sigma;\Phi_0, c_0}(v , d).
\end{equation}
Upon introducing the set
\begin{equation}
\mathcal{Z}_{\Theta} = \big\{ (v, d) \in H^2 \times \Real:
\norm{v}_{H^2} +  \abs{d} \le \min\{ 1 , \Theta \sigma^2 \} \big\}
 \subset H^2 \times \Real
\end{equation}
and applying Corollary \ref{cor:swv:msigma:bnds:phi:zero},
we see that for any $(v, d) \in \mathcal{Z}_{\Theta}$
we have
\begin{equation}
\norm{
  \mathcal{M}_{\sigma;\Phi_0, c_0}(v , d)
}_{L^2} \le
 K \big( \Theta^4 \sigma^4 + \sigma^2 \big)
 = K \sigma^2 \big( \Theta^2 \sigma^2 + 1 \big) ,
\end{equation}
while for any two pairs
$(v_A, d_A) \in \mathcal{Z}_{\Theta}$
and $(v_B , d_B) \in \mathcal{Z}_{\Theta}$
we have
\begin{equation}
\norm{
  \mathcal{M}_{\sigma;\Phi_0, c_0}(v_A , d_A)
  - \mathcal{M}_{\sigma;\Phi_0, c_0}(v_B , d_B)
}_{L^2} \le
K \sigma^2 \big[ 1 + 2 \Theta\big]
 \big[
   \norm{v_A - v_B}_{H^1}
   + \abs{d_A - d_B}
 \big].
\end{equation}
In particular, choosing
$\Theta$ to be sufficiently large
and $\delta_{\sigma} > 0$
to be sufficiently small,
we see that the map
$-\mathcal{L}_{\mathrm{inv}}
\mathcal{M}_{\sigma;\Phi_0, c_0}$
is a contraction on $\mathcal{Z}_{\Theta}$
for all $0 \le \sigma \le \delta_{\sigma}$.
\end{proof}

\begin{proof}[Proof of Proposition \ref{prp:mr:phase:shift}]
We first recall that
\begin{equation}
\langle \Phi_{\mathrm{ref}} , \psi'_{\mathrm{tw}} \rangle_{L^2}
= - \langle \Phi_{\mathrm{ref}}' , \psi_{\mathrm{tw}} \rangle_{L^2}
= - \langle \Phi_0' , \psi_{\mathrm{tw}} \rangle_{L^2}
 = -1 .
\end{equation}
Writing $u_0 = x_0 + \Phi_{\mathrm{ref}}$, this allows us to compute
\begin{equation}
\begin{array}{lcl}
\langle v_{\gamma} , \psi_{\mathrm{tw}} \rangle_{L^2}
& = &
    \langle x_0 + \Phi_{\mathrm{ref}} , T_{\gamma} \psi_{\mathrm{tw}} \rangle_{L^2}
  - \langle \Phi_{\sigma} , \psi_{\mathrm{tw}} \rangle_{L^2}
\\[0.2cm]
& = &
\langle x_0 + \Phi_{\mathrm{ref}} ,  \psi_{\mathrm{tw}} - \gamma \psi'_{\mathrm{tw}} + \mathcal{N}_{\mathrm{tw}}(\gamma) \rangle_{L^2}
  - \langle \Phi_{\sigma} , \psi_{\mathrm{tw}} \rangle_{L^2}
\\[0.2cm]
& = & \gamma + \langle x_0 + \Phi_{\mathrm{ref}} - \Phi_{\sigma} , \psi_{\mathrm{tw}} \rangle_{L^2} + \mathcal{E}_{\sigma}( x_0 , \gamma),
\end{array}
\end{equation}
in which we have
introduced the expression
\begin{equation}
\mathcal{E}_{\sigma}(x_0, \gamma)
 = - \gamma \langle x_0, \psi'_{\mathrm{tw}} \rangle_{L^2}
   + \langle x_0 + \Phi_{\mathrm{ref}} , \mathcal{N}_{\mathrm{tw}}(\gamma) \rangle_{L^2}.
\end{equation}

Using Lemma \ref{lem:swv:shifts:ntw},
we obtain the estimate
\begin{equation}
|\mathcal{E}_{\sigma}(x_0, \gamma)|
  \le C_1 \norm{x_0}_{L^2} \abs{ \gamma }
   + C_1 \big[ 1 + \norm{x_0}_{L^2} \big]
     \gamma^2 ,
\end{equation}
together with the Lipschitz bound
\begin{equation}
\begin{array}{lcl}
\norm{\mathcal{E}_{\sigma}(x_0, \gamma_A)
 -  \mathcal{E}_{\sigma}(x_0, \gamma_B) }_{L^2}
  & \le &
   C_2 \norm{x_0}_{L^2} \abs{ \gamma_A - \gamma_B }
\\[0.2cm]
& & \qquad
   + C_2 \big[ 1 + \norm{x_0}_{L^2} \big]
     \big[\abs{ \gamma_A } + \abs{\gamma_B} \big]
       \abs{\gamma_A - \gamma_B}.
\end{array}
\end{equation}
In particular, upon choosing $\delta_{\mathrm{fix}} > 0$ to be sufficiently small
and imposing the restriction
\begin{equation}
\label{eq:swv:phaseshift:small:cond}
\norm{x_0}_{L^2} +  \norm{x_0 + \Phi_{\mathrm{ref}} - \Phi_{\sigma} }_{L^2} <
  \delta_{\mathrm{fix}},
\end{equation}
we can define
$\gamma_0$ as the unique solution to the
fixed-point problem
\begin{equation}
- \gamma =
\langle x_0 + \Phi_{\mathrm{ref}} - \Phi_{\sigma} , \psi_{\mathrm{tw}} \rangle_{L^2}
 + \mathcal{E}_{\sigma}(x_0, \gamma)
\end{equation}
on the set
\begin{equation}
\Sigma_{x_0} = \{ \gamma: \abs{\gamma} \le 2 \norm{x_0 + \Phi_{\mathrm{ref}} - \Phi_{\sigma} }_{L^2}
  \norm{\psi_{\mathrm{tw}}}_{L^2} \}.
\end{equation}
By choosing $\delta_{\sigma} > 0$ and $\delta_0 > 0$
to be sufficiently small,
the bound \sref{eq:mr:prp:swv:bnd:phi:sigma}
allows us to conclude that \sref{eq:swv:phaseshift:small:cond}
is satisfied whenever
\sref{eq:swv:phaseshift:small:cond:u:only} holds.

For any $\gamma \in \Real$ we can compute
\begin{equation}
\begin{array}{lcl}
\norm{T_{-\gamma} \Phi_{\sigma} - \Phi_{\sigma}}^2_{L^2}
& = &
\int \big(  \Phi_{\sigma}(\xi + \gamma) - \Phi_{\sigma}(\xi) \big)^2 \, d \xi
\\[0.2cm]
& = &
\int \big[ \int_0^{\gamma} \Phi_{\sigma}'(\xi + s) \, ds \big]^2 \, d \xi
\\[0.2cm]
& \le &
\int \abs{\gamma} \int_0^{\gamma} \Phi_{\sigma}'(\xi + s)^2 \, ds \, d \xi
\\[0.2cm]
& = &
\abs{\gamma}^2  \int \Phi_{\sigma}'(\xi)^2 \, d \xi
\\[0.2cm]
& = & \abs{\gamma}^2  \norm{\Phi_{\sigma}'}^2_{L^2}.
\end{array}
\end{equation}
In particular, we obtain the bound
\begin{equation}
\begin{array}{lcl}
\norm{v_{\gamma_0}}_{L^2}
& = &
 \norm{ x_0 + \Phi_{\mathrm{ref}}
   - T_{\gamma_0} \Phi_{\sigma} }_{L^2}
\\[0.2cm]
& \le &
  \norm{ x_0 + \Phi_{\mathrm{ref}}
   -  \Phi_{\sigma} }_{L^2}
   + \norm{ T_{\gamma_0} \Phi_{\sigma}
     - \Phi_{\sigma} }_{L^2}
\\[0.2cm]
& \le &
  \norm{ x_0 + \Phi_{\mathrm{ref}}
   -  \Phi_{\sigma} }_{L^2}
  + C_3 \abs{ \gamma_0 }.
\end{array}
\end{equation}
The desired estimate \sref{eq:swv:bnds:on:shift:final}
hence follows from $\gamma_0 \in \Sigma_{x_0}$.
The final estimate \sref{eq:swv:bnds:on:shift:final:h1}
follows in a similar fashion, exploiting $\Phi_{\sigma}'' \in L^2$.
\end{proof}

\begin{proof}[Proof of Proposition \ref{prp:mr:expl:waves}]
For convenience, we introduce the notation
\begin{equation}
\alpha_{\sigma}
  =
  \big[ 1 + \frac{1}{2\rho} \sigma^2 \vartheta_{0}^2 \big]^{1/2}.
\end{equation}
Using the definitions \sref{eq:mr:id:for:mod:wave:exp:stb}
one easily verifies the identities
\begin{equation}
\Phi_{\sigma}'(\xi)
= \alpha_{\sigma}
  \Phi_0'(
    \alpha_{\sigma} \xi
  ),
\qquad \qquad
\Phi_{\sigma}''(\xi)
= \alpha_{\sigma}^2
  \Phi_0''(
    \alpha_{\sigma} \xi
  ),
\end{equation}
which yields
\begin{equation}
\begin{array}{lcl}
g\big(\Phi_{\sigma}(\xi)\big)
= g\big( \Phi_0( \alpha_{\sigma} \xi) \big)
= \vartheta_0 \Phi_0'\big(\alpha_{\sigma} \xi\big)
= \vartheta_0 \alpha_{\sigma}^{-1} \Phi'_{\sigma}(\xi),
\\[0.2cm]
\end{array}
\end{equation}
together with
\begin{equation}
f(\Phi_{\sigma} ) + c_{\sigma} \Phi_{\sigma}'
= - \alpha_{\sigma}^{-2} A_* \Phi_{\sigma}.
\end{equation}
Since
the cut-off functions in the definition of $b$
act as the identity for small $\sigma \ge 0$,
we obtain

\begin{equation}
\begin{array}{lcl}
b(\Phi_{\sigma}, \psi_{\mathrm{tw}})
  & = & - \vartheta_0 \alpha_{\sigma}^{-1},
\\[0.2cm]
\kappa_{\sigma}(\Phi_{\sigma}, \psi_{\mathrm{tw}} )
 & = & 1 + \frac{1}{2 \rho} \vartheta_{0}^2 \alpha_{\sigma}^{-2},
\end{array}
\end{equation}
which implies
\begin{equation}
\begin{array}{lcl}
\mathcal{J}_{\sigma}(\Phi_{\sigma}, c_{\sigma}, \psi_{\mathrm{tw}})
& = &
 [1 + \frac{1}{2\rho} \sigma^2 \vartheta_{0}^2 \alpha_{\sigma}^{-2} ]^{-1}
\big[ f( \Phi_{\sigma} ) + c_{\sigma} \Phi_{\sigma}'
  - \sigma^2 \vartheta_{0}^2 \alpha_{\sigma}^{-2}
      \Phi_{\sigma}'' \big]
\\[0.2cm]
& = &
-[1 + \frac{1}{2\rho} \sigma^2 \vartheta_{0}^2 \alpha_{\sigma}^{-2} ]^{-1}
\big[  \alpha_{\sigma}^{-2} A_* \Phi_\sigma
+ \frac{\sigma^2}{\rho} \vartheta_{0}^2 \alpha_{\sigma}^{-2}
      A_* \Phi_{\sigma} \big]
\\[0.2cm]
& = &
-[\alpha_{\sigma}^2
+ \frac{1}{2\rho} \sigma^2 \vartheta_{0}^2
]^{-1}
\big[  1 + \frac{\sigma^2}{\rho} \vartheta_0^2 \big]
   A_* \Phi_\sigma
\\[0.2cm]
& = & - A_* \Phi_{\sigma}.
\end{array}
\end{equation}
The claims now follow from the uniqueness
statement in Proposition \ref{prp:mr:swv:ex}.
\end{proof}

\section{Bounds on mild nonlinearities}
\label{sec:fnl}

In this section we set out to
obtain bounds on the
nonlinearities
$\overline{\mathcal{R}}_{\sigma;
 \Phi_{\sigma}, c_{\sigma}}$ and
 $\overline{\mathcal{S}}_{\sigma;\Phi_{\sigma}}$
defined in \sref{eq:md:overline:r:sigma}-\sref{eq:md:overline:s:sigma}.
In addition, we show that our choices
\sref{eq:mr:def:b:kappa} and \sref{eq:mr:def:a}
for $a_{\sigma}$ and $b$
prevent these nonlinearities from having
a component in the subspace of $L^2$
on which the semigroup $S(t)$ does not decay,
provided the cut-offs are not hit.

Our main result below shows that
the construction of $\Phi_{\sigma}$
has eliminated all $\O(1)$-terms
from the deterministic nonlinearity
$\overline{\mathcal{R}}$, leaving only
a small linear contribution together with the
expected higher order terms.
It is important to note here
that these higher order terms
depend at most quadratically on $\norm{v}_{H^1}$,
besides powers of $\norm{v}_{L^2}$.

In general, the stochastic nonlinearity
$\overline{\mathcal{S}}_{\sigma;\Phi_{\sigma}}$
will have an $\O(1)$-term, but we have
an explicit expression for this contribution so we also discuss the case when this contribution disappears.
In both cases, the higher order terms
depend at most linearly on $\norm{v}_{H^1}$.

\begin{prop}
\label{prp:fnl:bnds}
Consider the setting
of Proposition \ref{prp:mr:swv:ex}
and recall the definitions
\sref{eq:md:overline:r:sigma}
and \sref{eq:md:overline:s:sigma}.
Then there exists $K > 0$
so that for any
$0 \le \sigma \le \delta_{\sigma}$
and any
$v \in H^1$,
the following properties hold true.
\begin{itemize}
\item[(i)]{
We have the bound
\begin{equation}
\label{eq:fnl:rs:glb}
\begin{array}{lcl}
\norm{
\overline{\mathcal{R}}_{\sigma;
 \Phi_{\sigma}, c_{\sigma}}(v)
}_{L^2}
& \le &
 K \sigma^2 \norm{v}_{H^1}
+ K \norm{v}_{H^1}^2 \big[
  1
  + \norm{v}_{L^2}^2
  + \sigma^2 \norm{v}_{L^2}^3
\big].
\end{array}
\end{equation}
}
\item[(ii)]{
We have the estimate
\begin{equation}
\norm{
\overline{\mathcal{S}}_{\sigma;
 \Phi_{\sigma}}(v)
}_{L^2}
\le K\big[ 1 +  \norm{v}_{H^1}
  \big].
\end{equation}
}
\item[(iii)]{
  If the inequality
  \begin{equation}
    \norm{v}_{L^2}
      \le \min \{ 1 , [4 \norm{\psi_{\mathrm{tw}}}_{H^1}]^{-1} \}
  \end{equation}
  holds, then we have
  the identities
  \begin{equation}
    \langle \overline{\mathcal{R}}_{\sigma;
    \Phi_{\sigma}, c_{\sigma}}(v) , \psi_{\mathrm{tw}} \rangle_{L^2}
    =
    \langle \overline{\mathcal{S}}_{\sigma;
     \Phi_{\sigma}}(v) , \psi_{\mathrm{tw}} \rangle_{L^2} = 0.
  \end{equation}
}
\item[(iv)]{
  If the identity
  \begin{equation}
    \label{eq:fnl:cond:for:sp:case}
    g(\Phi_{\sigma}) = - b(\Phi_{\sigma}, \psi_{\mathrm{tw}} ) \Phi_{\sigma}'
  \end{equation}
  holds,
  then we have the bound
  \begin{equation}
     \norm{
        \overline{\mathcal{S}}_{\sigma;
       \Phi_{\sigma}}(v)
       }_{L^2}
        \le K  \norm{v}_{H^1} .
   \end{equation}
}
\end{itemize}
\end{prop}

In order to derive a compact expression
for $\overline{\mathcal{R}}_{\sigma;\Phi_{\sigma}, c_{\sigma}}$,
it is convenient to
recall the definition
\sref{eq:swv:def:m:sigma}
and
introduce the function
\begin{equation}
\overline{\mathcal{R}}_{\sigma;I}(v) =
\mathcal{M}_{\sigma;\Phi_{\sigma} , c_{\sigma} }
     (v, 0)
  - \mathcal{M}_{\sigma;\Phi_{\sigma} , c_{\sigma}}(0, 0).
\end{equation}
We note that the bounds in
Corollary \ref{cor:swv:msigma:bnds}
are directly applicable to this function.

\begin{lem}
Consider the setting
of Proposition \ref{prp:mr:swv:ex}.
Then for any $0 \le \sigma \le \delta_{\sigma}$
and $v \in H^1$,
we have the identity
\begin{equation}
\label{eq:fnl:id:for:r:sigma}
\overline{\mathcal{R}}_{\sigma;\Phi_{\sigma} , c_{\sigma}}(v)
= \overline{\mathcal{R}}_{\sigma;I}(v)
 - \Big[\chi_{\mathrm{low}}\big(\langle\partial_\xi[\Phi_{\sigma} + v ] ,
    \psi_{\mathrm{tw}} \rangle_{L^2} \big)\Big]^{-1}
   \langle
     \overline{\mathcal{R}}_{\sigma;I}(v),
     \psi_{\mathrm{tw}}
  \rangle_{L^2}
  [\Phi_{\sigma}' + v'].
\end{equation}
\end{lem}
\begin{proof}
Inspecting \sref{eq:swv:def:m:sigma}
and using the defining
property
\sref{eq:mr:prop:swv:eq}
for $(\Phi_{\sigma} , c_{\sigma})$,
we see that
\begin{equation}
-\mathcal{M}_{\sigma;\Phi_{\sigma} ,
   c_{\sigma}}(0, 0)
= A_* \Phi_{\sigma}
  + \mathcal{J}_0
    (\Phi_{\sigma}, c_\sigma ).
\end{equation}

Applying \sref{eq:swv:def:m:sigma}
once more,
we hence find
\begin{equation}
\begin{array}{lcl}
\mathcal{J}_{\sigma}
  ( \Phi_{\sigma} + v, c_{\sigma} , \psi_{\mathrm{tw}})
& = & \mathcal{J}_0
  (\Phi_\sigma, c_{\sigma})
+ [\mathcal{L}_{\mathrm{tw}} - A_*] v
+ \mathcal{M}_{\sigma;\Phi_\sigma , c_\sigma }
     (v, 0)
\\[0.2cm]
& = &
  [\mathcal{L}_{\mathrm{tw}} - A_* ] v
  -A_* \Phi_{\sigma}
  + \mathcal{M}_{\sigma;\Phi_{\sigma} ,
    c_{\sigma} }
     (v, 0)
  - \mathcal{M}_{\sigma;\Phi_{\sigma} ,
     c_{\sigma}}(0, 0)
\\[0.2cm]
& = &
 [\mathcal{L}_{\mathrm{tw}} - A_* ] v
  -A_* \Phi_{\sigma}
  + \overline{\mathcal{R}}_{\sigma;I}(v).
\end{array}
\end{equation}
Writing
\begin{equation}
\begin{array}{lcl}
\mathcal{I}_{\sigma}(v)
& = &
\langle \Phi_{\sigma} + v, A_* \psi_{\mathrm{tw}} \rangle_{L^2}
+ \langle \mathcal{J}_{\sigma}
  (\Phi_{\sigma} + v, c_{\sigma} , \psi_{\mathrm{tw}}),
 \psi_{\mathrm{tw}} \rangle_{L^2}
\end{array}
\end{equation}
and using $\mathcal{L}_{\mathrm{tw}}^{\mathrm{adj}} \psi_{\mathrm{tw}} = 0$,
we may compute
\begin{equation}
\begin{array}{lcl}
\mathcal{I}_{\sigma}(v)
&  = &
\langle \Phi_{\sigma} , A_* \psi_{\mathrm{tw}} \rangle_{L^2}
+
\langle v,
  [A_* - \mathcal{L}_{\mathrm{tw}}^{\mathrm{adj}} ] \psi_{\mathrm{tw}} \rangle_{L^2}
\\[0.2cm]
& & \qquad
+ \langle \mathcal{J}_{\sigma}
  (\Phi_{\sigma} + v, c_{\sigma} , \psi_{\mathrm{tw}}),
 \psi_{\mathrm{tw}} \rangle_{L^2}
\\[0.2cm]
& = &
\langle A_* \Phi_{\sigma}, \psi_{\mathrm{tw}} \rangle_{L^2}
+ \langle [A_* - \mathcal{L}_{\mathrm{tw}} ] v , \psi_{\mathrm{tw}}
  \rangle_{L^2}
\\[0.2cm]
& & \qquad
+ \langle \mathcal{J}_{\sigma}
  (\Phi_{\sigma} + v, c_{\sigma} , \psi_{\mathrm{tw}}),
 \psi_{\mathrm{tw}} \rangle_{L^2}
\\[0.2cm]
& = &
  \langle \overline{\mathcal{R}}_{\sigma;I}(v),
   \psi_{\mathrm{tw}} \rangle_{L^2} .
\end{array}
\end{equation}
In view
of the definition
\sref{eq:mr:def:a}
for $a_{\sigma}$,
we now obtain
\begin{equation}
\begin{array}{lcl}
\kappa_{\sigma}(\Phi_{\sigma}+ v, \psi_{\mathrm{tw}})^{-1}
a_{\sigma}\big(\Phi_{\sigma} + v ,
c_{\sigma}, \psi_{\mathrm{tw}} \big)
& = &
  - \Big[\chi_{\mathrm{low}}\big(\langle\partial_\xi[\Phi_{\sigma} + v ] ,
    \psi_{\mathrm{tw}} \rangle_{L^2} \big)\Big]^{-1}
  \mathcal{I}_{\sigma}(v)
\\[0.2cm]
& = &
  - \Big[\chi_{\mathrm{low}}\big(\langle\partial_\xi[\Phi_{\sigma} + v ] ,
    \psi_{\mathrm{tw}} \rangle_{L^2} \big)\Big]^{-1}
   \langle
      \overline{\mathcal{R}}_{\sigma;I}(v),
          \psi_{\mathrm{tw}}
   \rangle_{L^2} .
\end{array}
\end{equation}
In particular, the desired
identity
\sref{eq:fnl:id:for:r:sigma}
follows directly
from the definition
\sref{eq:md:overline:r:sigma}.
\end{proof}

\begin{lem}
\label{lem:fnl:rsi}
Consider the setting
of Proposition \ref{prp:mr:swv:ex}.
Then there exists $K > 0$
so that for any $v \in H^1$
and $0 \le \sigma \le \delta_{\sigma}$
we have the bound
\begin{equation}
\label{eq:fnl:rsi:glb}
\begin{array}{lcl}
\norm{\overline{\mathcal{R}}_{\sigma;I}(v)}_{L^2}
& \le &
K \sigma^2 \norm{v}_{H^1}
+ K \norm{v}_{H^1}^2
  \big[
    1 + \norm{v}_{L^2} + \sigma^2 \norm{v}_{L^2}^2
  \big] ,
\\[0.2cm]
\end{array}
\end{equation}
together with
\begin{equation}
\label{eq:fnl:rsi:glb:ip}
\begin{array}{lcl}
\abs{ \langle \overline{\mathcal{R}}_{\sigma;I}(v)  , \psi_{\mathrm{tw}} \rangle_{L^2} }
& \le &
K  \norm{v}_{L^2}
  \big[ \sigma^2 + \norm{v}_{L^2} \big]
+ K \norm{v}_{H^1}
  \big[
     \norm{v}_{L^2}^2 + \sigma^2 \norm{v}_{L^2}^3
  \big]    .
\end{array}
\end{equation}
\end{lem}
\begin{proof}
Applying Corollary
\ref{cor:swv:msigma:bnds},
we find
\begin{equation}
\begin{array}{lcl}
\norm{\overline{\mathcal{R}}_{\sigma;I}(v)}_{L^2}
& \le &
C_1 \big[ 1 + \norm{v}_{H^1}  \big]
    \norm{v}_{H^1}
     \norm{v}_{L^2}
\\[0.2cm]
& & \qquad
+ C_1 \big[
  \sigma^2 + \norm{v}_{H^1}
  \big] \norm{v}_{H^1}
\\[0.2cm]
& & \qquad
+ C_1 \sigma^2
    \norm{v}_{H^1}^2 \norm{v}_{L^2}
 \norm{v}_{L^2}
\\[0.2cm]
& & \qquad
  + C_1 \sigma^2
   \norm{v}_{H^1} \norm{v}_{L^2}
     \norm{v}_{H^1} ,
\\[0.2cm]
\end{array}
\end{equation}
together with
\begin{equation}
\begin{array}{lcl}
\abs{ \langle \overline{\mathcal{R}}_{\sigma;I}(v)  , \psi_{\mathrm{tw}} \rangle_{L^2} }
& \le &
C_2 \big[ 1 + \norm{v}_{H^1}  \big]
     \norm{v}_{L^2}
     \norm{v}_{L^2}
\\[0.2cm]
& & \qquad
+ C_2 \big[
  \sigma^2 +
  \norm{v}_{L^2}
    \big]
     \norm{v}_{L^2}
\\[0.2cm]
& & \qquad
+ C_2 \sigma^2
    \norm{v}_{H^1} \norm{v}_{L^2}^2
    \norm{v}_{L^2}
\\[0.2cm]
& & \qquad
  + C_2 \sigma^2
    \norm{v}_{L^2}^2
    \norm{v}_{H^1}.
\\[0.2cm]
\end{array}
\end{equation}
These expressions can be absorbed into
\sref{eq:fnl:rsi:glb} and \sref{eq:fnl:rsi:glb:ip}.
\end{proof}

\begin{lem}
\label{lem:fnl:lip:s:sigma}
Consider the setting
of Proposition \ref{prp:mr:swv:ex}.
Then there exists $K > 0$
so that for any
$0 \le \sigma \le \delta_{\sigma}$
and any
$v \in H^1$
we have the bound
\begin{equation}
\norm{
  \overline{\mathcal{S}}_{\sigma;\Phi_{\sigma}}(v)
  -
  \overline{\mathcal{S}}_{\sigma;\Phi_{\sigma}}(0)
}_{L^2}
\le
K \norm{v}_{H^1}.
\end{equation}
\end{lem}
\begin{proof}
Writing
\begin{equation}
\begin{array}{lcl}
\mathcal{I}
 & = & \overline{\mathcal{S}}_{\sigma;\Phi_{\sigma}}(v)
  -
  \overline{\mathcal{S}}_{\sigma;\Phi_{\sigma}}(0)
\\[0.2cm]
 & = &
    \kappa_{\sigma}(\Phi_{\sigma} + v , \psi_{\mathrm{tw}})^{-1/2}
  \Big[
     g( \Phi_{\sigma} + v )
     + b( \Phi_{\sigma} + v , \psi_{\mathrm{tw}})
        \partial_\xi [ \Phi_{\sigma} + v ]
  \Big]
\\[0.2cm]
& & \qquad
  -  \kappa_{\sigma}(\Phi_{\sigma}  , \psi_{\mathrm{tw}})^{-1/2}
  \Big[
     g( \Phi_{\sigma}  )
     + b( \Phi_{\sigma}  , \psi_{\mathrm{tw}})
        \partial_\xi [ \Phi_{\sigma}  ]
  \Big]
\\[0.2cm]

\end{array}
\end{equation}
and using
Lemma's
\ref{lem:prlm:ests:g},
\ref{lem:prlm:bnds:b},
\ref{lem:prlm:bnds:kappa}
and \ref{lem:prlm:nu},
we compute
\begin{equation}
\begin{array}{lcl}
\norm{\mathcal{I}
}_{L^2}
& \le &
\abs{\nu_{\sigma}^{(-1/2)}
  (\Phi_{\sigma} + v , \psi_{\mathrm{tw}})
- \nu_{\sigma}^{(-1/2)}(\Phi_\sigma,
   \psi_{\mathrm{tw}} )}
\big[
  \norm{g(\Phi_{\sigma} ) }_{L^2}
  + K_b \norm{\Phi_{\sigma}'  }_{L^2}
\big]
\\[0.2cm]
& & \qquad
+ K_{\kappa}
  \norm{g(\Phi_{\sigma} + v) -
    g(\Phi_\sigma) }_{L^2}
\\[0.2cm]
& & \qquad
+ K_{\kappa}
  \abs{b(\Phi_{\sigma} + v, \psi_{\mathrm{tw}})
     - b(\Phi_{\sigma} , \psi_{\mathrm{tw}}) }
  \norm{ \Phi_{\sigma}' }_{L^2}
\\[0.2cm]
& & \qquad
+ K_{\kappa} K_b \norm{v'}_{L^2}.
\end{array}
\end{equation}
Applying these results once more,
we find
\begin{equation}
\begin{array}{lcl}
\norm{\mathcal{I}
}_{L^2}
& \le &
C_1 \sigma^2
\norm{v}_{L^2}
+ C_1
  \norm{v}_{L^2}
+ C_1
    \norm{v}_{L^2}
+ C_1 \norm{v}_{H^1}
\\[0.2cm]
& \le &
C_2 \norm{v}_{H^1} ,
\end{array}
\end{equation}
as desired.
\end{proof}

\begin{proof}[
Proof of Proposition \ref{prp:fnl:bnds}]
To obtain (i),
we use \sref{eq:fnl:id:for:r:sigma}
together with
Lemma \ref{lem:fnl:rsi}
to compute
\begin{equation}
\begin{array}{lcl}
\norm{\overline{\mathcal{R}}_{\sigma;\Phi_{\sigma}, c_{\sigma}}(v)}_{L^2}
& \le &
\norm{\overline{\mathcal{R}}_{\sigma;I}(v)}_{L^2}
+ C_1 \abs{
  \langle \overline{\mathcal{R}}_{\sigma;I}(v),
    \psi_{\mathrm{tw}} \rangle_{L^2} }
    \big[ 1 + \norm{v}_{H^1} \big]
\\[0.2cm]
& \le &
C_2 \sigma^2 \norm{v}_{H^1}
+ C_2 \norm{v}_{H^1}^2
  \big[
    1 + \norm{v}_{L^2} + \sigma^2 \norm{v}_{L^2}^2
  \big]
\\[0.2cm]
& & \qquad
+C_2 \norm{v}_{L^2} \big[ \sigma^2 + \norm{v}_{L^2} \big]
    \big[ 1 + \norm{v}_{H^1} \big]
\\[0.2cm]
& & \qquad
+ C_2 \norm{v}_{H^1}
  \big[
     \norm{v}_{L^2}^2 + \sigma^2 \norm{v}_{L^2}^3
  \big]
     \big[ 1 + \norm{v}_{H^1} \big].
\\[0.2cm]
\end{array}
\end{equation}
These terms can all be absorbed into
\sref{eq:fnl:rs:glb}.

The bound (ii) follows directly
from Lemma \ref{lem:fnl:lip:s:sigma},
using the estimate
\begin{equation}
\norm{
  \overline{\mathcal{S}}_{\sigma;\Phi_{\sigma}}(v)
}_{L^2}
\le
 \norm{
   \overline{\mathcal{S}}_{\sigma;\Phi_{\sigma}}(0)
 }_{L^2}
 + \norm{
   \overline{\mathcal{S}}_{\sigma;\Phi_{\sigma}}(v)
  -
  \overline{\mathcal{S}}_{\sigma;\Phi_{\sigma}}(0)
 }_{L^2}
\end{equation}
and the a-priori bound
\begin{equation}
\norm{\overline{\mathcal{S}}_{\sigma;\Phi_{\sigma}}
  ( 0 ) }_{L^2}
\le C_3.
\end{equation}
The bound (iv) follows in the same fashion,
since
the condition
 \sref{eq:fnl:cond:for:sp:case}
implies that
\begin{equation}
\overline{\mathcal{S}}_{\sigma;\Phi_{\sigma}}(0)
  = 0 .
\end{equation}

Finally, (iii) follows from
the identities
\sref{eq:fnl:id:for:r:sigma}
and \sref{eq:prlm:id:for:b:without:cutoffs},
using the proof of
Lemma \ref{lem:prlm:b:without:cutoffs}
to show that the cut-off function
$\chi_{\mathrm{low}}$
in \sref{eq:fnl:id:for:r:sigma} acts
as the identity.
\end{proof}

\section{Nonlinear stability of mild solutions}
\label{sec:nls}

In this section we
prove Theorems \ref{thm:mr:orbital:stb} and \ref{thm:mr:exp:stb},
providing an  orbital and an exponential stability result
for the stochastic waves $(\Phi_{\sigma}, c_{\sigma})$ on timescales
of order $\sigma^{-2}$.
Recalling the function
\sref{eq:md:def:n},
our key statement is that
 $E \sup_{t} \overline{N}_{\e,\alpha}(t)$
can be bounded in terms of itself, the noise-strength $\sigma$
and the initial condition $\norm{\overline{V}(0)}_{H^1}^2$.
This requires a number of technical regularity estimates,
which we obtain in \S\ref{sec:nls:reg:ests:det}-\ref{sec:nls:reg:ests:st}.

In order to prevent cumbersome notation and to
highlight the broad applicability of our techniques here,
we do not refer to the specific functions $\overline{V}$
and the specific nonlinearities
$\overline{\mathcal{R}}_{\sigma;\Phi_{\sigma}, c_{\sigma}}$
here. Instead, we assume the
following general condition
concerning the form of our nonlinearities.
\begin{itemize}
\item[(hFB)]
We have $\nrm{B_{\mathrm{cn}}}_{L^2}=K_{B;\mathrm{cn}}<\infty$ and the maps
\begin{equation}
F_{\mathrm{lin}} : H^1 \to L^2,
\qquad
F_{\mathrm{nl}} : H^1 \to L^2,
\qquad
B_{\mathrm{lin}} : H^1 \to L^2
\end{equation}
satisfy the bounds
\begin{equation}
\begin{array}{lcl}
\nrm{F_{\mathrm{lin}}(v)}_{L^2}&\leq& K_{F;\mathrm{lin}}\nrm{v}_{H^1},
  \\[0.2cm]
\nrm{F_{\mathrm{nl}}(v)}_{L^2}&\leq & K_{F;\mathrm{nl}}\nrm{v}^2_{H^1}(1+\nrm{v}_{L^2}^{m}),
  \\[0.2cm]
\nrm{B_{\mathrm{lin}}(v)}_{L^2}&\leq&  K_{B;\mathrm{lin}}\nrm{v}_{H^1}
\end{array}
\end{equation}
for some $m > 0$.
In addition, there exists $\eta_0 > 0$
so that
\begin{equation}
\label{eq:nls:projs:are:zero}
\langle \sigma^2 F_{\mathrm{lin}}(v) + F_{\mathrm{nl}}(v) , \psi_{\mathrm{tw}} \rangle_{L^2}
 = 0,
\qquad
\langle  B_{\mathrm{cn}} + B_{\mathrm{lin}}(v) , \psi_{\mathrm{tw}} \rangle_{L^2}
 = 0
\end{equation}
whenever
$\nrm{v}_{L^2}\leq \eta_0$.
\end{itemize}

Using the nonlinearities above,
we can discuss the mild formulation
of the SPDE that we are interested in.
At present, we simply assume that a solution is a priori
available, but one can also set out to construct
such a solution directly.
\begin{itemize}
\item[(hSol)]
For any $T > 0$,
there exists a continuous $(\mathcal{F}_t)$-adapted process $V: \Omega \times [0, T] \to L^2$
for which we have the inclusions
\begin{equation}
V \in \mathcal{N}^2\big([0,T]; (\mathcal{F}_t) ; H^1\big),
\qquad
B_{\mathrm{lin}}(V) \in \mathcal{N}^2\big([0,T]; (\mathcal{F}_t) ; L^2\big).
\end{equation}
In addition, for almost all $\omega \in \Omega$
we have the inclusions
\begin{equation}
F_{\mathrm{lin}}\big(V(\cdot, \omega)\big) \in L^1([0,T];L^2),
\qquad
F_{\mathrm{nl}}\big(V(\cdot, \omega)\big) \in L^1([0,T];L^2)
\end{equation}
together with the identity
\begin{equation}
\label{eq:nls:hsol:id:for:v}
\begin{array}{lcl}
V(t) & = &S(t)V(0)+\sigma^2 \int_0^t S(t-s) F_{\mathrm{lin}}\big(V(s)\big) \, ds
                + \int_0^t S(t-s) F_{\mathrm{nl}}\big(V(s)\big) \, ds
\\[0.2cm]
& & \qquad
+\sigma \int_0^t S(t-s)B_{\mathrm{cn}} \, d\b_s
+ \sigma \int_0^t S(t-s)B_{\mathrm{lin}}\big(V(s)\big) \, d\b_s ,
\end{array}
\end{equation}
which holds for all $t \in [0, T]$.
Finally, we have $\langle V(0), \psi_{\mathrm{tw}} \rangle_{L^2} = 0$.
\end{itemize}

For any $\e > 0$
and $\alpha \ge 0$,
we recall the notation
\begin{equation}
\label{eq:nls:recdefN}
\begin{array}{lcl}
 N_{\e,\a}(t)=
 e^{\a t}\nrm{V(t)}_{L^2}^2
 +\int_0^te^{-\e(t-s)}e^{\a s}\nrm{V(s)}^2_{H^1}ds.
\end{array}
\end{equation}
For any $T > 0$ and $\eta > 0$,
we introduce the
$(\mathcal{F}_t)$-stopping time
\begin{equation}
\tau_{\e,\alpha}(T,\eta)
 = \inf\Big\{0 \leq t < T:
     N_{\e,\a}(t)
     > \eta
  \Big\} ,
\end{equation}
writing $\tau_{\e,\alpha}(T,\eta) = T$
if the set is empty.
Our two main results here, which we establish
in \S\ref{sec:nls:reg:ests:st}
provide 
bounds
on the expectation of $\sup_{0\leq t\leq \tau_{\e,\a}(T,\eta)}N_{\e,\alpha}(t)$.

\begin{prop}
\label{prp:nls:general}
Assume that (HA), (HTw), (HS), (H$\beta$),
(hSol) and (hFB) are satisfied.
Pick a constant $0 < \e < \beta$,
together with two sufficiently small
constants $\delta_{\eta} > 0$
and $\delta_{\sigma} > 0$.
Then there exists a constant
$K > 0$ so that for any $T > 1$,
any $0 < \eta \le \delta_{\eta}$
and any $0 \le \sigma \le \delta_{\sigma}  T^{-1/2}$
we have the bound
\begin{equation}
\label{eq:nls:prp:general:estimate}
\begin{array}{lcl}
E \sup_{0 \le t \le \tau_{\e,0}(T, \eta) } N_{\e,0} (t)
&\leq &
   K \Big[ 
     \nrm{V(0)}^2_{H^1}+ \sigma^2 T \Big] .
\end{array}
\end{equation}
\end{prop}

\begin{prop}
\label{prp:nls:exponential}
Assume that (HA), (HTw), (HS), (H$\beta$),
(hSol) and (hFB) are satisfied
and that $B_{\mathrm{cn}} = 0$.
Pick two constants $\e > 0$,
$\alpha \ge 0$ for which $\e + \frac{\alpha}{2} < \beta$,
together with two sufficiently small
constants $\delta_{\eta} > 0$
and $\delta_{\sigma} > 0$.
Then there exists a constant
$K > 0$ so that for any $T > 1$,
any $0 < \eta \le \delta_{\eta}$
and any $0 \le \sigma \le \delta_{\sigma} T^{-1/2}$
we have the bound
\begin{equation}
\label{eq:nls:prp:exponential:estimate}
\begin{array}{lcl}
E \sup_{0 \le t \le \tau_{\e,\a}(T,\eta) }
   N_{\e,\alpha}\big(  t\big)
&\leq &
   K 
   \nrm{V(0)}^2_{H^1} .
\end{array}
\end{equation}
\end{prop}

Exploiting the technique used in
Stannat \cite{Stannat},
these bounds can be turned into estimates
concerning the probabilities
\begin{equation}
p_{\e,\alpha}(T,\eta) = P\Big(
 \sup_{0 \leq t \leq T} \big[N_{\e,\a}(t)\big]
 > \eta
\Big).
\end{equation}
This allows our
main stability theorems to be established.

\begin{cor}
\label{cor:nls:general}
Consider the setting of
Proposition \ref{prp:nls:general}.
Then there exists a constant $K > 0$
so that for any $T > 1$,
any $0 < \eta \le \delta_{\eta}$
and any $0 \le \sigma \le \delta_{\sigma} T^{-1/2}$,
we have the bound
\begin{equation}
p_{\e,0}(T,\eta)\leq \eta^{-1} K \big[ \nrm{V(0)}_{H^1}^2+ \sigma^2 T \big] .
\end{equation}
\end{cor}
\begin{proof}
Upon computing
\begin{equation}
\begin{array}{lcl}
\eta p_{\e,0}(T,\eta)
& = & \eta P\big( \tau_{\e,0}(T, \eta) < T\big)
\\[0.2cm]
& = & E \Big[
      \mathbf{1}_{\tau_{\e,0}(T, \eta) < T}
      N_{\e,0}\big(\tau_{\e,0}(T, \eta)\big)
    \Big]
\\[0.2cm]
& \le & E N_{\e,0}\big( \tau_{\e,0}(T, \eta)  \big)
\\[0.2cm]
& \le & E \sup_{0 \le t \le \tau_{\e,0}(T, \eta) } N_{\e,0 } (t)
,
\end{array}
\end{equation}
the result follows from
\sref{eq:nls:prp:general:estimate}.
\end{proof}

\begin{cor}
\label{cor:nls:exponential}
Consider the setting of
Proposition \ref{prp:nls:exponential}.
Then there exists a constant $K > 0$
so that for any $T > 1$,
any $0 < \eta \le \delta_{\eta}$
and any $0 \le \sigma \le \delta_{\sigma} T^{-1/2}$
we have the bound
\begin{equation}
p_{\e,\alpha}(T,\eta)\leq \eta^{-1} K \nrm{V(0)}_{H^1}^2.
\end{equation}
\end{cor}
\begin{proof}
Upon computing
\begin{equation}
\begin{array}{lcl}
\eta p_{\e,\alpha}(T,\eta)
& = & \eta P\big( \tau_{\e,\alpha}(T, \eta) < T\big)
\\[0.2cm]
& = & E \Big[
      \mathbf{1}_{\tau_{\e,\alpha}(T, \eta) < T}
      N_{\e,\a}\big(\tau_{\e,\alpha}(T, \eta) \big)
    \Big]
\\[0.2cm]
& \le & E N_{\e,\a}\big(\tau_{\e,\alpha}(T, \eta)  \big)
\\[0.2cm]
& \le & E \sup_{0 \le t \le \tau_{\e,\a}(T, \eta)}
   N_{\e,\a}(t)
,
\end{array}
\end{equation}
the result follows from
\sref{eq:nls:prp:exponential:estimate}.
\end{proof}

\begin{proof}[Proof of Theorems \ref{thm:mr:orbital:stb} and \ref{thm:mr:exp:stb}]
On account of Propositions \ref{prp:mr:phase:shift}
and
\ref{prp:md:props:overline:v},
the map $\overline{V}$
defined in \sref{eq:mld:def:ovl:v}
satisfies the conditions of (hSol)
with $\big(\overline{\beta}_\tau , \overline{\mathcal{F}}_\tau \big)_{\tau \ge 0}$
as the relevant Brownian motion.  In addition, Proposition
\ref{prp:fnl:bnds} guarantees that (hFB) is satisfied.
The desired estimates now follow
from Corollaries \ref{cor:nls:general} and \ref{cor:nls:exponential},
using Proposition \ref{prp:md:prob:after:time:transf}
to reverse the time-transform.
\end{proof}

\subsection{Setup}
\label{sec:nls:setup}
In order to
establish Propositions
\ref{prp:nls:general}-\ref{prp:nls:exponential}
we need to estimate each of the terms
featuring in the identity
\sref{eq:nls:hsol:id:for:v}.
The regularity structure of the
semigroup $S(t)$ is crucial for our
purposes here, so we discuss this in some detail
using the terminology used in \cite[{\S}10]{hytonen2018analysis}.

In particular, for any $0 < \varphi < \pi$
we introduce the sector
\begin{equation}
\Sigma_{\varphi} = \{z \in \mathbb{C} \setminus \{ 0 \} :
\abs{ \mathrm{arg}(z) }  < \varphi \},
\end{equation}
in which we take $\mathrm{arg}(z) \in (-\pi, \pi)$.
We recall that a linear operator
$\mathcal{L}: D(\mathcal{L}) \subset X \to X$
on a Banach space $X$
is called sectorial if the spectrum of $\mathcal{L}$
is contained in $\overline{\Sigma}_{\omega}$
for some $0 < \omega(\mathcal{L}) < \frac{\pi}{2}$,
while the resolvent operators $R(z, \mathcal{L}) = (z - \mathcal{L})^{-1}$
satisfy the bound
\begin{equation}
\sup_{z \in \mathbb{C} \setminus  \overline{\Sigma}_{\omega(\mathcal{L})} }
  \norm{z R(z, \mathcal{L}) }_{\mathcal{L}(X,X)}
   < \infty .
\end{equation}

Our spectral assumptions (HS) combined
with the fact that $\mathcal{L}_{\mathrm{tw}}$
is a lower-order perturbation to
the diffusion operator $A_*$
guarantee that $-\mathcal{L}_{\mathrm{tw}}$
is sectorial. This means that
$\mathcal{L}_{\mathrm{tw}}$ generates
an analytic semigroup. In order
to isolate the behaviour caused by the neutral eigenmode,
we introduce the map
$Q: L^2 \to L^2$ that acts as
\begin{equation}
Q v = v - \langle v, \psi_{\mathrm{tw}} \rangle_{L^2} \Phi_0' .
\end{equation}
This projection allows us to formulate
several important estimates.
\begin{lem}[see \cite{lorenzi2004analytic}]
\label{lem:nls:sem:group:decay}
Assume that (HTw) and (HS) hold and consider
the analytic  semigroup $S(t)$
generated by $\mathcal{L}_{\mathrm{tw}}$.
Then there is a constant $M\geq 1$ for which
we have the bounds
\begin{equation}
\begin{array}{lcl}
\nrm{S(t)Q}_{\L(L^2,L^2)}&\leq& M e^{-\b t}, \hspace{1cm} 0<t<\infty ,
  \\[0.2cm]
\nrm{S(t)Q}_{\L(L^2,H^1)}&\leq& M t^{-\frac{1}{2}}, \hspace{1.1cm} 0<t\leq 2 ,
\\[0.2cm]
\nrm{S(t)Q}_{\L(L^2,H^1)}&\leq& M e^{-\b t}, \hspace{1.1cm} t\geq 1,
\\[0.2cm]
\norm{[\mathcal{L}_\mathrm{tw} - A_*] S(t)Q}_{\L(L^2,L^2) }
  & \le & M t^{-\frac{1}{2}}, \hspace{1.1cm} 0 < t\leq 2,
\\[0.2cm]
\norm{[\mathcal{L}^{\mathrm{adj}}_\mathrm{tw} - A_*] S(t)Q }_{\L(L^2,L^2) }
  & \le & M t^{-\frac{1}{2}}, \hspace{1.1cm} 0 < t\leq 2 .
\end{array}
\end{equation}
\end{lem}

In order to understand the combination $S(t) Q$
as an independent semigroup,
we introduce the spaces
\begin{equation}
L^2_Q = \{ v \in L^2: (I - Q) v = 0 \},
\qquad
H^2_Q = \{ v \in H^2: (I-Q) v = 0 \}
\end{equation}
and consider the operator $\mathcal{L}_{\mathrm{tw}}^Q: H^2_Q \to L^2_Q$
that arises upon restricting $\mathcal{L}_{\mathrm{tw}}$
to act on $H^2_Q$. Note that this is well-defined
since $\mathrm{Range} \big( \mathcal{L}_{\mathrm{tw}} \big)
   = L^2_Q$.
For any $\theta \in \Real$,
we now introduce the linear operators
\begin{equation}
B_{\theta} = -\big[ \mathcal{L}_{\mathrm{tw}} + \theta \big],
\qquad
B^Q_{\theta} = -\big[ \mathcal{L}_{\mathrm{tw}}^Q + \theta \big].
\end{equation}

\begin{lem}
\label{lem:nls:b:q:theta:sect}
Assume that (HTw) and (HS) hold
and pick any $0 \le \theta \le \beta$.
Then
the operator $B^Q_{\theta}$
is sectorial on $L^2_Q$ and the semigroup
generated by $-B^Q_{\theta}$ corresponds
with the restriction of $e^{\theta t} S(t)$
to $L^2_Q$.
\end{lem}
\begin{proof}
Note first that $\mathcal{L}_{\mathrm{tw}}^Q$ is bijective
since we have projected out the one-dimensional kernel.
For any $v \in L^2_Q$ and $\lambda$
in the resolvent set of $\mathcal{L}_{\mathrm{tw}}$,
we may compute
\begin{equation}
\begin{array}{lcl}
0 & = & (I - Q) \mathcal{L}_{\mathrm{tw}} R(\lambda, \mathcal{L}_{\mathrm{tw}}) v
\\[0.2cm]
 & = & (I-Q) \big[ -v + \lambda R(\lambda, \mathcal{L}_{\mathrm{tw}}) v \big]
\\[0.2cm]
 & = & \lambda (I - Q) R(\lambda, \mathcal{L}_{\mathrm{tw}}) v.
\end{array}
\end{equation}
which implies that
$R(\lambda, \mathcal{L}_{\mathrm{tw}}) v \in L^2_Q$.
In particular, the resolvent set of $\mathcal{L}_{\mathrm{tw}}$
is contained in the resolvent set of $\mathcal{L}_{\mathrm{tw}}^Q$.
The stated properties now follow in a standard fashion;
see, for example, \cite[Prop 3.1.5]{lorenzi2004analytic}.
\end{proof}

In order to define our final regularity concept,
we need to introduce the Hardy spaces
\begin{equation}
\begin{array}{lcl}
H^1(\Sigma_{\varphi} )
 & = & \{f: \Sigma_{\varphi} \to \mathbb{C} \hbox{ holomorphic for which }
\\[0.2cm]
& & \qquad \qquad
 \norm{f}_{H^1(\Sigma_{\varphi})}
  := \sup_{\abs{\nu} < \varphi }
      \int_0^\infty t^{-1} f( e^{i \nu} t ) \, dt < \infty \},
\\[0.2cm]
H^\infty(\Sigma_{\varphi} )
 & = & \{f: \Sigma_{\varphi} \to \mathbb{C} \hbox{ holomorphic for which }
\\[0.2cm]
& & \qquad \qquad
 \norm{f}_{H^\infty(\Sigma_{\varphi})} :=
   \sup_{z \in \Sigma_{\varphi} } \abs{f(z)} < \infty \} .
\end{array}
\end{equation}
If $\mathcal{L}$ is sectorial on a Banach space $X$,
then
for any $\omega(\mathcal{L}) < \varphi < \pi$
and any $h \in H^1(\Sigma_{\varphi})$ one can define
\begin{equation}
h(\mathcal{L}) = \frac{1}{2 \pi i}
  \int_{\partial \Sigma_{\nu} } R(z , \mathcal{L} ) h(z) \, dz
  \in \mathcal{L}(X, X)
\end{equation}
by picking an arbitrary $\nu \in (\omega(\mathcal{L}) , \varphi )$
and traversing the boundary in a downward fashion,
keeping the spectrum of $\mathcal{L}$ on the left.
It is however unclear if this integral converges
if we take $h \in  H^\infty(\Sigma_{\varphi})$.
The following result
states that this is indeed the case
for the sectorial operators discussed in
Lemma \ref{lem:nls:b:q:theta:sect}. Indeed, one can use a density argument to extend the conclusion 
to the whole space $H^\infty(\Sigma_{\varphi})$. 
Operators with this property
are said to admit a bounded $H^\infty$-calculus,
which is crucial for our
stochastic regularity estimates.

\begin{lem}
\label{lem:nls:h:inf:calc}
Assume that (HTw) and (HS) hold
and pick any $0 \le \theta \le \beta$.
There exists $\varphi \in (\omega(B^Q_{\theta}) , \frac{\pi}{2})$
together with a constant $K > 0$ so that for any
$h \in H^1(\Sigma_{\varphi}) \cap H^\infty(\Sigma_{\varphi})$
we have
\begin{equation}
\norm{ h(B^Q_{\theta}) } \le K \norm{h}_{H^\infty_{\varphi}} .
\end{equation}
\end{lem}
\begin{proof}
Since $\mathcal{L}_{\mathrm{tw}} - A_*$ is a first order
differential operator with continuous coefficients,
the perturbation theory described in \cite[{\S}8]{weis2006h}
can be applied to our setting.
In particular,
we can find constants $\Theta_0 \gg 1$
and $C_1 > 0$ together with an angle $\varphi_0 \in (\omega(B_{-\Theta_0}), \frac{\pi}{2} )$
for which
\begin{equation}
\norm{ h(B_{-\Theta_0}) }_{\mathcal{L}(L^2,L^2)} \le C_1 \norm{h}_{H^\infty_{\varphi_0}}
\end{equation}
holds for all $h \in H^1(\Sigma_{\varphi_0}) \cap H^\infty(\Sigma_{\varphi_0})$.
By restriction, we hence also have
\begin{equation}
\norm{ h(B^Q_{-\Theta_0}) }_{\mathcal{L}(L^2_Q,L^2_Q)} \le
\norm{ h(B_{-\Theta_0}) }_{\mathcal{L}(L^2,L^2)}
 C_1 \norm{h}_{H^\infty_{\varphi_0}}
\end{equation}
for all such $h$.
Fix two constants
\begin{equation}
\max\{ \omega(B^Q_{\theta} )  ,  \varphi_0 \} < \nu < \varphi < \frac{\pi}{2}
\end{equation}
and pick $h \in  H^1(\Sigma_{\varphi}) \cap H^\infty(\Sigma_{\varphi})$.
Using the resolvent identity,
we may compute
\begin{equation}
\begin{array}{lcl}
h(B^Q_{\theta} ) - h(B^Q_{-\Theta_0})
& = & \frac{1}{2\pi i}
 \int_{\partial\Sigma_\nu} h(z)
 \big[ R(z , B^Q_{\theta}) - R(z, B^Q_{-\Theta_0}) \big] \, dz
\\[0.2cm]
& = &
\frac{1}{2\pi i}
 \int_{\partial\Sigma_\nu} h(z)
 \big[ R(z , B^Q_{\theta}) - R(z - \theta - \Theta_0, B^Q_{\theta}) \big]
     \, dz
\\[0.2cm]
& = &
(\theta + \Theta_0)
\frac{1}{2\pi i} \int_{\partial\Sigma_\nu} h(z)
  R(z , B^Q_{\theta}) R(z - \theta - \Theta_0, B^Q_{\theta}) \ dz.
\end{array}
\end{equation}
Since zero is contained in the resolvent set of $B^Q_\theta$,
there exists $C_2 > 0$ for which
the estimate
\begin{equation}
\norm{R(z, B^Q_{\theta}) R(z - \theta - \Theta_0, B^Q_{\theta}) }_{\mathcal{L}(L^2_Q , L^2_Q)}
\le  \frac{C_2}{1 + \abs{z}^2}
\end{equation}
holds for all $z \in \partial \Sigma_{\nu}$.
This decays sufficently fast to ensure that
\begin{equation}
\norm{h(B^Q_{\theta} ) - h(B^Q_{-\Theta_0})}_{\mathcal{L}(L^2_Q , L^2_Q)} \le
C_3 \norm{h}_{H^\infty_{\varphi} }
\end{equation}
for some $C_3 > 0$ that does not depend on the choice of $h$.
The desired bound now follows from the inequality
\begin{equation}
\norm{h}_{H^\infty_{\varphi_0}} \le \norm{h}_{H^\infty_{\varphi}}.
\end{equation}
\end{proof}

Now that the formal framework has been set up, we are ready to
return to the quantity $N_{\e,\alpha}(t)$ defined in \sref{eq:nls:recdefN},
which is the main object of our interest.
For convenience, we use the shorthand notation
$\tau = \tau_{\e, \alpha}(T, \eta)$ ubiquitously
throughout the remainder of this section.
Writing $\nu = \alpha + \e$,
we introduce the splitting
\begin{equation}
\begin{array}{lcl}
N_{\e,\alpha;I}(t)
  & = & e^{\alpha t} \nrm{ V(t) }_{L^2}^2 ,
\\[0.2cm]
N_{\e,\alpha;II}(t)
 & = &
   \int_0^t e^{-\e(t-s)} e^{\alpha s}
       \nrm{V(s)}_{H^1}^2 \, ds
\\[0.2cm]
 & = &
  e^{-\e t}
    \int_0^t e^{\nu s}
       \nrm{V(s)}_{H^1}^2 \, ds.
\end{array}
\end{equation}

In order to understand $N_{\e,\alpha;I}$,
we introduce the expression
\begin{equation}
\begin{array}{lcl}
\mathcal{E}_{0}(t) &  = & S(t) Q V(0),
\\[0.2cm]
\end{array}
\end{equation}
together with the long-term integrals
\begin{equation}
\begin{array}{lclclcl}
\mathcal{E}^{\mathrm{lt}}_{F;\mathrm{lin}}(t) &  = &
  \int_0^{t-1} S(t  - s) Q F_{\mathrm{lin}}\big(V(s)\big)
  \mathbf{1}_{s < \tau} \, ds ,
\\[0.2cm]
\mathcal{E}^{\mathrm{lt}}_{F;\mathrm{nl}}(t) &  = &
  \int_0^{t-1} S(t  - s) Q F_{\mathrm{nl}}\big(V(s)\big)
  \mathbf{1}_{s < \tau} \, ds ,
\\[0.2cm]
\mathcal{E}^{\mathrm{lt}}_{B;\mathrm{lin}}(t)
  & = &
    \int_0^{t-1} S(t -s) Q B_{\mathrm{lin}}\big(V(s) \big)
     \mathbf{1}_{s < \tau} \,  d \beta_s,
\\[0.2cm]
\mathcal{E}^{\mathrm{lt}}_{B;\mathrm{cn}}(t)
  & = &
    \int_0^{t-1} S(t -s)Q B_{\mathrm{cn}}
    \mathbf{1}_{s < \tau} \, d \beta_s 
\\[0.2cm]
\end{array}
\end{equation}
and their short-term counterparts
\begin{equation}
\begin{array}{lclclcl}
\mathcal{E}^{\mathrm{sh}}_{F;\mathrm{lin}}(t) &  = &
  \int_{t-1}^{t} S(t - s) Q F_{\mathrm{lin}}\big(V(s)\big)
  \mathbf{1}_{s < \tau} \, ds ,
\\[0.2cm]
\mathcal{E}^{\mathrm{sh}}_{F;\mathrm{nl}}(t) &  = &
  \int_{t-1}^t S(t  - s) Q F_{\mathrm{nl}}\big(V(s)\big)
  \mathbf{1}_{s < \tau} \, ds ,
\\[0.2cm]
\mathcal{E}^{\mathrm{sh}}_{B;\mathrm{lin}}(t)
  & = &
    \int_{t-1}^t S(t -s) Q B_{\mathrm{lin}}\big(V(s) \big)
    \mathbf{1}_{s < \tau} \, d \beta_s,
\\[0.2cm]
\mathcal{E}^{\mathrm{sh}}_{B;\mathrm{cn}}(t)
  & = &
    \int_{t-1}^t S(t -s)Q B_{\mathrm{cn}}
    \mathbf{1}_{s < \tau} \, d \beta_s .
\\[0.2cm]
\end{array}
\end{equation}

Here we use the convention that integrands are set to zero for $s < 0$.
For convenience, we also write
\begin{equation}
\mathcal{E}_{F;\#}(t) = \mathcal{E}^{\mathrm{lt}}_{F;\#}(t) + \mathcal{E}^{\mathrm{sh}}_{F;\#}(t)
\end{equation}
for $\# \in \{\mathrm{lin},\mathrm{nl}\}$
and
\begin{equation}
\mathcal{E}_{B;\#}(t) = \mathcal{E}^{\mathrm{lt}}_{B;\#}(t) + \mathcal{E}^{\mathrm{sh}}_{B;\#}(t)
\end{equation}
for $\# \in \{\mathrm{lin},\mathrm{cn}\}$.

Turning to the terms in \sref{eq:nls:hsol:id:for:v}
that are relevant
for evaluating $N_{\e,\alpha;II}$,
we introduce the expression
\begin{equation}
\begin{array}{lcl}
\mathcal{I}_{\nu,\delta;0}(t) &  = &
   \int_0^t e^{\nu s}
    \norm{ S(\delta) \mathcal{E}_0(s) }_{H^1}^2 \, ds ,
\\[0.2cm]
\end{array}
\end{equation}
together with
\begin{equation}
\begin{array}{lcl}
\mathcal{I}^{\#}_{\nu,\delta;F;\mathrm{lin}}(t) &  = &
   \int_0^t e^{\nu s}
    \norm{ S(\delta) \mathcal{E}^{\#}_{F;\mathrm{lin}}(s) }_{H^1}^2 \, ds ,
\\[0.2cm]
\mathcal{I}^{\#}_{\nu,\delta;F;\mathrm{nl} }(t) &  = &
  \int_0^t e^{\nu s}
    \norm{ S(\delta) \mathcal{E}^{\#}_{F;\mathrm{nl}}(s) }_{H^1}^2 \, ds ,
\\[0.2cm]
\mathcal{I}^{\#}_{\nu,\delta;B;\mathrm{lin}}(t)
  & = &
    \int_0^t e^{\nu s}
    \norm{ S(\delta) \mathcal{E}^{\#}_{B;\mathrm{lin}}(s) }_{H^1}^2 \, ds  ,
\\[0.2cm]
 \mathcal{I}^{\#}_{\nu,\delta;B;\mathrm{cn}}(t)
  & = &
    \int_0^t e^{\nu s}
    \norm{ S(\delta) \mathcal{E}^{\#}_{B;\mathrm{cn}}(s) }_{H^1}^2 \, ds
\\[0.2cm]
\end{array}
\end{equation}
for $\# \in \{\mathrm{lt}, \mathrm{sh}\}$. The extra $S(\delta)$ factor will be used
to ensure that all the integrals we encounter are well-defined.
We emphasize that all our estimates are uniform in $0 < \delta < 1$,
allowing us to take $\delta \downarrow 0$.
The estimates concerning $\mathcal{I}^{\mathrm{sh}}_{\nu,\delta;F;\mathrm{nl}}$
and $\mathcal{I}^{\mathrm{sh}}_{\nu,\delta;B;\mathrm{lin}}$
in Lemmas \ref{lem:nls:f:nl:st} and \ref{lem:nls:b:lin:st:i}
are particularly delicate in this respect, as a direct
application of the bounds in Lemma \ref{lem:nls:sem:group:decay}
would result in expressions that diverge as $\delta \downarrow 0$.

\subsection{Deterministic regularity estimates}
\label{sec:nls:reg:ests:det}
In this part we set out to analyze the deterministic
integrals in \sref{eq:nls:hsol:id:for:v}.
The main complication is that we only have integrated
control over the squared $H^1$-norm of $V$. This is particularly
delicate for
$\mathcal{I}^{\mathrm{sh}}_{\nu, \delta;F;\mathrm{nl}}$,
where the nonlinearity itself is quadratic in $V$.

\begin{lem}
\label{lem:nls:e:zero}
Fix $T > 0$ and assume that (HA), (HTw), (HS), (H$\beta$),
(hSol) and (hFB) all hold.
Pick two constants $\e > 0$, $\alpha \ge 0$
for which $ \e + \frac{\alpha}{2} < \beta$
and write $\nu = \alpha + \e$.
Then for any
$0 \le \delta < 1$
and any $0 \le t \le T$,
we have the bound
\begin{equation}
\begin{array}{lcl}
e^{\alpha t} \norm{\mathcal{E}_0(t)}_{L^2}^2
 & \le &  M^2e^{-\e t} \norm{V(0)}^2_{L^2},
\end{array}
\end{equation}
together with
\begin{equation}
\begin{array}{lcl}
 e^{-\e t} \mathcal{I}_{\nu,\delta;0}(t)
 & \le &  \frac{M^2}{2 \beta - \nu}e^{-\e t} \norm{V(0)}^2_{H^1},
\end{array}
\end{equation}
\end{lem}
\begin{proof}
We compute
\begin{equation}
\begin{array}{lcl}
e^{\a t }\nrm{\mathcal{E}_{0}(t)}_{L^2}^2
  &\leq  & M^2 e^{ \alpha t} e^{-2 \beta t}
  \nrm{V(0)}_{L^2}^2  \\[0.2cm]
&\leq & M^2e^{-\e t}\nrm{V(0)}_{L^2}^2,
\end{array}
\end{equation}
together with
\begin{equation}
\begin{array}{lcl}
 e^{-\e t} \mathcal{I}_{\nu,\delta;0}(t)
& \le &
M^2 e^{-\e  t} \int_0^{t}e^{\nu s} e^{-2 \beta( s + \delta ) } \norm{V(0)}_{H^1}^2 \, ds
\\[0.2cm]
& \le &
   \frac{M^2}{2\b-\nu}e^{-\e t}\nrm{V(0)}^2_{H^1} .
\end{array}
\end{equation}
\end{proof}

\begin{lem}
Fix $T > 0$ and assume that (HA), (HTw), (HS), (H$\beta$),
(hSol) and (hFB) all hold.
Pick two constants $\e > 0$, $\alpha \ge 0$
for which $ \e + \frac{\alpha}{2} < \beta$
and write $\nu = \alpha + \e$.
Then for any
$0 \le \delta < 1$
and any $0 \le t \le \tau$,
we have the bound
\begin{equation}
\begin{array}{lcl}
e^{\alpha t} \norm{\mathcal{E}_{F;\mathrm{lin}}(t)}_{L^2}^2
 & \le &  K_{F;\mathrm{lin}}^2 \frac{M^2}{2 \beta - \nu} N_{\e,\alpha;II}(t),
\end{array}
\end{equation}
together with
\begin{equation}
\begin{array}{lcl}
 e^{-\e t} \mathcal{I}^{\mathrm{lt}}_{\nu,\delta;F;\mathrm{lin}}(t)
 & \le &
  K_{F;\mathrm{lin}}^2 \frac{M^2}{2(\beta  + \frac{\alpha}{2} - \nu)
    \e} N_{\e, \alpha;II}(t).
\end{array}
\end{equation}
\end{lem}
\begin{proof}
We first observe that
\begin{equation}
\begin{array}{lcl}
\nrm{\mathcal{E}_{F;\mathrm{lin}}(t)}_{L^2}^2
& \le &
   K_{F;\mathrm{lin}}^2 M^2  \left(
     \int_0^{t} e^{-\beta(t - s) } \norm{V(s)}_{H^1} \, ds \right)^2 ,
\\[0.2cm]
\norm{S(\delta)\mathcal{E}^{\mathrm{lt}}_{F;\mathrm{lin}}(t)}_{H^1}^2
& \le &
   K_{F;\mathrm{lin}}^2 M^2  \left(
     \int_0^{t} e^{-\beta(t - s) } \norm{V(s)}_{H^1} \, ds \right)^2 .
\end{array}
\end{equation}
This allows us to compute
\begin{equation}
\begin{array}{lcl}
e^{\a t }\nrm{\mathcal{E}_{F;\mathrm{lin}}(t)}_{L^2}^2
  & \le &
  K_{F;\mathrm{lin}}^2 M^2 e^{\alpha t} \left(
     \int_0^{t} e^{- (\beta - \frac{\nu}{2} )(t - s) } e^{-\frac{\nu}{2} ( t - s) } \norm{V(s)}_{H^1} \, ds \right)^2
\\[0.2cm]
& \le &
  K_{F;\mathrm{lin}}^2 \frac{M^2}{2 \beta - \nu} e^{\alpha t}
     \int_0^{t} e^{- \nu ( t - s) } \norm{V(s)}_{H^1}^2 \, ds
\\[0.2cm]
& = &
  K_{F;\mathrm{lin}}^2 \frac{M^2}{2 \beta - \nu} N_{\e,\alpha;II}(t) .
\end{array}
\end{equation}

Exploiting the inequality
$2\beta - \nu >  \e$,    
we write
\begin{equation}
\gamma_2 = \frac{\e + \nu}{2 \beta} < 1
\end{equation}
and observe that
\begin{equation}
2 \gamma_2 \beta - \nu = \e.
\end{equation}
Upon fixing
$\gamma_1 = 1 - \gamma_2$,
we readily see that
\begin{equation}
2 \gamma_1 \beta
 = 2 \beta - \e - \nu
 = 2\big(\beta + \frac{\alpha}{2} - \nu \big).
\end{equation}
This allows us to compute
\begin{equation}
\begin{array}{lcl}
 e^{-\e t} \mathcal{I}^{\mathrm{lt}}_{\nu,\delta ; F;\mathrm{lin}}(t)
  & \le &
  K_{F;\mathrm{lin}}^2 M^2 e^{-\e t}
  \int_0^t e^{\nu s}
  \left(
     \int_0^{s} e^{- \beta(s - s') }
         \norm{V(s')}_{H^1} \, ds' \right)^2  \, ds
\\[0.2cm]
& \le &
  K_{F;\mathrm{lin}}^2 M^2 e^{-\e t}
  \int_0^t e^{\nu s}
  \left(
     \int_0^{s} e^{- 2 \gamma_1 \beta(s - s') } \, ds' \right)
  \left(\int_0^s e^{-2 \gamma_2 \beta(s - s') }
         \norm{V(s')}^2_{H^1} \, ds' \right)  \, ds
\\[0.2cm]
&\le&
  K_{F;\mathrm{lin}}^2 \frac{M^2}{2 \gamma_1 \beta} e^{-\e t}
     \int_0^{t} e^{\nu s}
       \int_0^s e^{-2 \gamma_2 \beta(s - s')} \norm{V(s')}_{H^1}^2 \, ds' \, ds
\\[0.2cm]
& = &
    K_{F;\mathrm{lin}}^2 \frac{M^2}{2 \gamma_1 \beta} e^{-\e t}
     \int_0^{t} \int_{s'}^t e^{\nu s}
        e^{-2 \gamma_2 \beta(s - s')} \norm{V(s')}_{H^1}^2 \,  \, ds \, ds'
\\[0.2cm]
& = &
      K_{F;\mathrm{lin}}^2 \frac{M^2}{2 \gamma_1 \beta} e^{-\e t}
     \int_0^{t}  \Big[ \int_{s'}^t  e^{- (2 \gamma_2 \beta - \nu) s }  \, ds \Big]
        e^{ 2 \gamma_2 \beta s'} \norm{V(s')}_{H^1}^2 \,  \, ds'
\\[0.2cm]
& \le &
        K_{F;\mathrm{lin}}^2 \frac{M^2}{(2 \gamma_1 \beta)( 2 \gamma_2 \beta - \nu)}
         e^{-\e t}
     \int_0^{t}  e^{- (2 \gamma_2 \beta - \nu )s'}
        e^{ 2 \gamma_2 \beta s'} \norm{V(s')}_{H^1}^2 \,  \, ds'
\\[0.2cm]
& = &
  K_{F;\mathrm{lin}}^2 \frac{M^2}{(2 \gamma_1 \beta)( 2 \gamma_2 \beta - \nu)}
     e^{-\e t}
     \int_0^{t} e^{\nu s'}  \norm{V(s')}_{H^1}^2 \,  \, ds'
\\[0.2cm]
& = &
  K_{F;\mathrm{lin}}^2 \frac{M^2}{2(\beta  + \frac{\alpha}{2} - \nu)
    \e} N_{\e, \alpha;II}(t).
\\[0.2cm]
\end{array}
\end{equation}
\end{proof}

\begin{lem}
Fix $T > 0$ and assume that (HA), (HTw), (HS),
(H$\beta$),
(hSol) and (hFB) all hold.
Pick two constants $\e > 0$, $\alpha \ge 0$
and write $\nu = \alpha + \e$.
Then for
any $0 \le \delta < 1$
and any $0 \le t \le \tau$,
we have the bound
\begin{equation}
\begin{array}{lcl}
e^{-\e t} \mathcal{I}^{\mathrm{sh}}_{\nu,\delta;F;\mathrm{lin}}(t)
 & \le & 4 e^{\nu} M^2 K_{F;\mathrm{lin}}^2
 N_{\e,\alpha;II}(t) .
\end{array}
\end{equation}
\end{lem}
\begin{proof}
Using Cauchy-Schwarz, we compute
\begin{equation}
\begin{array}{lcl}
 e^{-\e t} \mathcal{I}^{\mathrm{sh}}_{\nu,\delta;F;\mathrm{lin}}(t)
& \le &
M^2 K_{F;\mathrm{lin}}^2
e^{-\e t}\int_0^t e^{\nu s}
  \left(\int_{s-1}^{s} \frac{1}{\sqrt{s+\d-s'}} \nrm{V(s')}_{H^1} ds' \right)^2 \, ds
  \\[0.2cm]
& \le &
M^2 K_{F;\mathrm{lin}}^2
e^{-\e t}\int_0^t e^{\nu s}
  \left(\int_{s-1}^{s} \frac{1}{\sqrt{s+\d-s'}} \, ds' \right)
  \left( \int_{s-1}^{s}\frac{1}{\sqrt{s+\d-s'}} \nrm{V(s')}_{H^1}^2 ds' \right) \, ds
  \\[0.2cm]
& \le &
2 M^2 K_{F;\mathrm{lin}}^2
e^{-\e t}\int_0^t e^{\nu s}
   \left( \int_{s-1}^{s}\frac{1}{\sqrt{s+\d-s'}} \nrm{V(s')}_{H^1}^2 ds' \right) \, ds
  \\[0.2cm]
& = &
2 M^2 K_{F;\mathrm{lin}}^2
e^{-\e t}  \int_0^t  \Big[\int_{s'}^{\min\{ t,s'+1 \} }
                      e^{\nu s} \frac{1}{\sqrt{s+\d-s'}}   \, ds
                   \Big]   \nrm{V(s')}_{H^1}^2       \, ds '
\\[0.2cm]
& \le &
4 e^{\nu} M^2 K_{F;\mathrm{lin}}^2
e^{-\e t}  \int_0^t  e^{\nu s'}   \nrm{V(s')}_{H^1}^2       \, ds '
\\[0.2cm]
& = &
4 e^{\nu} M^2 K_{F;\mathrm{lin}}^2
 N_{\e,\alpha;II}(t).
\\[0.2cm]
\end{array}
\end{equation}
\end{proof}

\begin{lem}
Fix $T > 0$ and assume that (HA), (HTw), (HS), (H$\beta$),
(hSol) and (hFB) all hold.
Pick two constants $\e > 0$, $\alpha \ge 0$
for which $ \e + \frac{\alpha}{2} < \beta$
and write $\nu = \alpha + \e$.
Then for any
$\eta > 0$,
any $0 \le \delta < 1$
and any $0 \le t \le \tau$,
we have the bound
\begin{equation}
\begin{array}{lcl}
e^{\alpha t} \norm{\mathcal{E}_{F;\mathrm{nl}}(t)}_{L^2}^2
 & \le &  \eta K_{F;\mathrm{nl}}^2 M^2( 1 + \eta^m)^2
  N_{\e,\alpha;II}(t),
\end{array}
\end{equation}
together with
\begin{equation}
\begin{array}{lcl}
 e^{-\e t} \mathcal{I}^{\mathrm{lt}}_{\nu,\delta;F;\mathrm{nl}}(t)
 & \le &    \eta K_{F;\mathrm{nl}}^2 ( 1 + \eta^m)^2
   \frac{M^2}{ \beta + \frac{\alpha}{2} - \nu }
      N_{\e,\alpha;II}(t) .
\end{array}
\end{equation}
\end{lem}
\begin{proof}
We first notice that
\begin{equation}
\begin{array}{lcl}
\nrm{\mathcal{E}_{F;\mathrm{nl}}(t)}_{L^2}^2
& \le &
   K^2_{F;\mathrm{nl}}(1 + \eta^m)^2 M^2  \left(
     \int_0^{t} e^{-\beta(t - s) } \norm{V(s)}^2_{H^1} \, ds
     \right)^2 ,
\\[0.2cm]
\norm{S(\delta)\mathcal{E}^{\mathrm{lt}}_{F;\mathrm{nl}}(t)}_{H^1}^2
& \le &
     K^2_{F;\mathrm{nl}}(1 + \eta^m)^2 M^2  \left(
     \int_0^{t} e^{-\beta(t - s) } \norm{V(s)}^2_{H^1} \, ds
     \right)^2 .
\end{array}
\end{equation}
Using $\beta > \nu - \frac{1}{2} \alpha
= \frac{1}{2} \alpha + \e$,
we compute
\begin{equation}
\begin{array}{lcl}
  \int_0^{t} e^{-\beta(t - s) } \norm{V(s)}_{H^1}^2 \, ds
& = &
  e^{\frac{\alpha}{2} t}
  \int_0^{t} e^{-\beta(t - s) } e^{-\frac{\alpha}{2} t} \norm{V(s)}_{H^1}^2 \, ds
\\[0.2cm]
& \le &
  e^{\frac{\alpha}{2} t}
  \int_0^{t} e^{-\beta(t - s) } e^{-\frac{\alpha}{2} (t-s)} \norm{V(s)}_{H^1}^2 \, ds
\\[0.2cm]
& \le &
  e^{\frac{\alpha}{2} t}
  \int_0^{t} e^{-\nu(t - s) } \norm{V(s)}_{H^1}^2 \, ds.
\end{array}
\end{equation}
This yields the desired bound
\begin{equation}
\begin{array}{lcl}
e^{\a t }\nrm{\mathcal{E}_{F;\mathrm{nl}}(t)}_{L^2}^2
  &\leq  &
K_{F;\mathrm{nl}}^2(1 + \eta^m )^2 M^2 e^{\alpha t} \left(
     \int_0^{t} e^{-\beta(t - s) } \norm{V(s)}_{H^1}^2 \, ds \right)^2
 \\[0.2cm]
& \le &
  K_{F;\mathrm{nl}}^2(1 + \eta^m )^2 M^2 e^{2 \alpha t} \left(
     \int_0^{t} e^{-\nu(t - s) }  \norm{V(s)}_{H^1}^2 \, ds \right)^2
\\[0.2cm]
& \le &
  K_{F;\mathrm{nl}}^2(1 + \eta^m )^2 M^2 \eta N_{\e,\alpha;II}(t) .
\\[0.2cm]
\end{array}
\end{equation}
In a similar spirit, we compute
\begin{equation}
\begin{array}{lcl}
 e^{-\e t} \mathcal{I}^{\mathrm{lt}}_{\nu,\delta ; F;\mathrm{nl}}(t)
  & \le &
  K_{F;\mathrm{nl}}^2 (1 + \eta^m )^2 M^2 e^{-\e t}
  \int_0^t e^{\nu s}
  \left(
     \int_0^{s} e^{- \beta(s - s') }
         \norm{V(s')}^2_{H^1} \, ds' \right)^2  \, ds
\\[0.2cm]
& \le &
  K_{F;\mathrm{nl}}^2 (1 + \eta^m )^2 M^2 e^{-\e t}
  \int_0^t e^{\nu s}
  e^{\frac{\alpha }{2}s }
  \left(
     \int_0^{s} e^{- \nu(s - s') }
         \norm{V(s')}^2_{H^1} \, ds' \right)
\\[0.2cm]
& & \qquad \qquad \times
\left(
     \int_0^{s} e^{- \beta(s - s') }
         \norm{V(s')}^2_{H^1} \, ds' \right)  \, ds
\\[0.2cm]
& \le &
  \eta K_{F;\mathrm{nl}}^2 (1 + \eta^m)^2 M^2 e^{-\e t}
  \int_0^t e^{\nu s} e^{-\frac{\alpha}{2}s }
     \int_0^{s} e^{- \beta(s - s') }
         \norm{V(s')}^2_{H^1} \, ds' \,  \, ds
\\[0.2cm]
& = &
   \eta K_{F;\mathrm{nl}}^2 (1 + \eta^m )^2 M^2 e^{-\e t}
     \int_0^{t}   \Big[ \int_{s'}^t  e^{- ( \frac{\alpha}{2} - \nu + \beta) s }  \, ds \Big]
        e^{ \beta s'} \norm{V(s')}_{H^1}^2 \,  \, ds'
\\[0.2cm]
& \le &
   \eta K_{F;\mathrm{nl}}^2 (1 + \eta^m )^2 \frac{M^2}{ \beta + \frac{\alpha}{2} - \nu }
      e^{-\e t}
     \int_0^{t}   e^{- ( \frac{\alpha}{2} - \nu + \beta) s' }  e^{ \beta s'} \norm{V(s')}_{H^1}^2 \,  \, ds'
\\[0.2cm]
& = &
   \eta K_{F;\mathrm{nl}}^2 (1 + \eta^m )^2 \frac{M^2}{ \beta + \frac{\alpha}{2} - \nu }
      e^{-\e t}
     \int_0^{t}  e^{\nu s'} e^{-  \frac{\alpha}{2}  s' }  \norm{V(s')}_{H^1}^2 \,  \, ds'
\\[0.2cm]
& \le &
   \eta K_{F;\mathrm{nl}}^2 (1 + \eta^m )^2 \frac{M^2}{ \beta + \frac{\alpha}{2} - \nu }
      N_{\e,\alpha;II}(t).
\\[0.2cm]
\end{array}
\end{equation}
\end{proof}

\begin{lem}
\label{lem:nls:f:nl:st}
Fix $T > 0$ and assume that (HA), (HTw), (HS), (H$\beta$),
(hSol) and (hFB) all hold.
Pick two constants $\e > 0$, $\alpha \ge 0$
and write $\nu = \alpha + \e$.
Then for any
$\eta > 0$,
any $0 \le \delta < 1$
and any $0 \le t \le \tau$,
we have the bound
\begin{equation}
\begin{array}{lcl}
e^{-\e t} \mathcal{I}^{\mathrm{sh}}_{\nu,\delta;F,\mathrm{nl}}(t) & \le &
 \eta M^2 K_{F;\mathrm{nl}}^2 (1 + \eta^m)^2 (1 + \rho^{-1})
   e^{3 \nu} (4  + \nu)
   N_{\e,\alpha;II}(t) .
\end{array}
\end{equation}
\end{lem}
\begin{proof}
We start by observing that
\begin{equation}
\norm{v}_{H^1}^2 = \norm{v}_{L^2}^2 + \rho^{-1} \norm{A_*^{1/2} v}_{L^2}^2.
\end{equation}
In addition, for any $w \in L^2$, $\vartheta > 0$, $\vartheta_A \ge 0$
and $\vartheta_B \ge 0$
we have
\begin{equation}
\label{eq:nls:deriv:semigroup}
\begin{array}{lcl}
\frac{d}{d \vartheta} \langle S(\vartheta + \vartheta_A) w,
   S(\vartheta + \vartheta_B) w \rangle_{L^2}
& = &
  \langle \mathcal{L}_\mathrm{tw} S(\vartheta + \vartheta_A) w ,
     S(\vartheta +\vartheta_B) w \rangle_{L^2}
\\[0.2cm]
& & \qquad
+ \langle  S(\vartheta + \vartheta_A) w ,
    \mathcal{L}_\mathrm{tw}S(\vartheta + \vartheta_B) w \rangle_{L^2}
\\[0.2cm]
& = &
  \langle  S(\vartheta + \vartheta_A) w ,
      \mathcal{L}_\mathrm{tw}^\mathrm{adj} S(\vartheta + \vartheta_B) w \rangle_{L^2}
\\[0.2cm]
& & \qquad
+ \langle  S(\vartheta + \vartheta_A) w ,
   \mathcal{L}_\mathrm{tw}  S(\vartheta + \vartheta_B) w \rangle_{L^2}
\\[0.2cm]
& = &
  \langle  S(\vartheta + \vartheta_A) w ,
     \big[ \mathcal{L}_\mathrm{tw}^\mathrm{adj} - A_*\big]
       S(\vartheta + \vartheta_B) w \rangle_{L^2}
\\[0.2cm]
& & \qquad
+ \langle  S(\vartheta + \vartheta_A) w ,  \
  \big[ \mathcal{L}_\mathrm{tw} - A_*\big]  S(\vartheta + \vartheta_B) w
     \rangle_{L^2}
\\[0.2cm]
& & \qquad
  - 2 \langle A_*^{1/2} S(\vartheta + \vartheta_A) w,
    A_*^{1/2} S(\vartheta + \vartheta_B) w \rangle_{L^2} .
\end{array}
\end{equation}

Assume for the moment that $\delta > 0$.
For convenience,
we introduce the expression
\begin{equation}
\mathcal{E}_{s,s',s''; \mathcal{H}}
= \langle S(s+\delta -s')Q F_{\mathrm{nl}}\big(V(s')\big) ,
  S(s + \delta - s'') Q F_{\mathrm{nl}}\big(V(s'')\big)
\rangle_{\mathcal{H}},
\end{equation}
where we allow $\mathcal{H} \in \{L^2, H^1 \}$.
Exploiting \sref{eq:nls:deriv:semigroup}
and the fact that $\delta > 0$,
we obtain the bound
\begin{equation}
\begin{array}{lcl}
\mathcal{E}_{s,s',s'';H^1}
& \le & M^2 K_{F;\mathrm{nl}}^2 (1 + \eta^m)^2 \norm{V(s')}_{H^1}^2
  \norm{V(s'')}_{H^1}^2
\\[0.2cm]
& & \qquad
 +  M^2 K_{F;\mathrm{nl}}^2 (1 + \eta^m)^2
 \rho^{-1} \frac{1}{\sqrt{s + \delta - s''} }
    \norm{V(s')}_{H^1}^2 \norm{V(s'')}_{H^1}^2
\\[0.2cm]
& & \qquad
  - \rho^{-1} \frac{1}{2}
  \frac{d}{ds}
    \mathcal{E}_{s,s',s'';L^2}
\end{array}
\end{equation}
for the values of $(s,s',s'')$ that are relevant below.
Upon introducing the integrals
\begin{equation}
\begin{array}{lcl}
\mathcal{I}_{I}
 & = &
     e^{-\e t } \int_0^{t}
   e^{\nu s} \int_{s-1}^s \int_{s-1}^s
       \big[ 1 + \frac{1}{\sqrt{s + \delta - s''} } ]
          \norm{V(s')}_{H^1}^2 \norm{V(s'')}_{H^1}^2
             \, ds'' \, ds' \, ds ,
\\[0.2cm]
\mathcal{I}_{II}
& = &
 e^{-\e t } \int_0^{t}
   e^{\nu s} \int_{s-1}^s \int_{s-1}^s
    \frac{d}{ds}   \mathcal{E}_{s,s',s'';L^2} ds'' \, ds' \, ds ,
\\[0.2cm]
\end{array}
\end{equation}
we hence readily obtain the estimate
\begin{equation}
\begin{array}{lcl}
e^{-\e t} \mathcal{I}^{\mathrm{sh}}_{\nu,\delta;F;\mathrm{nl}}(t)
& \le &
  (1 + \rho^{-1}) M^2 K^2_{F;\mathrm{nl}}(1 + \eta^m)^2 \mathcal{I}_{I}
  - \frac{1}{2} \rho^{-1} \mathcal{I}_{II}.
\end{array}
\end{equation}
Changing the order of the integrals, we find
\begin{equation}
\begin{array}{lcl}
 \mathcal{I}_{I} & =&
e^{- \e t}  \int_0^{t}
  \int_{s' - 1}^{\min\{t, s' + 1 \} }
  \Big[
 \int_{\max\{s', s''\} }^{\min\{ t, s'+1, s''+1 \} }
                      e^{\nu s}
                      \big[
                         1 + \frac{1}{\sqrt{s+\d-s''}} \big]  \, ds
                   \Big]   \nrm{V(s')}_{H^1}^2  \nrm{V(s'')}_{H^1}^2       \, ds'' \, ds'
\\[0.2cm]
& \le &
3
 e^{- \e t}  \int_0^{t} e^{ \nu s' }  e^{ \nu }
   \nrm{V(s')}_{H^1}^2
   \int_{s' - 1}^{\min\{t, s' + 1 \} }
      \nrm{V(s'')}_{H^1}^2
      \, ds'' \, ds'
\\[0.2cm]
& \le &
3
 e^{- \e t}  \int_0^{t}  e^{\nu s'}  e^{ \nu }
   \nrm{V(s')}_{H^1}^2
   e^{2 \nu}\int_{s' - 1}^{\min\{t, s' + 1 \} }
      e^{-\nu(\min\{t, s'+1\} -s'')}\nrm{V(s'')}_{H^1}^2
      \, ds'' \, ds'
\\[0.2cm]
& \le &
3 \eta  e^{3 \nu}
 e^{- \e t}  \int_0^{t}  e^{\nu s'}
   \nrm{V(s')}_{H^1}^2
  e^{-\alpha \min\{t, s'+1 \} }
      \, ds'
\\[0.2cm]
& \le &
  3 \eta  e^{3 \nu}
  N_{\e,\alpha;II}(t) .
\end{array}
\end{equation}
In a similar fashion,
we may use an integration by parts to write
\begin{equation}
\begin{array}{lcl}
 \mathcal{I}_{II} & =&
  e^{-\e t}  \int_0^{t}
  \int_{s'-1}^{\min\{t , s'+1 \} }
 \Big[\int_{\max\{s',s''\}}^{\min\{ t,s'+1, s''+1 \} }
                      e^{\nu s}
                       \frac{d}{ds} \mathcal{E}_{s,s',s'';L^2}   \, ds
                   \Big]        \, ds '' \, ds'
 \\[0.2cm]
 & = &
   \mathcal{I}_{II;A}
   + \mathcal{I}_{II;B}
   + \mathcal{I}_{II;C},
\end{array}
\end{equation}
in which
 we have introduced
\begin{equation}
\begin{array}{lcl}
\mathcal{I}_{II;A}
 & = &
    e^{-\e t}  \int_0^{t}  \int_{s'-1}^{\min\{t, s'+1\} }
     e^{\nu \min\{ t , s' + 1, s''+1 \} }
       \mathcal{E}_{\min\{t, s'+1,s''+1\}, s',s'';L^2} \, ds'' \, ds' ,
 \\[0.2cm]
\mathcal{I}_{II;B}
 & = &
   - e^{-\e t}  \int_0^{t}  \int_{s'-1}^{\min\{t, s'+1\} }
     e^{\nu \max\{ s', s''\} }
       \mathcal{E}_{\max\{s',s''\},s',s'';L^2} \, ds'' \, ds' ,
\\[0.2cm]
\mathcal{I}_{II;C}
 & = &
   - e^{-\e t}  \int_0^{t}
  \int_{s'-1}^{\min\{t , s'+1 \} }
 \Big[\int_{\max\{s',s''\}}^{\min\{ t,s'+1, s''+1 \} }
                      \nu e^{\nu s}
    \mathcal{E}_{s,s',s'';L^2}    \, ds
                   \Big]        \, ds '' \, ds'  .
 \\[0.2cm]
\end{array}
\end{equation} Note here that $\mathcal{I}_{II;B}$ is well defined because $\delta>0$. 
A direct inspection of these terms yields the bound
\begin{equation}
\begin{array}{lcl}
\abs{\mathcal{I}_{II} }
& \le &
  e^{ \nu} (2 + \nu)
  M^2 K_{F;\mathrm{nl}}^2 (1 + \eta^m)^2
     e^{-\e t}
    \int_0^{t }
       e^{\nu s'}
        \norm{V(s')}_{H^1}^2
       \int_{s'-1}^{\min\{t, s'+1 \} }
         \norm{V(s'')}_{H^1}^2 \, ds'' \, d s'
\\[0.2cm]
& \le &
  e^{ \nu} (2 + \nu)
  M^2 K_{F;\mathrm{nl}}^2 (1 + \eta^m)^2
     e^{-\e t}
    \int_0^{t }
       e^{\nu s'} \norm{V(s')}_{H^1}^2
       e^{2 \nu}
    \\[0.2cm]
    & & \qquad \qquad \times
       \int_{s'-1}^{\min\{t, s'+1 \} }
          e^{-\nu(\min\{t,s'+1\} - s'')} \norm{V(s'')}_{H^1}^2 \, ds'' \, d s'
\\[0.2cm]
& \le &
  \eta e^{3 \nu} (2 + \nu)
  M^2 K_{F;\mathrm{nl}}^2 (1 + \eta^m)^2
     e^{-\e t}
     \int_0^{t }
       e^{\nu s'} \norm{V(s')}_{H^1}^2
       e^{-\alpha \min\{t, s'+1\} }
        \, d s'
\\[0.2cm]
& \le &
\eta e^{3 \nu} (2 + \nu)
  M^2 K_{F;\mathrm{nl}}^2 (1 + \eta^m)^2
  N_{\e,\alpha;II}(t).
\end{array}
\end{equation}
It hence remains to consider the case $\delta = 0$.
We may apply Fatou's lemma
to conclude
\begin{equation}
\begin{array}{lcl}
\mathcal{I}^{\mathrm{sh}}_{\nu,0;F;\mathrm{nl}} (t)
& = &
  \int_0^{t}
e^{\nu s}
(\lim_{\d\to0}\nrm{S(\d)
  \mathcal{E}^{\mathrm{sh}}_{B;\mathrm{lin}}(s) }_{H^1})^2
    \mathbf{1}_{s < \tau} \, ds
\\[0.2cm]
&\leq & \liminf_{\d\to0}
   \mathcal{I}^{\mathrm{sh}}_{\nu,\delta;F;\mathrm{nl}} (t).
\\[0.2cm]
\end{array}
\end{equation}
The result now follows from the fact that
the bounds obtained above do not depend on  $\delta$.
\end{proof}

\subsection{Stochastic regularity estimates}
\label{sec:nls:reg:ests:st}

We are now ready to discuss the stochastic integrals
in \sref{eq:nls:hsol:id:for:v}.
These require special care because they cannot be bounded
in a pathwise fashion, unlike the deterministic integrals above.
Expectations of suprema are particularly delicate in this respect.
Indeed, the powerful Burkholder-Davis-Gundy inequalities cannot
be directly applied to the stochastic convolutions that arise in our
mild formulation. However, the $H^\infty$-calculus obtained
in Lemma \ref{lem:nls:h:inf:calc}
allows us to use the following mild version,
which is the source of the extra $T$ factors
that appear in our estimates.

\begin{lem}
\label{lem:nls:mild:regularity}
Fix $T > 0$ and assume that (HA), (HTw), (HS) and (H$\beta$)
all hold.
There exists a constant $K_{\mathrm{cnv}} > 0$ so that for any
$W \in \mathcal{N}^2([0,T];(\mathcal{F})_{t} ; L^2)$
and any $0 \le \alpha \le 2 \beta$ we have
\begin{equation}
E \sup_{0 \le t \le T} e^{\alpha t}
  \norm{ \int_0^t S( t - s) Q W(s) \, d \beta_s }_{L^2}^2
  \le K_{\mathrm{cnv}} E \int_0^T e^{\alpha s} \norm{W(s)}_{L^2}^2 \, d s.
\end{equation}
\end{lem}
\begin{proof}
Lemma \ref{lem:nls:h:inf:calc} implies that the generator
$B^Q_{\alpha} = \mathcal{L}_{\mathrm{tw}} + \frac{1}{2}\alpha$ of the
semigroup $e^{\frac{1}{2} \alpha t} S(t)$
on $L^2_Q$ satisfies assumption (H) in \cite{veraar2011note}.
On account of the identity
\begin{equation}
e^{\alpha t}
  \norm{ \int_0^t S( t - s) Q W(s) \, d \beta_s }_{L^2}^2
= \norm{ \int_0^t e^{\frac{1}{2} \alpha(t-s)}S( t - s) Q e^{\frac{1}{2} \alpha s} W(s) \, d \beta_s }_{L^2}^2,
\end{equation}
the desired estimate is now an immediate consequence of
\cite[Thm. 1.1]{veraar2011note}.
\end{proof}

\begin{lem}
\label{lem:nls:b:lin:sup:e}
Fix $T > 0$ and assume that (HA), (HTw), (HS),
(H$\beta$),
(hSol) and (hFB) all hold.
Then for any pair of constants
$\e > 0$
and $0 \le \alpha \le 2 \beta$ we have the bound
\begin{equation}
\begin{array}{lcl}
E \sup_{0 \le t \le \tau} e^{\alpha t } \norm{\mathcal{E}_{B;\mathrm{lin}}(t )}_{L^2}^2
 & \le &  (T + 1)  K_{\mathrm{cnv}} K_{B;\mathrm{lin}}^2 e^{\e}
   E \sup_{0 \le t \le \tau}   N_{\e,\a;II}(t) .
\end{array}
\end{equation}
\end{lem}
\begin{proof}
Using Lemma \ref{lem:nls:mild:regularity}
we compute
\begin{equation}
\label{eq:nls:bnd:on:e:b:lin:dir:from:mild:fmr}
\begin{array}{lcl}
E \sup_{0 \le t \le \tau} e^{\alpha t } \norm{\mathcal{E}_{B;\mathrm{lin}}(t )}_{L^2}^2
  & \le   &
  E \sup_{0 \le t \le T} e^{\alpha t } \norm{\mathcal{E}_{B;\mathrm{lin}}(t )}_{L^2}^2
\\[0.2cm]
& = &
    E \sup_{0 \le t \le T} e^{\alpha t }
      \norm{ \int_0^{t} S(t-s) Q B_{\mathrm{lin}}\big(V(s)\big)
        \mathbf{1}_{s < \tau} \, d \beta_s
       }_{L^2}^2
\\[0.2cm]
& \le &
  K_{\mathrm{cnv}} E
      \int_{0}^T
        e^{ \alpha s}
        \norm{
           B_{\mathrm{lin}}\big(V(s)\big)
            \mathbf{1}_{s < \tau}
         }_{L^2}^2 \, ds
\\[0.2cm]
& \le &
  K_{\mathrm{cnv}} K_{B;\mathrm{lin}}^2
  E \int_0^\tau
      e^{ \alpha s} \norm{V(s)}_{H^1}^2 \, ds .
\end{array}
\end{equation}
By dividing up the integral, we obtain
\begin{equation}
\begin{array}{lcl}
\int_0^\tau
      e^{ \alpha s} \norm{V(s)}_{H^1}^2 \, ds
& \le &
    e^{\e} \int_0^{1 }
     e^{ -\e(1 - s) } e^{ \alpha s}
       \norm{V(s)}_{H^1}^2 \mathbf{1}_{s < \tau}  \, ds
\\[0.2cm]
& & \qquad
  +
     e^{\e} \int_1^{2 }
        e^{ -\e(2 - s) } e^{ \alpha s}
        \norm{V(s)}_{H^1}^2 \mathbf{1}_{s < \tau}  \, ds
\\[0.2cm]
& & \qquad
  + \cdots
+
   e^{\e} \int_{\lfloor T \rfloor}^{\lfloor T \rfloor + 1}
       e^{ -\e(\lfloor T \rfloor +1 - s) } e^{ \alpha s}
       \norm{V(s)}_{H^1}^2 \mathbf{1}_{s < \tau}  \, ds
\\[0.2cm]
& \le &
  (T + 1)   e^{\e}
     \sup_{0 \le t \le T + 1}
      \int_0^t e^{-\e(t - s)} e^{\alpha s}
        \norm{V(s)}_{H^1}^2 \mathbf{1}_{s < \tau} \, ds
\\[0.2cm]
& \le &
  (T + 1)   e^{\e}
     \sup_{0 \le t \le \tau}
      \int_0^t e^{-\e(t - s)} e^{\alpha s}
        \norm{V(s)}_{H^1}^2  \, ds
\\[0.2cm]
& = &
  (T + 1)   e^{\e}
     \sup_{0 \le t \le \tau}
      N_{\e,\a;II}(t) ,
\end{array}
\end{equation}
which yields the desired bound upon taking expectations.
\end{proof}

\begin{lem}
\label{lem:nls:b:cn:sup:e}
Fix $T > 0$ and assume that (HA), (HTw), (HS),
(H$\beta$),
(hSol) and (hFB) all hold.
Then 
we have the bound
\begin{equation}
\begin{array}{lcl}
E \sup_{0 \le t \le \tau} \norm{\mathcal{E}_{B;\mathrm{cn}}(t )}_{L^2}^2
 & \le &  T K_{\mathrm{cnv}} K_{B;\mathrm{cn}}^2 .
\end{array}
\end{equation}
\end{lem}
\begin{proof}
This bound follows directly from \sref{eq:nls:bnd:on:e:b:lin:dir:from:mild:fmr}
by picking
$\alpha = 0$
and making the substitutions
\begin{equation}
  K_{B;\mathrm{lin}} \mapsto K_{B;\mathrm{cn}},
  \qquad
  \qquad
    \nrm{V(s)}_{H^1}^2 \mapsto 1 .
\end{equation}
\end{proof}

We now set out to bound the expectation of the supremum of the remaining double integrals
$\mathcal{I}^{\#}_{\nu,\delta;B;\mathrm{lin}}(t)$
and $\mathcal{I}^{\#}_{\nu,\delta;B;\mathrm{cn}}(t)$ with $\# \in \{ \mathrm{lt} , \mathrm{sh} \}$.
This is performed in Lemmas \ref{lem:nls:st:sup:on:b:lin:i} and \ref{lem:nls:st:sup:on:b:cn:i}, 
but we first compute several time independent bounds for the expectation of the integrals themselves.

\begin{lem}
Fix $T > 0$ and assume that (HA), (HTw), (HS), (H$\beta$),
(hSol) and (hFB) all hold.
Pick constants $\e > 0$,
$\alpha \ge 0$
and write $\nu = \alpha + \e$.
Then for any $0 \le \delta < 1$
and $0 \le t \le T$,
we have the identities
\begin{equation}
\label{eq:nls:ito:isometry:i:lt}
\begin{array}{lcl}
E  \mathcal{I}^{\mathrm{lt}}_{\nu,\delta;B;\mathrm{lin}}(t )
& = &
 E  \int_0^{t}
   e^{ \nu s } \int_{0}^{s-1 }
          \norm{  S(s+\delta-s') Q B_{\mathrm{lin}}\big(V(s')\big) }_{L^2}^2
            \mathbf{1}_{s' < \tau}\, ds' \, ds ,
\\[0.2cm]
E \mathcal{I}^{\mathrm{lt}}_{\nu,\delta;B;\mathrm{cn}}(t)
& = &
 E  \int_0^{t}
   e^{ \nu s } \int_{0}^{s-1}
          \norm{  S(s+\delta-s') Q B_{\mathrm{cn}} }_{L^2}^2
          \mathbf{1}_{s' < \tau} \, ds' \, ds ,
\end{array}
\end{equation}
together with their short-time counterparts
\begin{equation}
\label{eq:nls:ito:isometry:i:sh}
\begin{array}{lcl}
E  \mathcal{I}^{\mathrm{sh}}_{\nu,\delta;B;\mathrm{lin}}(t)
& = &
 E  \int_0^{t}
   e^{ \nu s } \int_{s-1  }^{s }
          \norm{  S(s+\delta-s') Q B_{\mathrm{lin}}\big(V(s')\big) }_{L^2}^2
          \mathbf{1}_{s' < \tau} \, ds' \, ds ,
\\[0.2cm]
E   \mathcal{I}^{\mathrm{sh}}_{\nu,\delta;B;\mathrm{cn}}(t)
& = &
 E  \int_0^{t}
   e^{ \nu s } \int_{s-1}^s
          \norm{  S(s+\delta-s') Q B_{\mathrm{cn}} }_{L^2}^2
          \mathbf{1}_{s' < \tau} \, ds' \, ds .
\end{array}
\end{equation}
\end{lem}
\begin{proof}
The linearity of the expectation operator,
the It\^o Isometry (see e.g. \cite[\S2.3]{Concise}) and the integrability of the integrands
imply that
\begin{equation}
\begin{array}{lcl}
E \mathcal{I}^{\mathrm{lt}}_{B;\mathrm{lin}}(t )
& = &
  E
   \int_0^{t  }
        e^{\nu s}
        \norm{ \int_0^{s-1}
                   S(s+\delta-s') Q B_{\mathrm{lin}}\big(V(s')\big)
                     \mathbf{1}_{s' < \tau} d \beta_{s'} }_{H^1}^2
                 \, ds
\\[0.2cm]
& = & E   \int_0^{t}  e^{\nu s}
   \int_0^{s-1}
          \norm{
                   S(s+\delta-s') Q B_{\mathrm{lin}}\big(V(s')\big) }_{L^2}^2
                   \mathbf{1}_{s' < \tau} \, ds' \, ds .
\\[0.2cm]
\end{array}
\end{equation}
The remaining expressions follow in a similar fashion.
\end{proof}


\begin{lem}
\label{lem:nls:b:lin:lt:i}
Fix $T > 0$ and assume that (HA), (HTw), (HS),
(H$\beta$),
(hSol) and (hFB) all hold.
Pick constants $\e > 0$,
$\alpha \ge 0$
for which $\e + \a  < 2\beta$
and write $\nu = \alpha + \e$.
Then for
any $0 \le  \delta < 1$
and any $0 \le t \le T$,
we have the bound
\begin{equation}
\begin{array}{lcl}
E e^{-\e t} \mathcal{I}^{\mathrm{lt}}_{\nu,\delta;B;\mathrm{lin}}(t)
 & \le &    \frac{M^2}{2 \beta - \nu} K_{B;\mathrm{lin}}^2
  E  N_{\e;\alpha;II}(t \wedge \tau).
\end{array}
\end{equation}
\end{lem}
\begin{proof}
Using \sref{eq:nls:ito:isometry:i:lt}
and switching the integration order,
we obtain
\begin{equation}
\begin{array}{lcl}
E  e^{-\e t}\mathcal{I}^{\mathrm{lt}}_{\nu,\delta;B;\mathrm{lin}}(t)
  &\leq  &
  M^2 K_{B;\mathrm{lin}}^2 E  e^{-\e t}  \int_0^{t}
  e^{ \nu s }
    \int_0^{s \wedge \tau}  e^{-2 \beta(s-s')}\norm{V(s')}_{H^1}^2 \, ds' \, ds
\\[0.2cm]
& = &
  M^2 K_{B;\mathrm{lin}}^2 E  e^{-\e t }  \int_0^{t \wedge \tau}
    \Big[\int_{s'}^{t }
      e^{  -(2 \beta -\nu )s } \, d s
   \Big]
      e^{2 \beta s'}\norm{V(s')}_{H^1}^2  \, ds'
\\[0.2cm]
& \le &
  \frac{M^2}{2 \beta - \nu} K_{B;\mathrm{lin}}^2
  E  e^{ - \e t }  \int_0^{t\wedge \tau}
      e^{-(2 \beta - \nu)s'}
       e^{2 \beta s'}\norm{V(s')}_{H^1}^2  \, ds'
\\[0.2cm]
& = &
  \frac{M^2}{2 \beta - \nu} K_{B;\mathrm{lin}}^2
  E  e^{-\e t }  \int_0^{t\wedge \tau}
   e^{ \nu s'}  \norm{V(s')}_{H^1}^2  \, ds'
\\[0.2cm]
& \le &
  \frac{M^2}{2 \beta - \nu} K_{B;\mathrm{lin}}^2
    E  e^{-\e t\wedge \tau }  \int_0^{t\wedge \tau}
   e^{ \nu s'}  \norm{V(s')}_{H^1}^2  \, ds'
\\[0.2cm]
& = &
   \frac{M^2}{2 \beta - \nu} K_{B;\mathrm{lin}}^2
  E  N_{\e,\alpha;II}(t \wedge \tau).
\end{array}
\end{equation}
\end{proof}

\begin{lem}
\label{lem:nls:b:lin:st:i}
Fix $T > 0$ and assume that (HA), (HTw), (HS),
(H$\beta$),
(hSol) and (hFB) all hold.
Pick two constants $\e > 0$, $\alpha \ge 0$
and write $\nu = \alpha + \e$.
Then for
any $0 \le \delta < 1$
and any $0 \le t \le T$,
we have the bound
\begin{equation}
\begin{array}{lcl}
E
e^{-\e t} \mathcal{I}^{\mathrm{sh}}_{\nu,\delta;B;\mathrm{lin}}(t)
 & \le &     K_{B;\mathrm{lin}}^2 M^2 (1 + \rho^{-1}) e^{\nu}( 4  + \nu )
  E  N_{\e,\alpha;II}(t \wedge \tau).
\end{array}
\end{equation}
\end{lem}
\begin{proof}
We only consider the case $\delta > 0$ here,
noting that the limit $\delta \downarrow 0$
can be handled as in the proof of
Lemma \ref{lem:nls:f:nl:st}.
Applying the identity \sref{eq:nls:deriv:semigroup}
with $\vartheta_A = \vartheta_B$,
we obtain the bound
\begin{equation}
\begin{array}{lcl}
\norm{S(s+\delta -s')Q B_{\mathrm{lin}}\big(V(s')\big) }_{H^1}^2
& \le & M^2 K_{B;\mathrm{lin}}^2 \norm{V(s')}_{H^1}^2
\\[0.2cm]
& & \qquad
 +  M^2 K_{B;\mathrm{lin}}^2 \rho^{-1} \frac{1}{\sqrt{s + \delta - s'} }
    \norm{V(s')}_{H^1}^2
\\[0.2cm]
& & \qquad
  - \rho^{-1} \frac{1}{2}
  \frac{d}{ds}
    \norm{ S(s + \delta -s') Q B_{\mathrm{lin}}\big(V(s')\big) }_{L^2}^2
\end{array}
\end{equation}
for the values of $(s,s')$ that are relevant below.
Upon writing
\begin{equation}
\begin{array}{lcl}
\mathcal{I}_{I}
 & = &
    E e^{-\e t  } \int_0^{t}
   e^{  \nu s } \int_{s  - 1}^{s }
       \big[ 1 + \frac{1}{\sqrt{s + \delta - s'} } ]
          \norm{V(s')}_{H^1}^2  \mathbf{1}_{s' < \tau} \, ds' \, ds ,
\\[0.2cm]
\mathcal{I}_{II}
& = &
E e^{-\e t } \int_0^{t}
   e^{\nu s } \int_{s - 1 }^{s }
    \frac{d}{ds}
    \norm{  S(s + \delta -s') Q B_{\mathrm{lin}}\big(V(s')\big) }_{L^2}^2  \mathbf{1}_{s' < \tau} ds'  \, ds ,
\\[0.2cm]
\end{array}
\end{equation}
we obtain the estimate
\begin{equation}
\begin{array}{lcl}
E e^{-\e t} \mathcal{I}^{\mathrm{sh}}_{\nu,\delta;B;\mathrm{lin}}(t)
& \le &
  (1 + \rho^{-1}) M^2 K^2_{B;\mathrm{lin}} \mathcal{I}_{I}
  - \frac{1}{2} \rho^{-1} \mathcal{I}_{II} .
\end{array}
\end{equation}
Changing the integration order, we obtain
\begin{equation}
\begin{array}{lcl}
 \mathcal{I}_{I} & =&
E e^{ - \e t}  \int_0^{t}
 \Big[\int_{s'}^{\min\{ t,s'+1 \} }
                      e^{\nu s}
                      \big[
                         1 + \frac{1}{\sqrt{s+\d-s'}} \big]  \, ds
                   \Big]   \nrm{V(s')}_{H^1}^2  \mathbf{1}_{s' < \tau}      \, ds '
\\[0.2cm]
& \le &
3  e^{\nu}
E e^{ - \e t }  \int_0^{t} e^{ \nu s' }  \nrm{V(s')}_{H^1}^2 \mathbf{1}_{s' < \tau}   \, ds '
\\[0.2cm]
& \le &
3  e^{\nu}
E e^{ - \e t \wedge \tau }  \int_0^{t \wedge \tau} e^{ \nu s' }  \nrm{V(s')}_{H^1}^2
   \, ds '
\\[0.2cm]
& = &
3  e^{\nu}
 E N_{\e,\alpha;II}(t\wedge\tau).
\\[0.2cm]
\end{array}
\end{equation}
Integrating by parts, we arrive at the identity
\begin{equation}
\begin{array}{lcl}
 \mathcal{I}_{II} & =&
 E e^{ - \e t}  \int_0^{t}
 \Big[\int_{s'}^{\min\{ t,s'+1 \} }
                      e^{\nu s}
                       \frac{d}{ds}
    \norm{  S(s + \delta -s') Q B_{\mathrm{lin}}\big(V(s')\big) }_{L^2}^2   \, ds
                   \Big]    \mathbf{1}_{s' < \tau}    \, ds '
 \\[0.2cm]
 & = &
   \mathcal{I}_{II;A}
   + \mathcal{I}_{II;B}
   + \mathcal{I}_{II;C},
\end{array}
\end{equation}
in which we have introduced the expressions
\begin{equation}
\begin{array}{lcl}
\mathcal{I}_{II;A}
 & = &
   E e^{ - \e t}  \int_0^{t}
     e^{\nu \min\{ t , s' + 1 \} }
       \norm{  S(\min\{t , s'+1\} + \delta -s')
         Q B_{\mathrm{lin}}\big(V(s')\big) }_{L^2}^2  \mathbf{1}_{s' < \tau} \, ds' ,
 \\[0.2cm]
\mathcal{I}_{II;B}
 & = &
   - E e^{ - \e t}  \int_0^{t}
     e^{\nu s' }
       \norm{  S(\delta) Q B_{\mathrm{lin}}\big(V(s')\big) }_{L^2}^2
         \mathbf{1}_{s' < \tau} \, ds' ,
\\[0.2cm]
\mathcal{I}_{II;C}
 & = &
   - E e^{- \e t}  \int_0^{t}
 \Big[\int_{s'}^{\min\{ t,s'+1 \} }
                      \nu e^{\nu s}
    \norm{  S(s + \delta -s') Q B_{\mathrm{lin}}\big(V(s')\big) }_{L^2}^2   \, ds
                   \Big]   \mathbf{1}_{s' < \tau}     \, ds ' .
 \\[0.2cm]
\end{array}
\end{equation}
Inspecting these expressions, we readily obtain the bound
\begin{equation}
\begin{array}{lcl}
\abs{\mathcal{I}_{II} }
& \le &
  e^{\nu} (2 + \nu)
  M^2 K_{B;\mathrm{lin}}^2
    E e^{ - \e t}
    \int_0^{t}
       e^{\nu s'}
         \norm{V(s')}_{H^1}^2 \mathbf{1}_{s' < \tau} \, ds'
\\[0.2cm]
& \le &
  e^{\nu} (2 + \nu)
  M^2 K_{B;\mathrm{lin}}^2
    E N_{\e, \alpha;II}(t \wedge \tau).
\end{array}
\end{equation}
\end{proof}

\begin{lem}
\label{lem:nls:b:cn:i}
Fix $T > 0$ and assume that (HA), (HTw), (HS),
(H$\beta$),
(hSol) and (hFB) all hold.
Pick a constant $0 < \e <  2 \beta$.
Then for
any $0 \le \delta < 1$,
and any $0 \le t \le T$,
we have the bounds
\begin{equation}
\begin{array}{lcl}
E e^{-\e t} \mathcal{I}^{\mathrm{lt}}_{\e,\delta;B;\mathrm{cn}}(t)
 & \le & \frac{M^2}{(2 \beta - \e)\e} K_{B;\mathrm{cn}}^2 ,
\\[0.2cm]
E  e^{-\e t} \mathcal{I}^{\mathrm{sh}}_{\e,\delta;B;\mathrm{cn}}(t)
 & \le &
 K_{B;\mathrm{cn}}^2 \frac{M^2}{\e} (1 + \rho^{-1}) e^{\e}
   ( 4  + \e ).
\end{array}
\end{equation}
\end{lem}
\begin{proof}
Using the fact that
\begin{equation}
e^{-\e t} \int_0^{t}
  e^{ \e s } \,ds
\le \frac{1}{\e},
\end{equation}
these bounds can be obtained
from Lemmas \ref{lem:nls:b:lin:lt:i}
and \ref{lem:nls:b:lin:st:i}
by picking $\alpha = 0$
and making the substitutions
\begin{equation}
  K_{B;\mathrm{lin}} \mapsto K_{B;\mathrm{cn}},
  \qquad
  \qquad
  E N_{\e,0;II}(t\wedge\tau) \mapsto \frac{1}{\e}.
\end{equation}
\end{proof}

\begin{lem}
\label{lem:nls:st:sup:on:b:lin:i}
Fix $T > 0$ and assume that (HA), (HTw), (HS),
(H$\beta$),
(hSol) and (hFB) all hold.
Pick two constants $\e > 0$, $\alpha \ge 0$
for which $\e + \alpha < 2 \beta$
and write $\nu = \alpha + \e$.
Then we have the bounds
\begin{equation}
\begin{array}{lcl}
E \sup_{0 \le t \le \tau}
    e^{-\e t} \mathcal{I}^{\mathrm{lt}}_{\nu,0;B; \mathrm{lin}}(t)
& \le &
  e^{\e} (T + 1) \frac{M^2}{2 \beta - \nu} K_{B;\mathrm{lin}}^2
  E  \sup_{0 \le t \le \tau} N_{\e;\alpha;II}(t) ,
\\[0.2cm]
E \sup_{0 \le t \le \tau}
    e^{-\e t} \mathcal{I}^{\mathrm{sh}}_{\nu,0;B; \mathrm{lin}}(t)
& \le &
  e^{\e} (T + 1)
   K_{B;\mathrm{lin}}^2 M^2 (1 + \rho^{-1}) e^{\nu}( 4  + \nu )
  E  \sup_{0 \le t \le \tau} N_{\e;\alpha;II}(t) .
\end{array}
\end{equation}
\end{lem}
\begin{proof}
By splitting the integration interval we obtain
\begin{equation}
\begin{array}{lcl}
\sup_{0 \le t \le \tau}
    e^{-\e t} \mathcal{I}^{\mathrm{lt}}_{\nu,0;B; \mathrm{lin}}(t)
  & \le   &
  \sup_{0 \le t \le T}
    e^{-\e t} \mathcal{I}^{\mathrm{lt}}_{\nu, 0;B; \mathrm{lin}}(t)
\\[0.2cm]
& = &
   \sup_{0 \le t \le T}
    e^{-\e t}
      \int_0^t e^{ \nu s}
        \norm{\mathcal{E}^{\mathrm{lt}}_{B;\mathrm{lin}}(s)}_{H^1}^2 \, ds
\\[0.2cm]
& \le &
    e^{\e}   e^{-\e} \int_0^{1 }
     e^{ \nu s}
       \norm{\mathcal{E}^{\mathrm{lt}}_{B;\mathrm{lin}}(s)}_{H^1}^2 \, ds
\\[0.2cm]
& & \qquad
  + e^{\e}  e^{-2 \e}\int_1^{2 }
       e^{ \nu s}
       \norm{\mathcal{E}^{\mathrm{lt}}_{B;\mathrm{lin}}(s)}_{H^1}^2 \, ds
\\[0.2cm]
& & \qquad
  + e^{\e}  e^{- (\lfloor T \rfloor + 1) \e}
     \int_{\lfloor T \rfloor}^{\lfloor T \rfloor + 1}
        e^{ \nu s}
       \norm{\mathcal{E}^{\mathrm{lt}}_{B;\mathrm{lin}}(s)}_{H^1}^2 \, ds
\\[0.2cm]
& = & e^{\e}
\big[
     e^{-\e }\mathcal{I}^{\mathrm{lt}}_{\nu,0;B; \mathrm{lin}}(1)
  +  e^{-2\e }\mathcal{I}^{\mathrm{lt}}_{\nu,0;B; \mathrm{lin}}(2)
\\[0.2cm]
& & \qquad
  + \ldots
  +  e^{-\e (\lfloor T \rfloor + 1)} \mathcal{I}^{\mathrm{lt}}_{\nu,0;B; \mathrm{lin}}
       (\lfloor T \rfloor + 1)
\big] .
\\[0.2cm]
\end{array}
\end{equation}
Applying Lemma \ref{lem:nls:b:lin:lt:i},
we hence see
\begin{equation}
\begin{array}{lcl}
E \sup_{0 \le t \le \tau}
    e^{-\e t}\mathcal{I}^{\mathrm{lt}}_{\nu, 0;B; \mathrm{lin}}(t)
& \le &
e^{\e} \frac{M^2}{2 \beta - \nu} K_{B;\mathrm{lin}}^2
E
\big[
  N_{\e;\alpha;II}(1 \wedge \tau)
  + \ldots + N_{\e;\alpha;II}((\lfloor T \rfloor + 1) \wedge \tau)
\big]
\\[0.2cm]
& \le &
(T + 1) e^{\e} \frac{M^2}{2 \beta - \nu} K_{B;\mathrm{lin}}^2
\mathrm{sup}_{0 \le t \le \tau} N_{\e;\alpha;II}(t) .
\end{array}
\end{equation}
The same procedure works for the second estimate.
\end{proof}

\begin{lem}
\label{lem:nls:st:sup:on:b:cn:i}
Fix $T > 0$ and assume that (HA), (HTw), (HS),
(H$\beta$),
(hSol) and (hFB) all hold.
Pick a constant $\e > 0$, $\alpha \ge 0$
for which $\e  < 2 \beta$.
Then we have the bounds
\begin{equation}
\begin{array}{lcl}
E \sup_{0 \le t \le \tau}
    e^{-\e t} \mathcal{I}^{\mathrm{lt}}_{\nu,0;B; \mathrm{cn}}(t)
& \le &
  e^{\e} (T + 1)
\frac{M^2}{(2 \beta - \e)\e} K_{B;\mathrm{cn}}^2 ,
\\[0.2cm]
E \sup_{0 \le t \le \tau}
    e^{-\e t} \mathcal{I}^{\mathrm{sh}}_{\nu,0;B; \mathrm{cn}}(t)
& \le &
  e^{\e} (T + 1)
K_{B;\mathrm{cn}}^2 \frac{M^2}{\e} (1 + \rho^{-1}) e^{\e}
   ( 4  + \e ) .
\end{array}
\end{equation}
\end{lem}
\begin{proof}
Following the procedure in the proof of Lemma \ref{lem:nls:st:sup:on:b:lin:i},
these bounds can be obtained
from the estimates in Lemma \ref{lem:nls:b:cn:i} .
\end{proof}

\begin{proof}[Proof of Proposition \ref{prp:nls:general}]
Pick $T >0$ and $0 < \eta < \eta_0$
and write $\tau = \tau_{\e,\alpha}(T, \eta)$.
Since the identities \sref{eq:nls:projs:are:zero}
with $v = V(t\wedge\tau)$ hold for all $0 \le t \le T$,
we may compute
\begin{equation}
\begin{array}{lcl}
E \sup_{0 \le t \le \tau} N_{\e,0;I}(t)
& \le &
5 E
\sup_{0 \le t \le \tau}
\Big[
\norm{\mathcal{E}_0(t)}_{L^2}^2
+ \sigma^4 \norm{\mathcal{E}_{F;\mathrm{lin}}(t)}_{L^2}^2
+ \norm{\mathcal{E}_{F;\mathrm{nl}}(t)}_{L^2}^2
\\[0.2cm]
& & \qquad
+ \sigma^2 \norm{\mathcal{E}_{B;\mathrm{lin}}(t)}_{L^2}^2
+ \sigma^2 \norm{\mathcal{E}_{B;\mathrm{cn}}(t)}_{L^2}^2
\Big]
\end{array}
\end{equation}
by applying Young's inequality.
The inequalities in
Lemma's \ref{lem:nls:e:zero}-\ref{lem:nls:st:sup:on:b:cn:i}
now imply that
\begin{equation}
\begin{array}{lcl}
E \sup_{0 \le t \le \tau} N_{\e,0;I}(t)
& \le &
C_1 \big[ \norm{V(0)}_{H^1}^2
  + \sigma^2 T
  + \Big( \eta
  + \sigma^2 T + \sigma^4 \Big)
    E \sup_{0 \le t \le \tau}  N_{\e,0;II}(t)
\big].
\end{array}
\end{equation}
In addition,
we note that
\begin{equation}
\begin{array}{lcl}
E \sup_{0 \le t \le \tau} N_{\e,0;II}(t)
& \le &
9 E \sup_{0 \le t \le \tau}
e^{-\e t}
\Big[
  \mathcal{I}_{\nu,0;0}(t)
+ \sigma^4 \mathcal{I}^{\mathrm{lt}}_{\nu,0;F;\mathrm{lin}}(t)
+ \sigma^4 \mathcal{I}^{\mathrm{sh}}_{\nu,0;F;\mathrm{lin}}(t)
\\[0.2cm]
& & \qquad
+ \mathcal{I}^{\mathrm{lt}}_{\nu,0;F;\mathrm{nl}}(t)
+ \mathcal{I}^{\mathrm{sh}}_{\nu,0;F;\mathrm{nl}}(t)
\\[0.2cm]
& & \qquad
+ \sigma^2 \mathcal{I}^{\mathrm{lt}}_{\nu,0;B;\mathrm{lin}}(t)
+ \sigma^2 \mathcal{I}^{\mathrm{sh}}_{\nu,0;B;\mathrm{lin}}(t)
\\[0.2cm]
& & \qquad
+\sigma^2 \mathcal{I}^{\mathrm{lt}}_{\nu,0;B;\mathrm{cn}}(t)
+ \sigma^2 \mathcal{I}^{\mathrm{sh}}_{\nu,0;B;\mathrm{cn}}(t)
\Big] .
\end{array}
\end{equation}
The inequalities in
Lemma's \ref{lem:nls:e:zero}-\ref{lem:nls:st:sup:on:b:cn:i}
now imply that
\begin{equation}
\begin{array}{lcl}
E \sup_{0\le t \le \tau} N_{\e,0;II}(t)
& \le &
C_2 \big[ 
    \norm{V(0)}_{H^1}^2
    + \sigma^2T
  + \Big( \eta
  + \sigma^2 T + \sigma^4 \Big) \sup_{0\le t \le \tau}
    N_{\e,0;II}(t)
\big].
\end{array}
\end{equation}
In particular, we see that
\begin{equation}
E \sup_{0 \le t \le \tau} N_{\e,0}(t)
\le C_3 \big[ 
  \norm{V(0)}_{H^1}^2
  + \sigma^2 T
  + (\eta + \sigma^2 T + \sigma^4)
      E \sup_{0 \le t \le \tau} N_{\e,0}(t)
\big] .
\end{equation}
The desired bound hence follows
by appropriately restricting the size of
$\eta + \sigma^2 T + \sigma^4$.
\end{proof}

\begin{proof}[Proof of Proposition \ref{prp:nls:exponential}]
Ignoring the contributions
arising from $B_{\mathrm{cn}}$,
we can follow the proof of
Proposition \ref{prp:nls:general}
to obtain the bound
\begin{equation}
E \sup_{0 \le t \le \tau} N_{\e,\alpha}( t)
\le C_4 \big[   
  \norm{V(0)}_{H^1}^2
  + (\eta + \sigma^2 T + \sigma^4)
      E \sup_{0 \le t \le \tau} N_{\e,\alpha}(t)
  \big].
\end{equation}
The desired estimate hence follows
by appropriately restricting the size of
$\eta + \sigma^2 T + \sigma^4$.
\end{proof}

\bibliographystyle{klunumHJ}
\bibliography{ref}

\begin{thebibliography}{000}

\bibitem{alili2005representations}
L. Alili, P. Patie and J.~L. Pedersen (2005), Representations of the first
  hitting time density of an Ornstein-Uhlenbeck process.
\newblock {\em Stochastic Models} {\bf 21}(4), 967--980.

\bibitem{HJHNLS}
M. Beck, H.~J. Hupkes, B. Sandstede and K. Zumbrun (2010), Nonlinear
  {S}tability of {S}emidiscrete {S}hocks for {T}wo-{S}ided {S}chemes.
\newblock {\em SIAM J. Math. Anal.} {\bf 42}, 857--903.

\bibitem{beck2010nonlinear}
M. Beck, B. Sandstede and K. Zumbrun (2010), Nonlinear stability of
  time-periodic viscous shocks.
\newblock {\em Archive for rational mechanics and analysis} {\bf 196}(3),
  1011--1076.

\bibitem{BHM}
H. Berestycki, F. Hamel and H. Matano (2009), Bistable traveling waves around
  an obstacle.
\newblock {\em Comm. Pure Appl. Math.} {\bf 62}(6), 729--788.

\bibitem{brassesco1995brownian}
S. Brassesco, A. De~Masi and E. Presutti (1995), Brownian fluctuations of the
  interface in the D=1 Ginzburg-Landau equation with noise.
\newblock {\em Ann. Inst. H. Poincar{\'e} Probab. Statist} {\bf 31}(1),
  81--118.

\bibitem{bressloff2015nonlinear}
P.~C. Bressloff and Z.~P. Kilpatrick (2015), Nonlinear Langevin equations for
  wandering patterns in stochastic neural fields.
\newblock {\em SIAM Journal on Applied Dynamical Systems} {\bf 14}(1),
  305--334.

\bibitem{Bressloff}
P.~C. Bressloff and M.~A. Webber (2012), Front propagation in stochastic neural
  fields.
\newblock {\em SIAM Journal on Applied Dynamical Systems} {\bf 11}(2),
  708--740.

\bibitem{Cartwright2019}
M. Cartwright and G.~A. Gottwald (2019), A collective coordinate framework to
  study the dynamics of travelling waves in stochastic partial differential
  equations.
\newblock {\em Physica D: Nonlinear Phenomena}.

\bibitem{chen2015traveling}
C.-N. Chen and Y. Choi (2015), Traveling pulse solutions to FitzHugh--Nagumo
  equations.
\newblock {\em Calculus of Variations and Partial Differential Equations} {\bf
  54}(1), 1--45.

\bibitem{Chow}
P.-L. Chow (2014), {\em Stochastic partial differential equations}.
\newblock CRC Press.

\bibitem{cornwell2017opening}
P. Cornwell (2017), Opening the Maslov Box for Traveling Waves in Skew-Gradient
  Systems.
\newblock {\em arXiv preprint arXiv:1709.01908}.

\bibitem{cornwell2017existence}
P. Cornwell and C.~K. Jones (2018), On the Existence and Stability of Fast
  Traveling Waves in a Doubly Diffusive FitzHugh--Nagumo System.
\newblock {\em SIAM Journal on Applied Dynamical Systems} {\bf 17}(1),
  754--787.

\bibitem{DaPratomild}
G. Da~Prato, A. Jentzen and M. R{\"o}ckner (2010), A mild {I}t{\^o} formula for
  {SPDE}s.
\newblock {\em arXiv preprint arXiv:1009.3526}.

\bibitem{NunnoAdvMathFinance2011}
G. di~Nunno and B.~O. (editors) (2011), {\em Advanced Mathematical Methods for
  Finance}.
\newblock Springer.

\bibitem{Evans}
L. Evans (1998), {\em Partial differential equations}.
\newblock American Mathematical Society, Providence, R.I.

\bibitem{Fife1977}
P.~C. Fife and J.~B. McLeod (1977), The approach of solutions of nonlinear
  diffusion equations to travelling front solutions.
\newblock {\em Arch. Ration. Mech. Anal.} {\bf 65}(4), 335--361.

\bibitem{Climate}
C.~L.~E. Franzke, T.~J. O'Kane, J. Berner, P.~D. Williams and V. Lucarini
  (2015), Stochastic climate theory and modeling.
\newblock {\em Wiley Interdisciplinary Reviews: Climate Change} {\bf 6}(1),
  63--78.

\bibitem{funaki1995scaling}
T. Funaki (1995), The scaling limit for a stochastic PDE and the separation of
  phases.
\newblock {\em Probability Theory and Related Fields} {\bf 102}(2), 221--288.

\bibitem{gardner1982existence}
R.~A. Gardner (1982), Existence and stability of travelling wave solutions of
  competition models: a degree theoretic approach.
\newblock {\em Journal of Differential equations} {\bf 44}(3), 343--364.

\bibitem{Gawarecki}
L. Gawarecki and V. Mandrekar (2010), {\em Stochastic differential equations in
  infinite dimensions: with applications to stochastic partial differential
  equations}.
\newblock Springer Science \& Business Media.

\bibitem{gowda2015early}
K. Gowda and C. Kuehn (2015), Early-warning signs for pattern-formation in
  stochastic partial differential equations.
\newblock {\em Communications in Nonlinear Science and Numerical Simulation}
  {\bf 22}(1), 55--69.

\bibitem{hale1999stability}
J. Hale, L.~A. Peletier and W.~C. Troy (1999), Stability and instability in the
  Gray-Scott model: the case of equal diffusivities.
\newblock {\em Applied mathematics letters} {\bf 12}(4), 59--65.

\bibitem{Hamster2018Uneq}
C. Hamster and H. Hupkes (2018), Stability of Travelling Waves for Systems of
  Reaction-Diffusion Equations with Multiplicative Noise.
\newblock {\em preprint}.

\bibitem{HJHSTB2D}
A. Hoffman, H. Hupkes and E. Van~Vleck (2015), Multi-dimensional Stability of
  Waves Travelling through Rectangular Lattices in Rational Directions.
\newblock {\em Transactions of the American Mathematical Society} {\bf
  367}(12), 8757--8808.

\bibitem{HJHOBST2D}
A. Hoffman, H. Hupkes and E. Van~Vleck (2017), {\em Entire {S}olutions for
  {B}istable {L}attice {D}ifferential {E}quations with {O}bstacles}.
\newblock American Mathematical Society.

\bibitem{hytonen2018analysis}
T. Hyt{\"o}nen, J. Van~Neerven, M. Veraar and L. Weis (2018), {\em Analysis in
  Banach Spaces: Volume II: Probabilistic Methods and Operator Theory},
  Vol.~67.
\newblock Springer.

\bibitem{iannelli2015introduction}
M. Iannelli and A. Pugliese (2015), {\em An Introduction to Mathematical
  Population Dynamics: Along the Trail of Volterra and Lotka}, Vol.~79.
\newblock Springer.

\bibitem{Inglis}
J. Inglis and J. MacLaurin (2016), A general framework for stochastic traveling
  waves and patterns, with application to neural field equations.
\newblock {\em SIAM Journal on Applied Dynamical Systems} {\bf 15}(1),
  195--234.

\bibitem{jacod2006calcul}
J. Jacod (2006), {\em Calcul stochastique et problemes de martingales}, Vol.
  714.
\newblock Springer.

\bibitem{jeanblanc2009mathematical}
M. Jeanblanc, M. Yor and M. Chesney (2009), {\em Mathematical methods for
  financial markets}.
\newblock Springer Science \& Business Media.

\bibitem{KAP1997}
T. Kapitula (1997), Multidimensional {S}tability of {P}lanar {T}ravelling
  {W}aves.
\newblock {\em Trans. Amer. Math. Soc.} {\bf 349}, 257--269.

\bibitem{kruger2017multiscale}
J. Kr{\"u}ger and W. Stannat (2017), A multiscale-analysis of stochastic
  bistable reaction--diffusion equations.
\newblock {\em Nonlinear Analysis} {\bf 162}, 197--223.

\bibitem{Kuske2017}
R. Kuske, C. Lee and V. Rottsch{\"a}fer (2017), Patterns and coherence
  resonance in the stochastic Swift-Hohenberg equation with Pyragas control:
  The Turing bifurcation case.
\newblock {\em Physica D: Nonlinear Phenomena}.

\bibitem{Lang}
E. Lang (2016), A multiscale analysis of traveling waves in stochastic neural
  fields.
\newblock {\em SIAM Journal on Applied Dynamical Systems} {\bf 15}(3),
  1581--1614.

\bibitem{LangStannat2016l2}
E. Lang and W. Stannat (2016), L2-stability of traveling wave solutions to
  nonlocal evolution equations.
\newblock {\em Journal of Differential Equations} {\bf 261}(8), 4275--4297.

\bibitem{lee1994experimental}
K.-J. Lee, W.~D. McCormick, J.~E. Pearson and H.~L. Swinney (1994),
  Experimental observation of self-replicating spots in a reaction--diffusion
  system.
\newblock {\em Nature} {\bf 369}(6477), 215.

\bibitem{LiuRockner}
W. Liu and M. R{\"o}ckner (2010), {SPDE} in {H}ilbert space with locally
  monotone coefficients.
\newblock {\em Journal of Functional Analysis} {\bf 259}(11), 2902--2922.

\bibitem{Lord2012}
G. Lord and V. Th{\"u}mmler (2012), Computing stochastic traveling waves.
\newblock {\em SIAM Journal on Scientific Computing} {\bf 34}(1), B24--B43.

\bibitem{lorenzi2004analytic}
L. Lorenzi, A. Lunardi, G. Metafune and D. Pallara (2004), Analytic semigroups
  and reaction-diffusion problems.
\newblock In: {\em Internet Seminar}, Vol. 2005.
\newblock p. 127.

\bibitem{MasciaZumbrun02}
C. Mascia and K. Zumbrun (2002), Pointwise {G}reen's function bounds and
  stability of relaxation shocks.
\newblock {\em Indiana Univ. Math. J.} {\bf 51}(4), 773--904.

\bibitem{ORG1998}
M. Or-Guil, M. Bode, C.~P. Schenk and H.~G. Purwins (1998), Spot {B}ifurcations
  in {T}hree-{C}omponent {R}eaction-{D}iffusion {S}ystems: {T}he {O}nset of
  {P}ropagation.
\newblock {\em Physical Review E} {\bf 57}(6), 6432.

\bibitem{DaPratoZab}
G. Prato and J. Zabczyk (1992), {\em Stochastic equations in infinite
  dimensions}.
\newblock Cambridge University Press, Cambridge New York.

\bibitem{Concise}
C. Pr{\'e}v{\^o}t and M. R{\"o}ckner (2007), {\em A concise course on
  stochastic partial differential equations}, Vol. 1905.
\newblock Springer.

\bibitem{revuz2013continuous}
D. Revuz and M. Yor (2013), {\em Continuous martingales and Brownian motion},
  Vol. 293.
\newblock Springer Science \& Business Media.

\bibitem{ricciardi1988first}
L.~M. Ricciardi and S. Sato (1988), First-passage-time density and moments of
  the Ornstein-Uhlenbeck process.
\newblock {\em Journal of Applied Probability} {\bf 25}(1), 43--57.

\bibitem{rinzel1973traveling}
J. Rinzel and J.~B. Keller (1973), Traveling wave solutions of a nerve
  conduction equation.
\newblock {\em Biophysical journal} {\bf 13}(12), 1313--1337.

\bibitem{Sattinger}
D.~H. Sattinger (1976), On the stability of waves of nonlinear parabolic
  systems.
\newblock {\em Advances in Mathematics} {\bf 22}(3), 312--355.

\bibitem{SCHENK1997}
C.~P. Schenk, M. Or-Guil, M. Bode and H.~G. Purwins (1997), Interacting
  {P}ulses in {T}hree-component {R}eaction-{D}iffusion {S}ystems on
  {T}wo-{D}imensional {D}omains.
\newblock {\em Physical Review Letters} {\bf 78}(19), 3781.

\bibitem{Shardlow}
T. Shardlow (2005), Numerical simulation of stochastic PDEs for excitable
  media.
\newblock {\em Journal of computational and applied mathematics} {\bf 175}(2),
  429--446.

\bibitem{Stannat}
W. Stannat (2013), Stability of travelling waves in stochastic {N}agumo
  equations.
\newblock {\em arXiv preprint arXiv:1301.6378}.

\bibitem{stannat2014stability}
W. Stannat (2014), Stability of travelling waves in stochastic bistable
  reaction-diffusion equations.
\newblock {\em arXiv preprint arXiv:1404.3853}.

\bibitem{tonnelier2003piecewise}
A. Tonnelier and W. Gerstner (2003), Piecewise linear differential equations
  and integrate-and-fire neurons: insights from two-dimensional membrane
  models.
\newblock {\em Physical Review E} {\bf 67}(2), 021908.

\bibitem{van2010front}
P. van Heijster, A. Doelman, T.~J. Kaper and K. Promislow (2010), Front
  interactions in a three-component system.
\newblock {\em SIAM Journal on Applied Dynamical Systems} {\bf 9}(2), 292--332.

\bibitem{veraar2011note}
M. Veraar and L. Weis (2011), A note on maximal estimates for stochastic
  convolutions.
\newblock {\em Czechoslovak mathematical journal} {\bf 61}(3), 743.

\bibitem{vinals1991numerical}
J. Vi{\~n}als, E. Hern{\'a}ndez-Garc{\'\i}a, M. San~Miguel and R. Toral (1991),
  Numerical study of the dynamical aspects of pattern selection in the
  stochastic Swift-Hohenberg equation in one dimension.
\newblock {\em Physical Review A} {\bf 44}(2), 1123.

\bibitem{volpert1994traveling}
A.~I. Volpert, V.~A. Volpert and V.~A. Volpert (1994), {\em Traveling wave
  solutions of parabolic systems}, Vol. 140.
\newblock American Mathematical Soc.

\bibitem{weis2006h}
L. Weis (2006), The $H^\infty$ holomorphic functional calculus for sectorial
  operators - a survey.
\newblock In: {\em Partial differential equations and functional analysis}.
\newblock Springer, pp. 263--294.

\bibitem{zemskov2010wave}
E. Zemskov and I. Epstein (2010), Wave propagation in a FitzHugh-Nagumo-type
  model with modified excitability.
\newblock {\em Physical Review E} {\bf 82}(2), 026207.

\bibitem{Zhang}
J. Zhang, A. Holden, O. Monfredi, M. Boyett and H. Zhang (2009), Stochastic
  vagal modulation of cardiac pacemaking may lead to erroneous identification
  of cardiac â€œchaosâ€.
\newblock {\em Chaos: An Interdisciplinary Journal of Nonlinear Science} {\bf
  19}(2), 028509.

\bibitem{Zumbrun2009}
K. Zumbrun (2011), Instantaneous {S}hock {L}ocation and {O}ne-{D}imensional
  {N}onlinear {S}tability of {V}iscous {S}hock {W}aves.
\newblock {\em Quarterly of applied mathematics} {\bf 69}(1), 177--202.

\end{thebibliography}
\end{document}